\documentclass[12pt]{amsart}
\usepackage{graphicx} % Required for inserting images
\usepackage{booktabs}
\usepackage{float}
\usepackage{soul}
\usepackage[citecolor=Navy]{hyperref}
\usepackage[ruled,vlined]{algorithm2e}
\usepackage{amsfonts,amssymb,amsmath,amsthm}
\usepackage{algorithmic}
\usepackage{amssymb,amscd,amsxtra,calc}
\usepackage{cmmib57}
\usepackage{graphicx}
\usepackage{mathrsfs,enumerate,xcolor}

\usepackage[all]{xy}
\usepackage{multirow} % Required for multirows
\usepackage{psfrag,xmpmulti,amscd,color,pstricks, import}
\usepackage{tikz-cd} 
\newtheorem{thm}{Theorem}[section]

\newtheorem{lem}[thm]{Lemma}

\theoremstyle{definition}
\newtheorem{defn}[thm]{Definition}
\newtheorem{rem}[thm]{Remark}
\newtheorem*{rem*}{Remark}

\numberwithin{equation}{section}

\definecolor{OrangeRed}{cmyk}{0,0.6,1,0}            % half magenta only, full yellow
\definecolor{DarkBlue}{cmyk}{1,1,0,0.20}
\definecolor{DarkGreen}{cmyk}{1,0,0.6,0.2}
\definecolor{myblue}{rgb}{0.66,0.78,1.00}
\definecolor{Violet}{cmyk}{0.79,0.88,0,0}
\definecolor{Lavender}{cmyk}{0,0.48,0,0}

\renewcommand{\epsilon}{\varepsilon}
\renewcommand{\phi}{\varphi}

\title[BNQN for finding roots of meromorphic functions]{Backtracking New Q-Newton's method for finding roots of meromorphic functions in 1 complex variable: Global convergence, and local stable/unstable curves}

\date{}

\author[Forn\ae ss]{John Erik Forn\ae ss}
\address{Department of Mathematics, NTNU, Norway }
\email{fornaess@gmail.com}

\author[Hu]{Mi Hu}
\address{Department of Mathematics, University of Oslo, Norway}
\email{humihqu@gmail.com}

\author[Truong]{Tuyen Trung Truong}
\address{Department of Mathematics, University of Oslo, Norway}
\email{tuyentt@math.uio.no}

\begin{document}
	    \maketitle

     	    \begin{abstract}
In this paper, we research more in depth properties of Backtracking New Q-Newton's method (recently designed by the third author),  when used to find roots of meromorphic functions. 

If $f=P/Q$, where $P$ and $Q$ are polynomials  in 1 complex variable z with $\deg (P)>\deg (Q)$, we show the existence of an exceptional set $\mathcal{E}\subset\mathbf{C}$, which is contained in a countable union of real analytic curves in $\mathbf{R}^2=\mathbf{C}$, so that the following statements A and B hold. Here, $\{z_n\}$ is the sequence constructed by BNQN with an initial point $z_0$ which is not a pole of $f$.

A) If $z_0\in\mathbf{C}\backslash\mathcal{E}$, then $\{z_n\}$ converges to a root of $f$.  

B) If $z_0\in \mathcal{E}$, then $\{z_n\}$ converges to a critical point - but not a root - of $f$. 

Experiments seem to indicate that in general, even when $f$ is a polynomial, the set $\mathcal{E}$ is not contained in a finite union of real analytic curves. We provide further results relevant to whether locally $\mathcal{E}$ is contained in a finite number of real analytic curves. A similar result holds for general meromorphic functions. Moreover, unlike previous work, here we do not require that the parameters of BNQN are random, or that the meromorphic function $f$ is generic.  

Based on the theoretical results, we explain (both rigorously and heuristically) of what observed in experiments with BNQN, in previous works by the authors. In particular, the dynamics of BNQN (an iterative method) seems to have some striking similarities to Newton's method (a continuous method) and the classical Poincar\'e-Bendixon theorem for differentiable real dynamical systems on the complex plane. This is the more interesting given that discrete versions of Newton's method (e.g. Relaxed Newton's method) does not behave this way.

	    \end{abstract}
\section{Introduction}

Finding roots of polynomials and more generally meromorphic functions in 1 complex variable is an active topic of interest in many sub-fields in mathematics, including: complex analysis, dynamical systems and number theory. For example, the Riemann hypothesis concerns the roots of the Riemann zeta function. A class of meromorphic functions, usually called special functions (e.g. Gamma function), is helpful also in other fields e.g. physics. 

A favourite method used in this topic is Newton's method, where to find roots of a polynomial or meromorphic function, one chooses a random initial point $z_0$ and construct a sequence: 
\begin{eqnarray*}
z_{n+1}=z_n-\frac{f(z_n)}{f'(z_n)}. 
\end{eqnarray*}

If $z_0$ is close enough to a root $z^*$ of $f$, then $\{z_n\}$ will converge to $z^*$. Moreover, if $z^*$ is a simple root of $f$, then the local rate of convergence is quadratic, i.e.  there is a positive constant $C>0$ such that $|z_{n+1}-z^*|\leq C|z_n-z^*|^2$ for all $n$. The classical result\cite{RefSE2} by Schroder on Newton's method applied to a polynomial of degree 2 was a starting point for the field of complex dynamics. For some good references on complex dynamics in general, the reader can consult \cite{RefAlex}, \cite{RefB}, \cite{RefBer}, \cite{RefCG}, \cite{RefD} and \cite{RefM}. 

On the other hand, Newton's method suffers from having no global convergence guarantee. Indeed, \cite{RefMc} shows that there is no algebraic algorithm which has global guarantee for finding roots of a generic polynomial of degree $4$ when starting from a random initial point $z_0$. (Here random means we choose these points outside a set of 0 measure.)

A variant of Newton's method, the so-called Relaxed Newton's method, which is defined as follows: 
\begin{eqnarray*}
z_{n+1}=z_n-\gamma \frac{f(z_n)}{f'(z_n)},
\end{eqnarray*}
where $\gamma\in\mathbf{C}$ is a small enough constant, has better global convergence guarantee for polynomials. In fact, \cite{{RefHK}} and \cite{RefMe} show that if certain conditions are satisfied for the polynomial $f$ (or more generally a rational function), then global convergence to roots is achieved. The idea is to relate this to Newton's flow method. (For a good survey on Relaxed Newton's method, Newton's flow and some other methods, the reader can consult \cite{RefBer}.) On the other hand, in practice  can have less than optimal rate of convergence. Indeed, there seems to be a striking similarity between the behaviour of BNQN (an iterative method) and that of Newton's method (a continuous method) and the classical Poincar\'e-Bendixon theorem \cite{RefPoincare}\cite{RefBendixon}, which will be discussed in more detail in Section 4.4. 

Yet, another variant is Random Relaxed Newton's method, with the update rule:  
\begin{eqnarray*}
z_{n+1}=z_n-\gamma _n \frac{f(z_n)}{f'(z_n)},
\end{eqnarray*}
where $\gamma _n\in \mathbf{C}$ is randomly chosen in an appropriate probability distribution. \cite{RefS} shows that
for each polynomial $f$ and a random initial point $z_0$, if the sequence $\{\gamma _n\}$ is appropriately randomly chosen (meaning outside a certain set of zero measure in an appropriate probablility space of all sequences $\{\gamma _n\}$) then there is global convergence guarantee to roots of $f$. 

However, there is no general enough global convergence for finding roots of meromorphic functions, concerning the above mentioned methods. To this end, we now discuss Backtracking New Q-Newton's method (BNQN), the main subject of study in this paper. 

The third author designed in \cite{RefT} a new variant of Newton's method, named Backtracking New Q-Newton's method (BNQN), which aims to minimise a function $F:\mathbf{R}^m\rightarrow \mathbf{R}$. For the sake of simplicity, we defer the definition of this method to Section 2. BNQN is shown to have the same rate of convergence as Newton's method near non-degenerate local minima, while having stronger global convergence guarantees (applicable to functions which have at most countably many critical points or satisfy certain Lojasiewicz gradient inequality, the reader is referred to \cite{RefT} for more detail). To find roots of meromorphic functions in 1 complex variable $z=x+iy$, one simply applies BNQN to the function $F:\mathbf{R}^2\rightarrow [0,+\infty]$ defined by $F(x,y)=|f(x+iy)|^2/2$. See Section 2 for a brief review of the method, in particular Theorem \ref{TheoremConvergence} on global convergence guarantee of BNQN when applied to a generic meromorphic function $f$ (i.e. both $f$ and $f'$ have only simple roots), and where the parameters of the method are randomly chosen. 

Recently, there have been (besides those in \cite{RefTT}, \cite{RefT}) further theoretical and experimental work concerning BNQN. In particular, the paper \cite{RefFHTW2} presents experimental results which show connections between the basin of attraction of BNQN and the Voronoi's diagrams (see \cite{RefV1}, \cite{RefV2}) of the roots, as well as to that of Newton's flow, and Newton's method for the function $g=f/f'$. The paper \cite{RefFHTW} establishes Schr\"oder theorem \cite{RefSE2} also for BNQN, and the paper \cite{RefThuanTuyen} proves a new equivalence of the Riemann hypothesis in terms of the dynamics of BNQN, as well as - based on experiments in using BNQN for the Riemann xi function - proposes a new approach towards the Riemann hypothesis. Takayuki Watanabe (in a conference talk in August 2024 at Osaka Metropolitan University) showed, via experiments, that in complex dimension 2, BNQN seems also behave better than Newton's method and Random Relaxed Newton's method. More precisely, when using BNQN method to find roots of the 2-dim  polynomial map $h:\mathbf{C}^2\rightarrow \mathbf{C}^2$  given by $h(z_1,z_2)=(z_1+z_1z_2,z_2^2-z_1^2)$ which has a multiple root at $(0,0)$, it is found experimentally that it seems $(0,0)$ is an attractor for the dynamics of BNQN, while not an attractor for the dynamics of both Newton's method and Random Relaxed Newton's method. It is known \cite{RefY} \cite{RefDDS} that for this polynomial map, also $(0,0)$ is not an attractor for the dynamics of Newton's method. The third author of this paper (in private correspondence, to appear somewhere else) has theoretically explained what observed in Takayuki Watanabe's experiments.   

The current paper researches more extensively on the global convergence of BNQN method for finding roots of meromorphic functions in 1 complex variable. In particular, we are able to resolve the following issues left open in \cite{RefTT}: 

Issue 1: Prove global convergence for a general meromorphic function $f$, which may have non-simple roots or critical points. 

Issue 2: Prove global convergence for general choice of parameters (and not just for randomly chosen parameters, again this means that the parameters must be chosen outside certain sets of measure 0) in BNQN. (In this regards, note that in Random Relaxed Newton's method \cite{RefS}, the sequence $\{\gamma _n\}$ is required to be random for global convergence for finding roots of polynomial functions.) Even so, note that the use of Armijo's Backtracking line search \cite{RefAr} in BNQN still adds random features into the algorithm. 

Issue 3: Describe in finer detail the exceptional set $\mathcal{E}$, beyond that (at least when $f$ is a generic meromorphic function) it has Lebesgue measure 0. 

\begin{defn}
 Here we define some notions to be used in the results. Fix $f$ a non-constant meromorphic function.
 
 1) We define: $\mathcal{P}(f)=\{z\in \mathbf{C}: |f(z)|=\infty\}$ the set of poles of $f$, $\mathcal{Z}(f)=\{z\in \mathbf{C}: f(z)=0\}$ the set of roots of $f$, and $\mathcal{C}(f)=\{z\in \mathbf{C}: f'(z)=0, f(z)\not=0\}$ the set of critical points but not roots of $f$.  

2) A set $A\subset\mathbf{C}$ is a real analytic curve if there is an open subset $U\subset\mathbf{C}$ so that $A$ is defined by a real analytic equation in $U$, and has dimension $1$. 

3) A set $\mathcal{E}\subset \mathbf{C}\backslash \mathcal{P}(f)$ is locally contained in a finite union of real analytic curves, if for each compact set $K\subset \mathbf{R}^2\backslash \mathcal{P}(f)$, there is  a finite open cover $\{\Omega _k\}$ of an open neighbourhood of $K$ and closed real analytic curves $R_k\subset \Omega _k$ (maybe empty) such that $\mathcal{E} \subset \bigcup _kR_k$.

4) A sublevel of $f$ is a set of the form  $\{z\in \mathbf{C}: |f(z)|\leq c\}$, for some $\infty >c\geq 0$. 

$f$ has compact sublevels if all sublevels are compact. Equivalently, $f$ has a pole at $\infty$, hence a rational function of the form $f(z)=P(z)/Q(z)$, where $P$ and $Q$ are polynomials with $\deg (P)>\deg (Q)$. In particular, to find roots of a polynomial $P(z)$, one can find roots of the rational function $f(z)=P(z)/P'(z)$, which has the advantage that all the roots are simple, where the local rate of convergence near roots for both Newton's method and BNQN are quadratic. 
\end{defn}

The first main result of this paper is as follows. 

\begin{thm} Let $f(z)$ be a non-constant meromorphic function. There is an exceptional set $\mathcal{E}\subset \mathbf{C}\backslash \mathcal{P}(f)$ - which is locally contained in a countable union of real analytic curves in $\mathbf{R}^2=\mathbf{C}$ - so that the following assertions A and B hold. Here $\{z_n\}$ is a sequence constructed by BNQN from an initial point $z_0\in \mathbf{C}\backslash \mathcal{P}(f)$. 

A) If $z_0\notin \mathcal{E}$, then either A1) $\{z_n\}$ converges to a point in $\mathcal{Z}(f)$, or A2) $\{z_n\}$ converges to $\infty$. 

B) On the other hand, if $z_0\in \mathcal{E}$ then $\{z_n\}$ converges to a point in $\mathcal{C}(f)$. 

If moreover f has compact sublevels, then possibility A2 does not happen. 
\label{Theorem1}\end{thm}

Experiments, for simple polynomials of degree $4$, seem to indicate that in general $\mathcal{E}$ cannot be contained in a finite union of real analytic curves, see Figures \ref{fig:FF3} and \ref{fig:FF4}. On the other hand, pictures seem to indicate that $\mathcal{E}$ should be locally contained in a finite union of real analytic curves. We note that the proof of Theorem \ref{Theorem1} constructs explicitly certain subsets of $\mathcal{E}$ which are contained in a finite number of real analytic curves and whose union is $\mathcal{E}$. Moreover, we have the following characterisation.   

\begin{thm} Assume that $f$ has compact sublevels, and $\mathcal{E}$ is the exceptional set given in Theorem \ref{Theorem1}. Then the following two statements are equivalent. 

1) Away from the poles of $f$, $\mathcal{E}$ is locally contained in a finite union of real analytic curves. More precisely, if $K\subset \mathbf{C}\backslash \mathcal{P}(f)$ is a compact set, then $\mathcal{E}\cap K$ is contained in a finite union of real analytic curves. 

2) Each point $z^*\in \mathcal{C}(f)$ has an open neighbourhood $U$ such that $\mathcal{E}\cap U$ is contained in a finite union of real analytic curves.

\label{Theorem1Bis}\end{thm}

The proof of Theorem \ref{Theorem1Bis} provides a uniform upper bound on the least number of iterations needed, for any point $z_0$ in a sublevel set $B_R=\{z: |f(z)|\leq R\}$ to land in an $\epsilon$ neighbourhood $D_{\epsilon}$ of $\mathcal{Z}(f)\cup \mathcal{C}(f)$, based on certain information of $f$ on $B_R\backslash D{\epsilon}$. Combining the above theorems and Theorem \ref{Theorem2} below, we provide in Section 4 some explanations (both rigorously and heuristically), in addition to \cite{RefFHTW}, several remarkable phenomena observed from experiments involving BNQN, in particular \cite{RefFHTW2}: the boundaries of basins of attraction for BNQN are more smooth than that of Newton's method and Random Relaxed Newton's method, basins of attraction for BNQN seem similar to that of Newton's flow and Voronoi's diagrams of the roots, and the apparent existence of ''channels to infinity" for every root of a polynomial.   

\begin{figure}
    \centering
    
    \includegraphics[width=10cm]{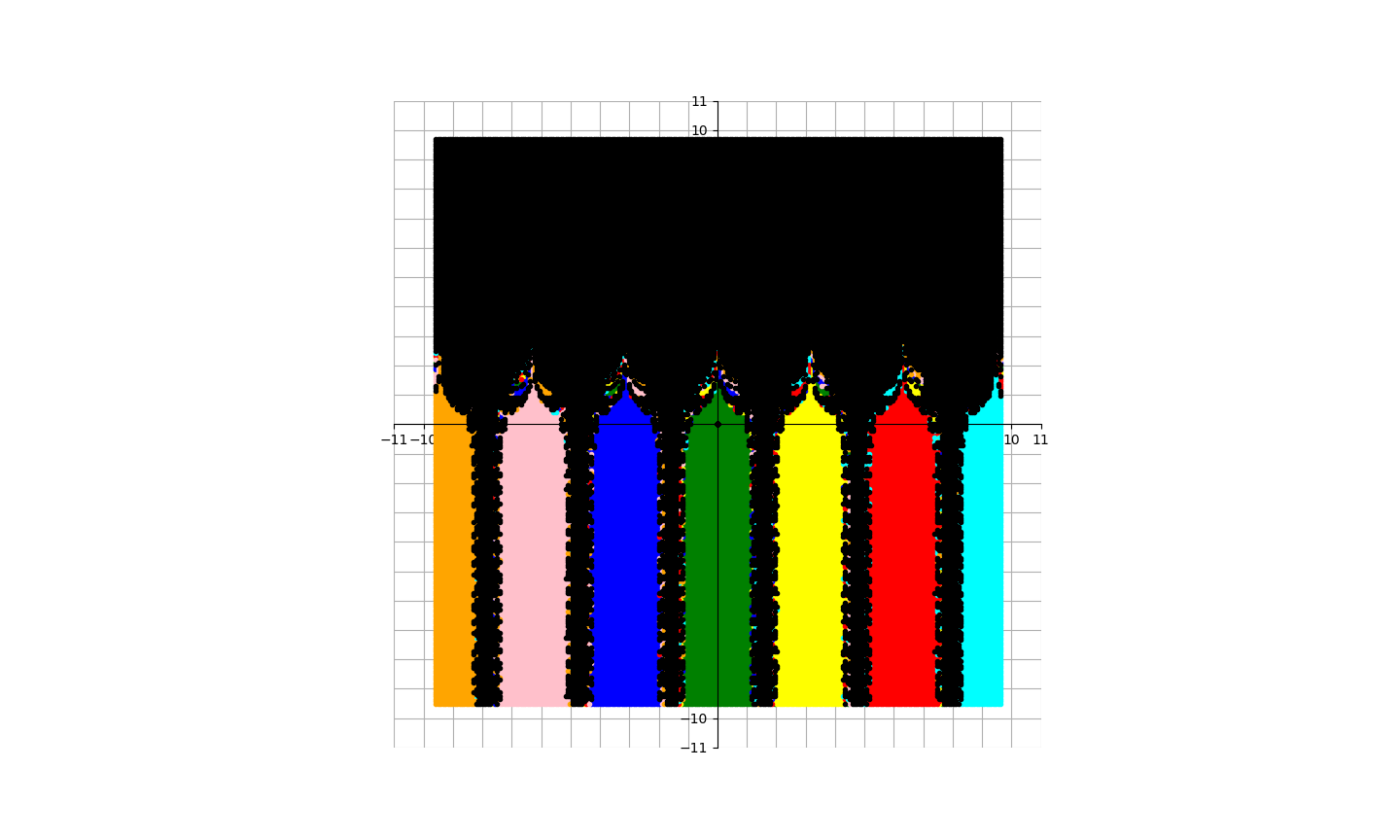}
    \includegraphics[width=10cm]{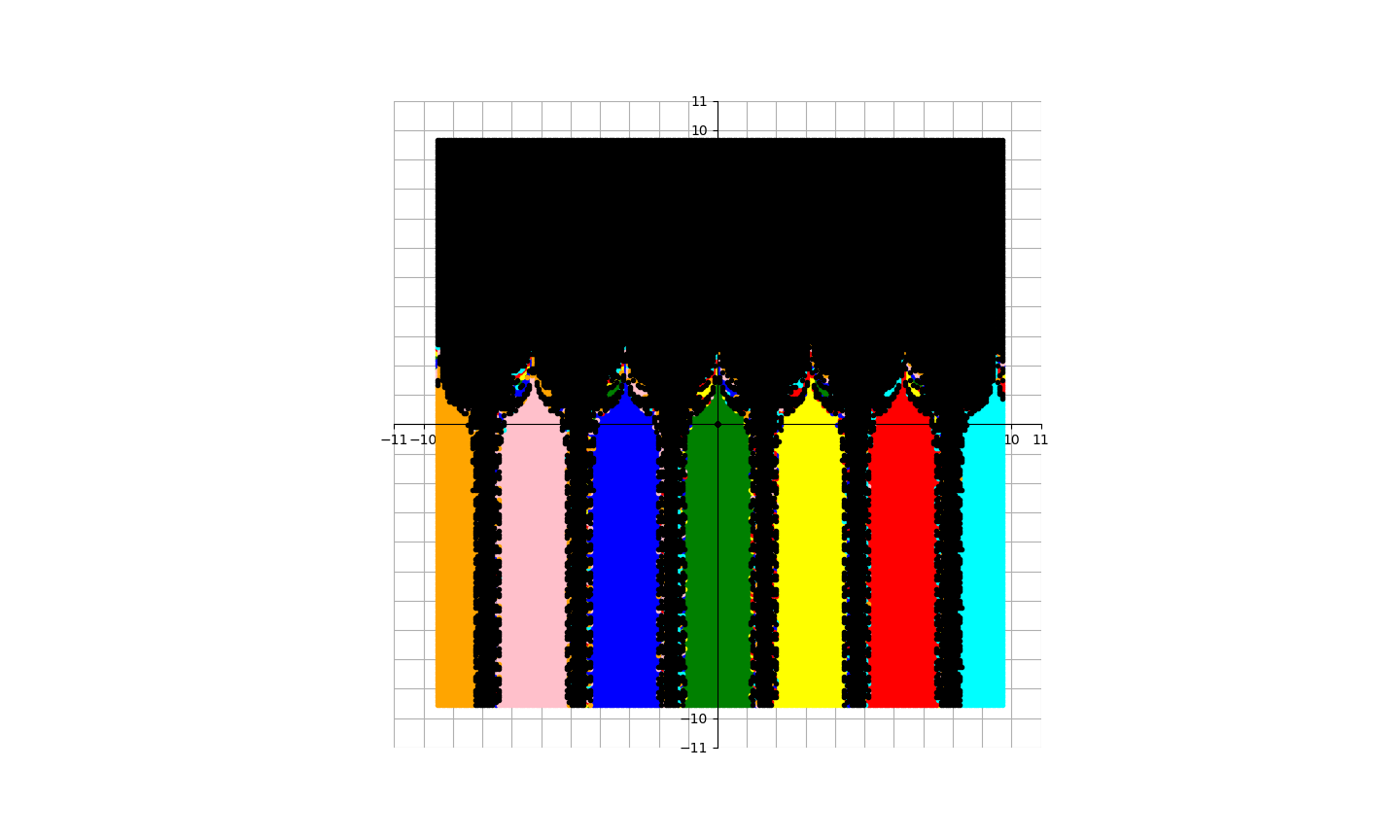}
    \caption{Basins of attraction for finding the 7 roots inside the domain $|z|<10$ of $f(z)=e^{2iz}-1$, an entire function with no critical points. Top picture: Newton's method, bottom picture: Random Relaxed Newton's methods. Points of the same colour  belong to basin of attraction of the same root. Points in black do not belong to the basin of attraction of any of the concerned 7 roots. Note that the 2 pictures look very similar.}
    \label{fig:BasinOfAttractionExponentialFunction}
\end{figure}

\begin{figure}
    \centering
    \includegraphics[width=10cm]{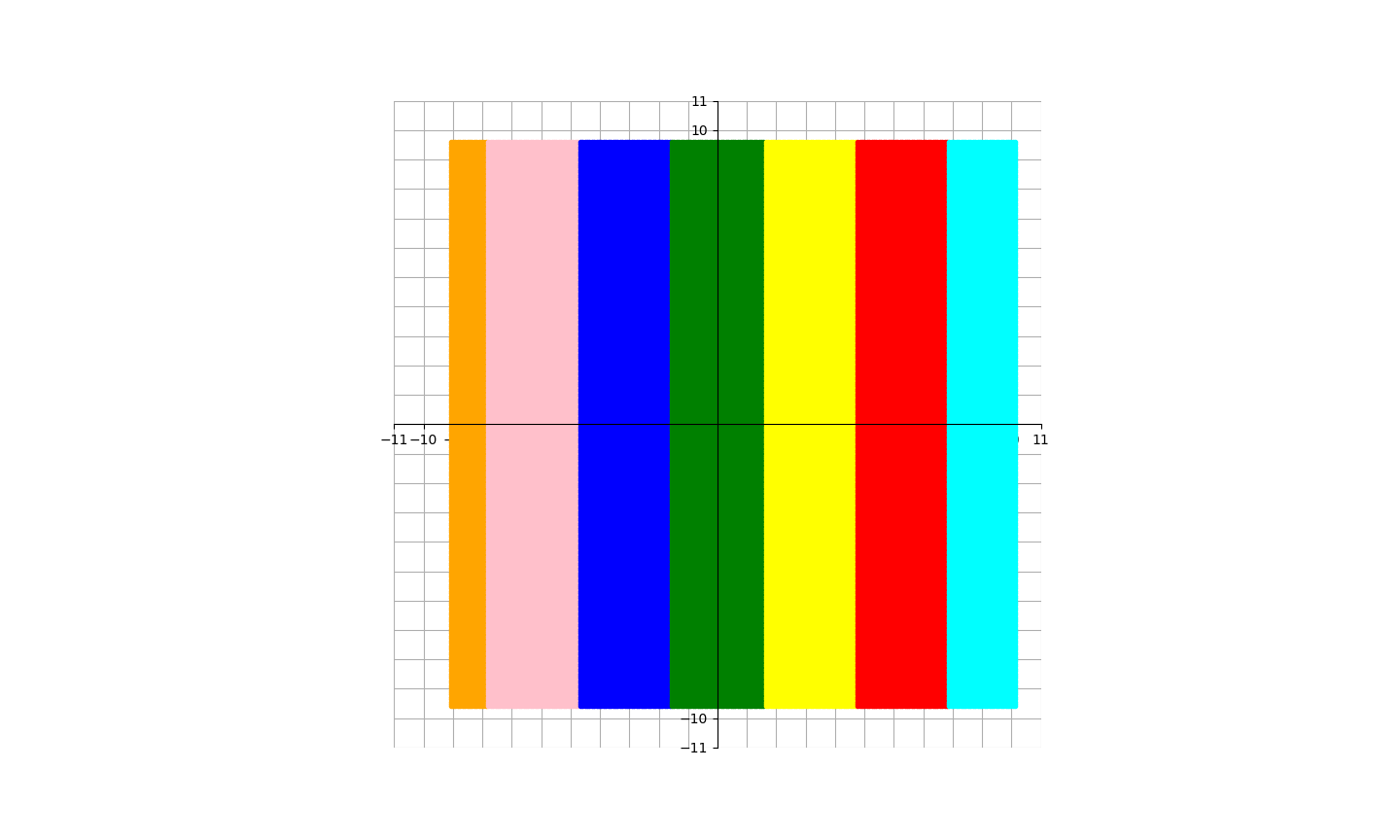}
   
    \includegraphics[width=10cm]{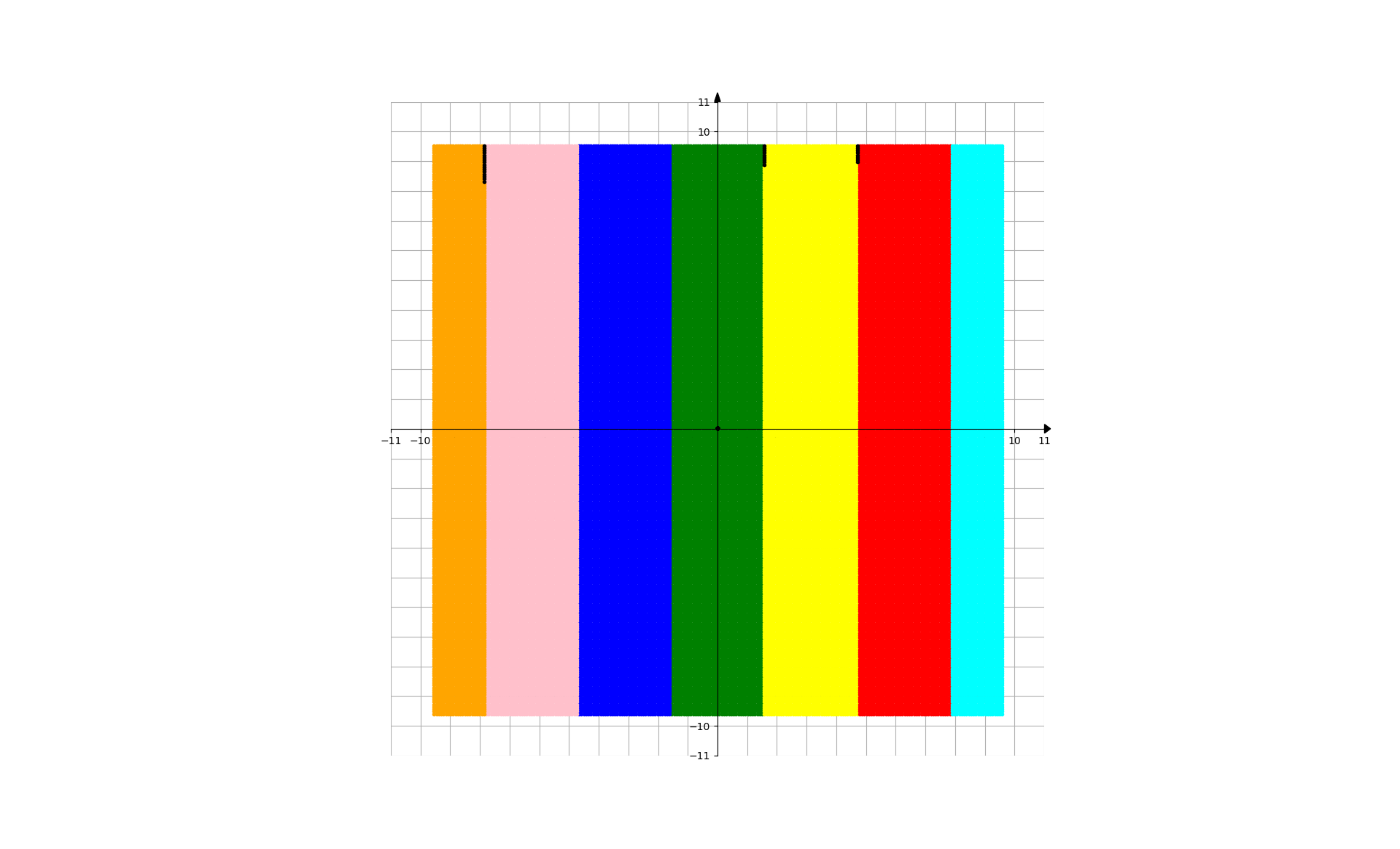}
    \caption{Basins of attraction for finding the 7 roots inside the domain $|z|<10$ of $f(z)=e^{2iz}-1$, an entire function with no critical points. Top picture: Voronoi's diagram of the 7 roots, bottom picture: BNQN. Points of the same colour  belong to basin of attraction of the same root. Points in black do not belong to the basin of attraction of any of the concerned 7 roots. Here $\theta =1$ in Algorithm \ref{table:alg0}. Note that the 2 pictures look very similar.}
    \label{fig:BasinOfAttractionExponentialFunction1}
\end{figure}

 Figures \ref{fig:BasinOfAttractionExponentialFunction} and \ref{fig:BasinOfAttractionExponentialFunction1} compare the performance of Newton's method, Random Relaxed Newton's method and BNQN for finding the 7 roots, within the domain $|z|\leq 10$, of $f(z)=e^{2iz}-1$, an entire function without critical points. We also draw a picture for Voronoi's diagram of the corresponding 7 roots.

\begin{figure}[!htb]
	\centering
\includegraphics[width=0.7\textwidth,height=0.4\textheight]{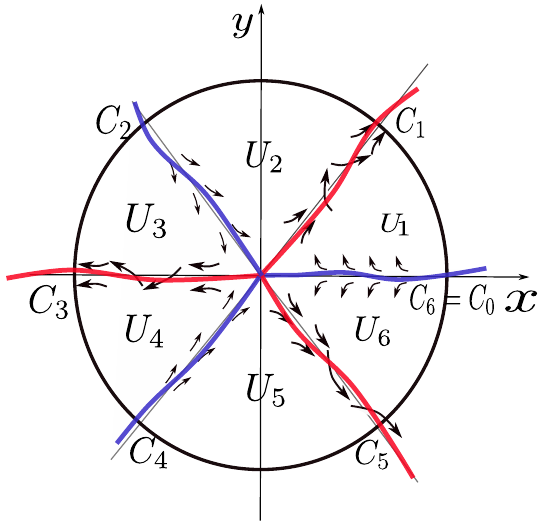}
	\caption{An illustration of Theorem \ref{Theorem2}: The orbits of BNQN  near a critical point of order $d=3$. Here the local dynamics is that of a 3-fold saddle point.}
	\label{fig3BNQNdSaddlePoint}
\end{figure}

The second main result of this paper is the following, which is the main tool to prove Theorem 
\ref{Theorem1}. For $r>0$, we use the usual notation $\mathbb{D}(0,r)=\{z\in \mathbf{C}: |z|<r\}$. We refer to Section 2 for some details on stable and unstable manifolds/curves of dynamical systems. In the statement of the theorem, $arg(z)$ is the branch of the argument function with values in $(0,2\pi ]$. See Figure \ref{fig3BNQNdSaddlePoint} for an illustration of Theorem \ref{Theorem2} when $d=3$. 

\begin{thm} Let $f(z)$ be a non-constant meromorphic function, and $z^*=0$ a critical point of $f$ of order $d$, which is not a root of $f$ (i.e. the Taylor series of $f$ in a neighborhood of $0$ has the following form $f(z)=a_0+a_dz^d+$ higher order terms, where $0\not= a_0,a_d\in \mathbf{C}$ and $d\geq 2$). For simplicity, we fix the form to be $f(z)=1+z^d+$ higher order terms, the general form can be treated similarly (see an explanation in the beginning of Section 3). 

The point $z^*$ is a d-fold saddle point for the dynamics of BNQN. More precisely, there are $r_0>0$, and continuous curves $C_1,\ldots ,C_{2d}$ in $\mathbb{D}(0,r_0)$, with the following properties: 

1. $C_j\cap C_{j'}=\{0\}$ if $j,j'=1,\ldots 2d$ and $j\not=j'$. 

2. For all $j=1,\ldots 2d$, $C_j$ can be extended to a closed smooth real analytic curve $\widehat{C_j}$ in $\mathbb{D}(0,r_0)$. Moreover, $\widehat{C_j}$ is tangent to the line $\{z\in \mathbf{C}: arg(z)=\frac{\pi j}{d}\}$ at $0$.  

3. For $j=1,\ldots d$, $C_{2j}$ are stable curves for the dynamics of BNQN, and $C_{2j-1}$ are unstable curves for the dynamics of BNQN. 

More precisely, if $z_0\in \mathbb{D}(0,r)$, and $\{z_n\}$ is the sequence constructed by BNQN with initial point $z_0$, then: 

3a) The curves $C_{2j}$ are invariant (locally) by the dynamics of BNQN, and the restriction of BNQN to $C_{2j}$ is attracting. Moreover,  $z_0\in  \bigcup _jC_{2j}$ iff $z_n\in \mathbb{D}(0,r_0)$ for all $n$. In this case, $\lim _{n\rightarrow \infty}z_n=z^*$. 

3b) The curves $C_{2j-1}$ are invariant (locally) by the dynamics of BNQN, and the restriction of BNQN to $C_{2j-1}$ is repelling. 

3c) For $j=1,\ldots ,d-1$, let $U_j$ be the open domain bound by the 2 stable curves $C_{2j}$ and $C_{2j+2}$, and contains the unstable curve $C_{2j+1}$. Define $U_{d}$ be the open domain bound by the 2 stable curves $C_{2d}$ and $C_{2}$, and contains the unstable curve $C_{1}$. If $z_0\in U_j$, then there is $N$ so that all $z_0,\ldots ,z_N$ belong to $U_j$, while $z_{N+1}\notin \mathbb{D}(0,r_0)$. 
\label{Theorem2}\end{thm}

\begin{rem}{\bf Some remarks on Theorem \ref{Theorem2}:} 

 1) The proofs for the cases $d=2$ and $d\geq 3$ are significantly different. 

{\bf For the case $d\geq 3$:}

This consists of the main bulk of this paper. In this case, the dynamics of BNQN is not $C^1$, but continuous, at the critical point $z^*$. (Recall that in general the dynamics of BNQN is defined everywhere, except at poles of $f$, but not continuous. Its fixed points are precisely the zeros and critical points of $f$.) The reason is that the associated map of BNQN involves division by a power of $|z|$, which is not a real analytic function near $z^*=0$. Moreover, both eigenvalues of the involved matrix $\nabla ^2F(z)+\delta _j||\nabla F(z)||^{\tau}Id$ (in the definition of BNQN in Section 2) are 0 at $0$. Therefore, the usual unstable/stable manifold results are not directly applicable. To proceed, we first show that in a small open neighborhood of $z^*$, the parameter $\delta _j$ in the definition of BNQN is $\delta _0$, and (more difficultly) the learning rate $\gamma _n$ is $1$. Then, we restrict the consideration to different sectors  $\frac{\pi j}{d}-\epsilon _0 <arg(z)<\frac{\pi j}{d}+\epsilon _0$, for some small number $\epsilon _0>0$. We then observe that an appropriate transformation makes the associated map into another map where stable/unstable manifold theorems are applicable. Pulling back to the origin coordinates and using some special properties of the dynamics (including the ''No jump lemma'' in Section 2, applicable to the small sectors around stable curves) we are able to complete the proof of Theorem \ref{Theorem2}. Potentially, there is also another stable curve $E_{2j-1}$ in the same sector as $C_{2j-1}$ (while $C_{2j-1}$ is the pullback of the unstable manifold, $E_{2j-1}$ is the pullback of the stable manifold), but we can exclude by some special consideration of the associated map near these small sectors. The complete arguments for the proof of the case $d\geq 3$ are rather complicated, hence we add - as a warming up aid - a simple proof for the simpler case of the linearization of the associated map. 

{\bf For the case $d=2$}: 
 
 In \cite{RefTT}, for BNQN, it was shown that in this case the point $z^*$ is a usual saddle point of the function $F$, and hence the classical theory of center-stable manifolds (see e.g. \cite{shub}) show the existence of $C^1$ manifolds $\widehat{C_1},\widehat{C_2}$ (in this case $C_1,C_3$ are parts of $\widehat{C_1}$, and $C_2,C_4$ are parts of $\widehat{C_3}$, and the curves $E_1,E_2,E_3,E_4$ are empty) satisfying parts 1 and 2, 3a and 3b of the theorem. The work in \cite{RefTT} - based on \cite{RefT},  which is  developed for a general $C^3$ cost function $F$ - shows only that the local stable/unstable manifolds are $C^1$, and the exceptional set $\mathcal{E}$ has Lebesgue measure 0. Moreover, it is required that the parameters in the algorithm BNQN are random.  Theorem \ref{Theorem2} improves upon these points, and also proves the extra information in part 3c. Here, like the case $d\geq 3$, the ''No jump lemma" mentioned above is used at stable curves. 

For the special case where $f(z)=1+z^2$ (indeed, any polynomial of degree $2$ with 2 distinct roots) is a polynomial of degree $2$, the global dynamics of BNQN has been treated in \cite{RefFHTW}, where an analog of the classical result of Schr\"oder for Newton's method is proven. In this case, it turns out that the curves $C_1$, $C_3$ are halves of the y axis, and $C_2$ and $C_4$ are halves of the x axis. All parts of Theorem \ref{Theorem2} hold. 

For the case $d=2$ of Theorem \ref{Theorem2}, the linearization of the associated map is $\Phi _1(x,y)=(0,\lambda y)$, and the transformation used in the case $d\geq 3$ does not help (and which is not needed anyway, since the map $\Phi$ is already real analytic). 

2)  The behaviour of BNQN described in Theorem \ref{Theorem2} is very similar to that of Newton's flow. The latter can be explained through the classical Poincar\'e-Bendixon theorem \cite{RefPoincare}\cite{RefBendixon}.  We will discuss more on this in Section 4.4.
\end{rem}

If $f(z)$ has a root $z^*=0$ of order $d\geq 1$ (i.e. in a neighbourhood of $0$, $f$ has Taylor expansion $f(z)=az^d+$higher order terms), then near $0$ the dynamics of Newton's method is
\begin{eqnarray*}
z_{n+1}=\frac{d-1}{d}z_n+O(|z_n|^2). 
\end{eqnarray*}
In particular, if $d=1$ (i.e. the root $z^*=0$ is simple), then we recover the well known fact that Newton's method has quadratic local rate of convergence near $z^*=0$. As mentioned, when $z^*=0$ is a simple root, then BNQN also has quadratic local rate of convergence near $z^*$. It is known (see \cite{RefThuanTuyen}) that every (multiple) root of $f$ is an attractor of the dynamics of BNQN. The last main result of this paper provides a quantitative description of this fact. 

\begin{thm} Let $f(z)$ be a non-constant meromorphic function and $z^*=0$ is a root of $f$ of order $d$. Then near $0$, the dynamics of BNQN  is as follows: 
\begin{eqnarray*}
z_{n+1}=\frac{2d-2}{2d-1}z_n+O(|z_n|^2)=\frac{d-1}{d-0.5}z_n+O(|z_n|^2).
\end{eqnarray*}
\label{Theorem3}\end{thm}

As a byproduct of the proof of Lemma \ref{LemmaApproximateFormulaDynamics}, since $z_k-v_k$ in that lemma is (up to an $O(|z_k|^2)$ error) the update rule for Newton's method for Optimization (for detail please see e.g. \cite{RefFHTW2}) and it is equal to $$\frac{(d-2)z_k}{d-1}+O(|z_k|^2),$$ it is rigorously seen here that for Newton's method for Optimization applied to $F(x,y)=|f(x+iy)|^2$, the critical points (no matter if it is a root or not) of $f(z)$ are also attractors, as observed experimentally in \cite{RefFHTW2}. Therefore, like Newton's method, Newton's method for Optimization also does not have global convergence guarantee for finding roots of polynomials in 1-complex variable.  

{\bf Plan of the paper.} In Section 2 we present some preliminary results. In Section 3 we prove Theorem \ref{Theorem2}, and in Section 4 we prove Theorems \ref{Theorem1}, Theorem \ref{Theorem1Bis} and \ref{Theorem3}. Also, in Section 4 we discuss how the results in this paper can be used to explain (both rigorously and heuristically) the observations in many experiments on BNQN, in particular those in \cite{RefFHTW2}, in particularly connections to Newton's flow and the classical Poincar\'e-Bendixon theorem. 

{\bf Acknowledgments.} We would like to thank Takayuki Watanabe for sharing the experiments concerning finding roots of the 2-dim polynomial map $h(z_1,z_2)=(z_1+z_1z_2,z_2^2-z_1^2)$, and for the comment on the non-existence of ''channel to infinity" in experiments for Newton's method applied to certain rational functions $f$ with compact supports. The first and second authors would like to thank the Department of Mathematics, University of Oslo, for hospitality. The second and third authors are partially supported by Research Council of Norway's Young research talents grant 300814, as well as Trond Mohn Foundation. 

\section{Some preliminary results} This section recalls the definition and some first properties of BNQN. We also collect a couple of facts concerning local stable/unstable manifolds for real analytic maps in dimension 2, needed in later proofs. 

\subsection{BNQN and some first properties} BNQN is originally designed in \cite{RefT}, as an improvement of the precursor New Q-Newton's method (NQN) in \cite{RefTT}. BNQN is a modification of Newton's method for Optimization, the latter being Newton's method applied to find critical points of a function $F:\mathbb{R}^m\rightarrow \mathbb{R}$. If $F:\mathbf{R}^m\rightarrow \mathbf{R}$ is a function, then the update rule for Newton's method for Optimization is:
$$z_{n+1}=z_n-(\nabla ^2F(z_n))^{-1}.\nabla F(z_n),$$
provided the Hessian matrix $\nabla ^2F(z_n)$ is invertible. 

Like Newton's method, Newton's method for Optimization also has quadratic local rate of convergence near non-degenerate local minima, but has no global convergence guarantee. Moreover, Newton's method for Optimization cannot avoid saddle points (see \cite{RefFHTW2} for experimental results, and see the comments at the end of Section 1 for a theoretical justification). NQN has also quadratic local rate of convergence near non-degenerate local minima, and in addition also can avoid saddle points. However, NQN has no global convergence guarantee. The design of BNQN is to preserve the good properties of NQN, while also has better global convergence guarantee for a large class of functions of interest (including functions which have at most countably many critical points). We refer the reader to the mentioned papers for more detail. 

In this paper we will use  the following modification of the definition in \cite{RefFHTW} (the latter being different to the original version in \cite{RefT} only at the introduction of the parameter $\theta$, which is used in Theorem \ref{Theorem2Linearization}). The only difference of the version in this paper to the version in \cite{RefFHTW} is as follows.  We require $\tau \in 2\mathbf{N}$ (where $\mathbf{N}=\{1,2,3,\ldots \}$ is the set of positive integers), instead of $\tau >0$. 

For a symmetric, square real matrix $A$, we define: 
  
  $sp(A)=$ the maximum among $|\lambda |$'s, where $\lambda  $ runs in the set of eigenvalues of $A$, this is usually called the spectral radius in the Linear Algebra literature;
  
  and 
  
  $minsp(A)=$ the minimum among $|\lambda |$'s, where $\lambda  $ runs in the set of eigenvalues of $A$, this number is non-zero precisely when $A$ is invertible.
  
 One can easily check the following more familiar formulas: $sp(A)=\max _{\|e\|=1}\|Ae\|$ and $minsp(A)=\min _{\|e\|=1}\|Ae\|$, using for example the fact that $A$ is diagonalisable.

For the reader's convenience, here we explain the main components in Algorithm \ref{table:alg0}. First, a modification of the Hessian matrix $\nabla ^2F(z_k)$ is done, by adding the term $\delta _j||\nabla F(z_k)||^{\tau}$. The first While loop in the algorithm makes sure that the new matrix is sufficiently invertible, which helps in establishing that any cluster point of the constructed sequence is a critical point of $F$. The $v_k$ can be thought of as the one obtained if one uses Newton's method (where the matrix $A_k$ is thought of as the Hessian matrix). The $w_k$ is obtained by changing the signs of negative values of the matrix $A_k$, this helps to avoid saddle points. The second While loop is Armijo's Backtracking line search \cite{RefAr}, which boosts global convergence guarantee.  

\medskip
{\color{blue}
 \begin{algorithm}[H]
\SetAlgoLined
\KwResult{Find a minimum of $F:\mathbf{R}^m\rightarrow \mathbf{R}$}
Given: $\{\delta_0,\delta_1,\ldots, \delta_{m}\} \subset \mathbf{R}$,\, $\tau \in 2\mathbf{N}$, $0\leq\theta $, and $0<\gamma _0\leq 1$;\\
Initialization: $z_0\in \mathbf{R}^m$\;
$\kappa:=\frac{1}{2}\min _{i\not=j}|\delta _i-\delta _j|$;\\
 \For{$k=0,1,2\ldots$}{ 
    $j=0$\\
  \If{$\|\nabla F(z_k)\|\neq 0$}{
   \While{$minsp(\nabla^2F(z_k)+\delta_j \|\nabla F(z_k)\|^{\tau}Id)<\kappa  \|\nabla F(z_k)\|^{\tau}$}{$j=j+1$}}
  
 $A_k:=\nabla^2F(z_k)+\delta_j \|\nabla F(z_k)\|^{\tau}Id$\\
$v_k:=A_k^{-1}\nabla F(z_k)=pr_{A_k,+}(v_k)+pr_{A_k,-}(v_k)$\\
$w_k:=pr_{A_k,+}(v_k)-pr_{A_k,-}(v_k)$\\
$\widehat{w_k}:=w_k/\max\{1,\theta \|w_k\|\}$.\\
(If $F$ has compact sublevels, then one can choose $ \widehat{w_k}=w_k$, or equivalently $\theta =0$.)\\

$\gamma :=1$\\
 \If{$\|\nabla F(z_k)\|\neq 0$}{
   \While{$F(z_k-\gamma \widehat{w_k})-F(z_k)>-\gamma \langle\widehat{w_k},\nabla F(z_k)\rangle/3$}{$\gamma =\gamma /3$}}

$z_{k+1}:=z_k-\gamma \widehat{w_k}$
   }
  \caption{Backtracking New Q-Newton's method (BNQN)} \label{table:alg0}
\end{algorithm}
}
\medskip

In \cite{RefFHTW} the following two theorems can be found. 

\begin{thm}[A convergence result for root finding] Let $f(z):\mathbf{C}\rightarrow \mathbf{P}^1$ be a non-constant meromorphic function. 
Define a function $F:\mathbf{R}^2\rightarrow [0,+\infty]$ by the formula $F(x,y)=|f(x+iy)|^2/2$. 

Given an initial point $z_0\in \mathbf{C}\backslash \mathcal{P}(f)$, we let $\{z_n\}$ be the sequence constructed by BNQN  applied to the function $F$ with initial point $z_0\in \mathbf{C}\backslash \mathcal{P}(f)$. 

A) Either $\{z_n\}$ converges to a root of $f(z)f'(z)$ (i.e. a point in $\mathcal{Z}(f)\cup \mathcal{C}(f)$), or it converges to $\infty$. If $f$ has compact sublevels, then the latter case cannot happen.  

B) Assume that $f$ is generic in the sense that both $f$ and $f'$ have only simple roots (in other words,  $\{z\in \mathbf{C}:f(z)f"(z)=f'(z)=0\}=\emptyset$). Assume also that the parameters in BNQN are randomly chosen. If  $\mathcal{E}$ is the exceptional set, consisting of initial points $z_0$ for which $\{z_n\}$ converges to a point in $\mathcal{C}(f)$, then $\mathcal{E}$ has Lebesgue measure $0$.   
\label{TheoremConvergence}\end{thm}

\begin{thm}[Invariant under appropriate linear changes of coordinates] The dynamics of BNQN is invariant under conjugation by linear operators of a certain form. The precise statement is as follows. 

Let $F:\mathbf{R}^m\rightarrow \mathbf{R}$ be an objective function. Let $A:\mathbf{R}^m\rightarrow \mathbf{R}^m$ be an invertible linear map, of the form $A=cR$, where $R$ is  {\bf unitary} (i.e. $RR^T=Id$, where $R^T$ is the transpose of $R$) and $c>0$. Define $G(z)=F(Az)$. 

Let $z_0\in \mathbf{R}^m$ be an initial point, and let $\{z_n\}$ be the sequence constructed by BNQN (with parameters $\delta _0,\delta _1,\ldots ,\delta _m, \theta, \tau$) applied to the function $F(z)$. 

Let $z_0'=A^{-1}z_0$, and let $\{z_n'\}$ be the sequence constructed by BNQN (with parameters $\delta _0',\delta _1',\ldots ,\delta _m', \theta ', \tau $) applied to the function $G(z)$. 

Assume that $\delta_i'=\delta_i c^{2-\tau}$ for all $0\leq i\leq m$, and $\theta '=c\theta $.

Then for all $n$, we have $z_n'=A^{-1}z_n.$  
\label{TheoremLinearConjugacyInvariant}
\end{thm}

\subsection{Some results on stable/unstable manifolds of real analytic maps in dimension 2} We present here some necessary results concerning stable/unstable manifolds, which are needed in later proofs. Here we consider $U\subset \mathbf{R}^2$ an open neighbourhood of $0$, and $H:U\mapsto \mathbf{R}^2$ a real analytic map having $(0,0)$ as an isolated fixed point. 

The next two theorems should be well known to experts. 

For the first result, there are proofs for invertible germs, but we cannot find a reference where it is proved for non-invertible germs. We will provide a sketch of proof for it, sufficiently self-contained, following those in \cite[Chapters 5 and 6]{shub} (differentiable case) and \cite[Section 6.4]{RefMNTU} (complex analytic case). 

\begin{thm}[Stable/Unstable manifolds for real analytic maps] Assume that $H$ is real analytic and the Jacobian $JH(0)$ has two eigenvalues $\lambda _1,\lambda _2\in \mathbf{R}$ satisfying $0\leq |\lambda _1|<1<|\lambda _2|$. Then there exists an open neighbourhood $V\subset \mathbf{R}^2$ of $0$ and two smooth real analytic curves $C_1$ and $C_2$ of $V$ with the following properties: 

1) $C_1$ is the stable manifold of $H$: is it is invariant by $H$ and the restriction of $H$ to $C_1$ is attracting. Moreover, for $z_0\in V$ the following two statements are equivalent: 1a) the sequence $\{H^n(z_0)\}$ stays in $V$, and 1b) $z_0\in C_1$. 

2) $C_2$ is the unstable manifold of $H$: the restriction of $H$ to $C_2$ is repelling. Moreover, if $H$ is invertible, then $C_2$ is the stable manifold of $H^{-1}$.   

\label{TheoremStableUnstableManifoldRealAnalytic}\end{thm}
\begin{proof}[Sketch of proof] The idea of the proof is to first complexify the map, prove the result in the complex setting using complex dynamics, then pullback to the real setting. 

We can assume that the stable vector is $(1,0)$ and the unstable vector is $(0,1)$. 

First, complexify the real variables $x,y$ to complex variables $X,Y$, and shrink the domain if necessary, we obtain a complex analytic map $\widehat{H}: \widehat{U}\rightarrow \mathbf{C}^2$. The Jacobian of $\widehat{H}$ at $(0,0)$ has the same eigenvalues and eigenvectors as $H$. By \cite[Theorem 6.4.3, case ii]{RefMNTU}, $\widehat{H}$ has an unstable manifold $\widehat{C}$, which is given by an equation of the form $X=\phi (Y)$. Take the intersection of this graph with $\mathbf{R}^2$ we obtain the unstable manifold $C_2$ for $H$, which is real analytic. Moreover, we can, after a linear change of coordinates, assume that the unstable manifold of $\widehat{H}$ is the $Y$ axis, and $\widehat{H}(0,Y)=(0,\lambda _2Y)$.  

Second, if $0<|\lambda _1|$, then the same argument (applied to $H^{-1}$) gives that the stable manifold $C_1$ which is real analytic. 

Third, now we consider the case $\lambda _1=0$. We will construct first the stable manifold of $\widehat{H}$ as  a function $Y=g(X)$. To this end, we first show that: for each $X_0$, there is a unique $Y_0$ such that the orbit of $(X_0,Y_0)$ stays inside the neighbourhood $V$ where $\widehat{H}$ is defined. Since $\widehat{H}(0,Y)=(0,\lambda _2Y)$ and $\lambda _1=0$, we can write $\widehat{H}(X,Y)=(X\varphi (X,Y), \lambda _2Y+X\psi (X,Y))$. Consider a closed rectangle $R=[a_1,a_2]\times [b_1,b_2]\subset V$. Let $I_1=\{X=X_0\}\cap R$. Then $\widehat{H}(I_1)$ will roughly be dilated by a factor of $|\lambda _2|>1$. Hence the size of  $I_2=\widehat{H}^{-1}(\widehat{H}(I_1)\cap R)\cap I_1$ is about $1/|\lambda _2|$ of the size of $I_1$. Iterate this process (using more iterations of $\widehat{H}$), we obtain a sequence of closed subsets $\ldots I_3\subset I_2\subset I_1$, whose size decreases to $0$. Therefore, their intersection is a unique point $Y_0$, which is the unique point on the line $\{X=X_0\}$. This is characterised as the unique point so that the orbit of $(X_0,Y_0)$ belong to the rectangle $R$. This gives us a function $Y=g(X)$. We need next to show that $g$ is complex analytic. To this end, for $0<\mu <1$ small, we consider a new map $\widehat{H}_{\mu}(X,Y)=(\mu X,0)+\widehat{H}(X,Y)$. For each $\mu$, by the Second step, $\widehat{H}_{\mu}$ has a stable complex manifold $\widehat{C_{1,\mu}}$ given as the graph of a function $Y=g_{\mu}(X)$. When $\mu$ converges to $0$, the sequence of these graphs uniformly converge to the graph for $g(X)$. Since $g_{\mu}(X)$ is holomorphic, it follows that $g(X)$ is also holomorphic. Then intersect with $\mathbf{R}^2$ gives that the stable manifold is real analytic.  
 \end{proof}

The second result is folklore, but no reference for it is found, and we provide a simple proof for completeness.  

\begin{thm}[No jump lemma] Assume $\Omega\subset\mathbf{R}^2$ is a small open neighbourhood of $0$, and $h:\Omega\rightarrow\mathbf{R}^2$ is a real analytic map. Assume that $Jh(0,0)$ has two eigenvalues $a,A\in \mathbf{R}$ such that $0\leq |a|<1<|A|$. 

1) If $A>0$, then the stable manifold partitions $\Omega $ into 2 parts, each invariant under $h$.  

2) If $a>0$, then the unstable manifold partitions $\Omega $ into 2 parts, each invariant under $h$. 

3) In particular, if both $a,A>0$, then the stable and unstable manifolds of $h$ partition $\Omega$ into 4 domains, each invariant under the dynamics of $h$. 

\label{TheoremNoJumpLemma}\end{thm}
\begin{proof}
Since the 2 eigenvalues $a,A$ of $Jh(0,0)$ are different, by applying a linear change of coordinates, we can assume that the following holds. $Jh(0,0)$ is a diagonal matrix with stable eigenvector $(1,0)$, stable eigenvalue $0<a<1$,  and unstable eigenvector $(0,1)$, unstable eigenvalue $A>1$.  Using Theorem \ref{TheoremStableUnstableManifoldRealAnalytic}, we can also assume that the stable manifold is of the form $y=u(x)$ for a real analytic function $u$ with $ u(0)=u'(0)=0$, and that the unstable manifold is of the form $x=v(y)$ for a real analytic function $v$ with $v(0)=v'(0)=0$. In this case, the stable and unstable manifolds partition $\Omega$ into the following quadrants divided by the stable and unstable manifolds (see Figure \ref{fig2}):
\begin{figure}[!htb]
	\centering
\includegraphics[width=0.5\textwidth,height=0.25\textheight]{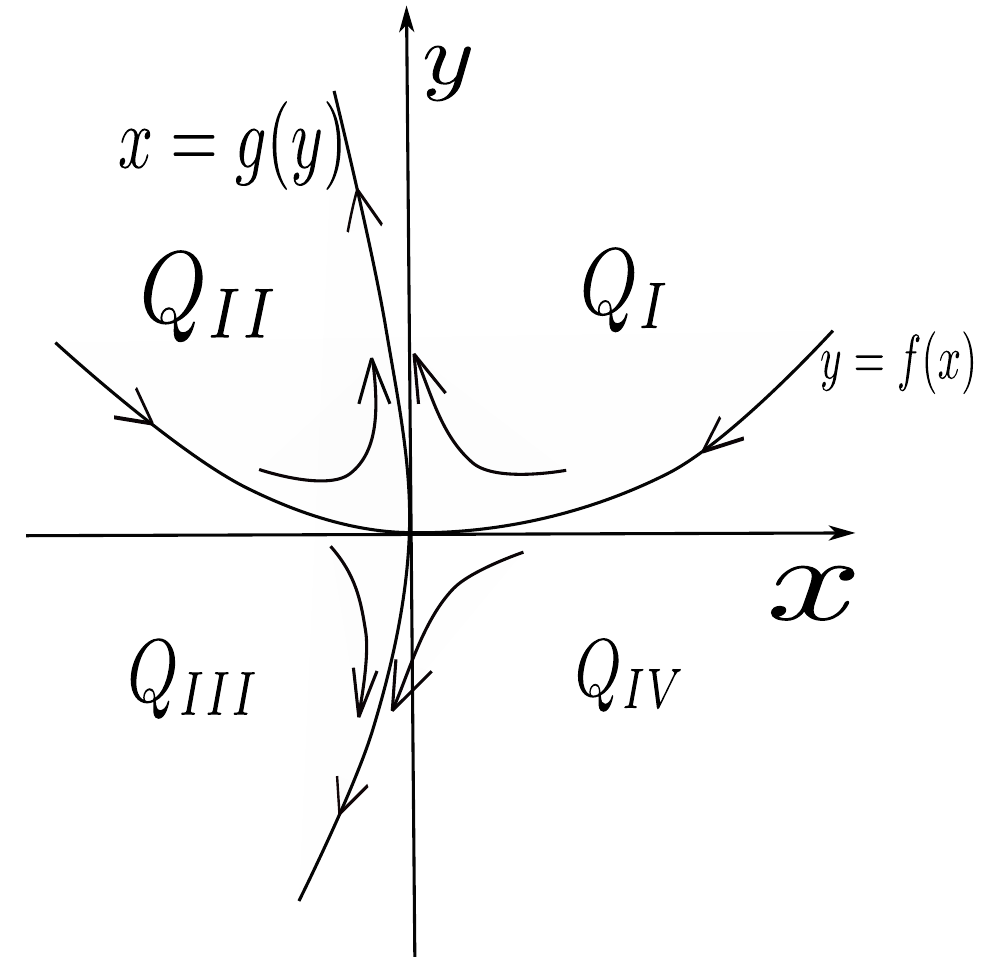}
	\caption{The stable and unstable manifolds.}
	\label{fig2}
\end{figure}
	\begin{itemize}
	\item $Q_{I}= \{ x>v(y),y>u(x)\}$
	\item $Q_{II}= \{ x<v(y),y>u(x)\}$
	\item $Q_{III}= \{ x>v(y),y<u(x)\}$
	\item $Q_{IV}= \{ x<v(y),y<u(x)\}$
	\end{itemize}

Define $W(x,y)=(x-v(y), y-u(x))=(x_1,y_1)$ Then $JW(0,0)$ is the identity matrix, $W$ is a local real analytic isomorphism which maps the stable manifold to the $x$ axis and the unstable manifold to the $y$ axis. Moreover, the images of the quadrants are the classical quadrants, $W(Q_{I})=\{x_1>0, y_1>0\}$ etc.  This conjugates $h$ to a diffeomorphism $H$, where $H= W\circ h\circ W^{-1}$. 
	
Then 
$$
	JH(0,0)= 
	\left[ \begin{array}{ccc}
		a&0 \\
		0&A\\
	\end{array} \right],
	$$

$H$ has the $x$ axis as its stable manifold and the $y$ axis as its unstable manifold. We only need to show that $H$ maps the classical quadrants to themselves. We can write $H$ in the following form: $H(x,y)= (H_1(x,y),H_2(x,y))$ where $H_1(0,y)=0$ and $H_2(x,0)=0$, $\frac{\partial H_1}{\partial x}>0$ and 
$\frac{\partial H_2}{\partial y}>0$.

1) Assume that $A>0$. If $(x,y)$ is such that $y>0$ then by the Intermediate Value Theorem, there are $0<y_1<y$  so that:  
$H_2(x,y)=H_2(x,y)-H_2(x,0)= \frac{\partial H_2}{\partial y}(x,y_1)(y-0)>0$.

This means that the 2 parts above and below the stable manifold are invariant by $H$.

2) and 3): The proofs are similar.  
\end{proof}

\section{Local behavior of BNQN near critical points but not roots of $f(z)$: Proof of Theorem \ref{Theorem2}} In this section we research the local behaviour of BNQN near a critical point $z^*=0$ but not a root of $f$, which accumulates in the proof of Theorem \ref{Theorem2}. 

Near $z^*=0$, $f$ has the following Taylor's series: 
\begin{eqnarray*}
f(z)=a_0+a_dz^d+h.o.t.,
\end{eqnarray*}
where $h.o.t.$ stands for ''higher order terms", with $0\not= a_0,a_d\in \mathbf{C}$. 

We can, either by using Theorem \ref{TheoremLinearConjugacyInvariant} or by mimicking the arguments in this section, assume throughout this section that $f(z)$ has the form: 
\begin{eqnarray*}
f(z)=1+z^d+h.o.t.
\end{eqnarray*}
near $z^*=0$. 

Then BNQN is applied to the function $F:\mathbf{R}^2\rightarrow [0,+\infty ]$, which is defined as $F(x,y)=|f(x+iy)|^2/2$. The critical points of $F$ are precisely roots of $f(z)f'(z)$, see \cite{RefTT}. 

Through out this section, we use both Cartesian and polar coordinates for a complex number: $z=x+iy=re^{i\theta}$. 

\subsection{Some local calculations} Here we present some useful local calculations near the critical point $z^*=0$ involving the function $F$ and its derivatives. 

We first note the following simple formulas, for each positive integer $N$:  
\begin{eqnarray*}
\frac{\partial z^N}{\partial x}&=&Nz^{N-1}=Nr^{N-1}e^{(N-1)\theta },\\
\frac{\partial z^N}{\partial y}&=&iNr^{N-1}e^{(N-1)\theta },\\
\frac{\partial ^2 z^N}{\partial x^2}&=&N(N-1)r^{N-2}e^{(N-2)\theta },\\
\frac{\partial ^2 z^N}{\partial x\partial y}&=&iN(N-1)r^{N-2}e^{(N-2)\theta },\\
\frac{\partial ^2 z^N}{\partial y^2}&=&-N(N-1)r^{N-2}e^{(N-2)\theta }. 
\end{eqnarray*}

Therefore, if $f(z)=1+z^d+h.o.t. = u(x,y)+iv(x,y)$ with $u,v$ real functions, we have $u=1+r^d\cos (d\theta )+O(r^{d+1})$, $v=r^d\sin (d\theta )+O(r^{d+1})$, and: 
\begin{eqnarray*}
\frac{\partial u}{\partial x}&=&dr^{d-1}\cos ((d-1)\theta )+O(r^{d}),\\
\frac{\partial v}{\partial x}&=&dr^{d-1}\sin ((d-1)\theta )+O(r^{d}),\\
\frac{\partial u}{\partial y}&=&-dr^{d-1}\sin ((d-1)\theta )+O(r^{d}),\\
\frac{\partial v}{\partial y}&=&dr^{d-1}\cos ((d-1)\theta )+O(r^{d}),\\
\frac{\partial ^2 u}{\partial x^2}&=&d(d-1)r^{d-2}\cos ((d-2)\theta )+O(r^{d-1}),\\
\frac{\partial ^2 v}{\partial x^2}&=&d(d-1)r^{d-2}\sin ((d-2)\theta )+O(r^{d-1}),\\
\frac{\partial ^2 u}{\partial x\partial y}&=&-d(d-1)r^{d-2}\sin ((d-2)\theta )+O(r^{d-1}),\\
\frac{\partial ^2 v}{\partial x\partial y}&=&d(d-1)r^{d-2}\cos ((d-2)\theta )+O(r^{d-1}),\\
\frac{\partial ^2 u}{\partial y^2}&=&-d(d-1)r^{d-2}\cos ((d-2)\theta )+O(r^{d-1}),\\
\frac{\partial ^2 v}{\partial y^2}&=&-d(d-1)r^{d-2}\sin ((d-2)\theta )+O(r^{d-1}). 
\end{eqnarray*}

The above calculations can also be found in \cite{RefTT}. 

Therefore, since $F=(u^2+v^2)/2$, the gradient and Hessian of $F$ can be approximated as follows, here $u_x=\partial u/\partial x$, $u_{xy}=\partial ^2u/\partial x\partial y$ and so on: 
$$
	\nabla F= (uu_x+vv_x,uu_y+vv_y)=
	\left[ \begin{array}{ccc}
		dr^{d-1}\cos ((d-1)\theta )+O(r^{d}) \\
		-dr^{d-1}\sin ((d-1)\theta )+O(r^{d})\\
	\end{array} \right],
	$$
and
\begin{eqnarray*}
	&&\nabla^2 F= \left[ \begin{array}{ccc}
		uu_{xx}+vv_{xx}+u_x^2+v_x^2&uu_{xy}+vv_{xy}+u_xu_y+v_xv_y \\
		uu_{xy}+vv_{xy}+u_xu_y+v_xv_y&uu_{yy}+vv_{yy}+u_y^2+v_y^2\\
	\end{array} \right]\\
    &=&d(d-1)r^{d-2}\left[ \begin{array}{ccc}
		\cos ((d-2)\theta )&-\sin ((d-2)\theta ) \\
		-\sin ((d-2)\theta )&-\cos ((d-2)\theta )\\
	\end{array} \right]+O(r^{d-1})
			.
\end{eqnarray*}

\begin{lem} The Hessian matrix $\nabla ^2F(z)$ has the following eigenvalue/eigenvector pairs: 

Eigenvalue $\lambda _1=d(d-1)r^{d-2}+O(r^{d-1})$ with eigenvector 
$$
	U_1= 
	\left[ \begin{array}{ccc}
		\cos ((d-2)\theta /2)+O(r) \\
		-\sin ((d-2)\theta /2)+O(r)\\
	\end{array} \right],
	$$

Eigenvalue $\lambda _2=-d(d-1)r^{d-2}+O(r^{d-1})$ with eigenvector 
$$
	U_2= 
	\left[ \begin{array}{ccc}
		\sin ((d-2)\theta /2)+O(r) \\
		\cos ((d-2)\theta /2)+O(r)\\
	\end{array} \right].
	$$

\label{LemmaEigenHessian}\end{lem}
\begin{proof}
This can be easily done by first working with the following matrix
$$
d(d-1)r^{d-2}\left[ \begin{array}{ccc}
		\cos ((d-2)\theta )&-\sin ((d-2)\theta ) \\
		-\sin ((d-2)\theta )&-\cos ((d-2)\theta )\\
	\end{array} \right],
$$
which is the main term of the Hesian matrix. 
\end{proof} 

\subsection{Approximate formulas for the map associated to BNQN} In this subsection we determine approximate formulas for the map associated to BNQN. 

\begin{lem} Assume that $|z|$ is small enough. Then in the definition of BNQN, see Algorithm \ref{table:alg0}, here $r=|z_k|$: 

\begin{eqnarray*}
A_k&=&\nabla ^2F(z_k)+\delta _0||\nabla F(z_k)||^{\tau}Id\\
&=&d(d-1)r^{d-2}\left[ \begin{array}{ccc}
		\cos ((d-2)\theta )&-\sin ((d-2)\theta ) \\
		-\sin ((d-2)\theta )&-\cos ((d-2)\theta )\\
	\end{array} \right]+O(r^{d-1}),\\
v_k&=&\frac{r}{d-1}\left[ \begin{array}{ccc}
		\cos (\theta ) \\
		\sin (\theta )\\
	\end{array} \right]+O(r^2)=\frac{z_k}{d-1}+O(r^2),\\
w_k&=&\frac{r}{d-1}\left[ \begin{array}{ccc}
		\cos ((d-1)\theta ) \\
		-\sin ((d-1)\theta )\\
	\end{array} \right]+O(r^2),\\
\widehat{w_k}&=&w_k,\\
z_{k+1} &=&z_k-w_k. 
\end{eqnarray*}
\label{LemmaApproximateFormulaDynamics}\end{lem}
\begin{proof}
To show that the matrix $A_k$ has the given form, first, by Algorithm \ref{table:alg0} we need to show that $minsp(\nabla ^2F(z_k)+\delta _0||\nabla F(z_k)||^{\tau}Id)\geq \kappa ||\nabla F(z_0)||^{\tau}$. By Lemma \ref{LemmaEigenHessian}, the Hessian matrix $\nabla ^2F(z_k)$ has eigenvalues $\pm d(d-1)r^{d-2}+O(r^{d-1})$, while by the calculation in Section 3.1 we have $||\nabla F(z_k)||^{\tau}=O(r^{\tau (d-1)})$. Therefore, since $\tau >1$, the needed inequality is satisfied when $r$ is small enough. Then $A_k=\nabla ^2F(z_k)+O(r^{\tau (d-1)})$, and the formula for $\nabla ^2F(z_k)$ in Section 3.1 gives the given approximation for $A_k$.   

Then, it is easy to check that
\begin{eqnarray*}
A_k^{-1}=\frac{1}{d(d-1)r^{d-2}}\left[ \begin{array}{ccc}
		\cos ((d-2)\theta )&-\sin ((d-2)\theta ) \\
		-\sin ((d-2)\theta )&-\cos ((d-2)\theta )\\
	\end{array} \right]+O(\frac{1}{r^{d-3}}),
\end{eqnarray*}
and
\begin{eqnarray*}
A_k^{-1}.\nabla F(z_k)=\frac{r}{d-1}\left[ \begin{array}{ccc}
		\cos ((d-2)\theta /2) \\
		-\sin ((d-1)\theta /2)\\
	\end{array} +O(r^2)\right]=\frac{z_k}{d-1}+O(r^2). 
\end{eqnarray*}
Since $v_k=A_k^{-1}.\nabla F(z_k)$, we obtain the asserted formula for $v_k$. ({\bf Remark}: Note that $z_k-v_k$ is - up to an $O(r^2)$ error - the update rule for Newton's method for Optimization applied to $F(x,y)$. Hence, since by the above formulas we have $$z_k-v_k=\frac{(d-2)z_k}{d-1}+O(|z_k|^2),$$ when $d\geq 2$ Newton's method for Optimization, applied to $F(x,y)$, will converge to the critical point $z^*=0$ if the initial point $z_0$ is close enough to $z^*$. This fact was observed experimentally in \cite{RefFHTW2}.) 

We observe that $||U_1||=1+O(r)$, $||U_2||=1+O(r)$, and since $U_1$ and $U_2$ correspond to different eigenvalues of a real symmetric matrix, we have $<U_1,U_2>=0$. Therefore, from the definition of $w_k$, we obtain
\begin{eqnarray*}
w_k&=&<A_k^{-1}.\nabla F(z_k),U_1>U_1-<A_k^{-1}.\nabla F(z_k),U_2>U_2+O(r^2)\\
&=&\frac{r}{d-1}\left[ \begin{array}{ccc}
		\cos ((d-1)\theta ) \\
		-\sin ((d-1)\theta )\\
	\end{array} \right]+O(r^2).
\end{eqnarray*}

By definition, 
\begin{eqnarray*}
\widehat{w_k}=\frac{w_k}{\max\{(1,\theta ||w_k||\}}. 
\end{eqnarray*}
Since $||w_k||=O(r)$ by the above paragraph, we obtain $\widehat{w_k}=w_k$, as wanted. 

By definition of BNQN, $z_{k+1}=z_k-\gamma _k\widehat{w_k}$, where $\gamma_k$ is the value obtained after finishing the second While loop in Algorithm \ref{table:alg0}. From the above paragraphs, to show that $z_{k+1}=z_k-w_k$, it suffices to show that $\gamma _k=1$. ( Note: when $d=2$, the point $z^*$ is a saddle point of $F$, and this assertion is proven in \cite{RefT} in the setting of $C^3$ functions in higher dimensions.) This amounts to showing the following: 

{\bf Claim:} For $r>0$ small enough, we have
$$F(z_k-w_k<F(z_k)-\frac{1}{3}<w_k,\nabla F(z_k)>+O(r^{d+1}).$$

Proof of Claim: 

The following formula for $F(z)$ is more convenient for the purpose here: 
$$F(z)=\frac{1}{2}+\mathcal{R}(z^d)+O(|z|^{d+1}),$$
where $\mathcal{R}(.)$ is the real part of a complex number. 

Observe that 
\begin{eqnarray*}
w_k&=&\frac{r}{d-1}\left[ \begin{array}{ccc}
		\cos ((d-1)\theta ) \\
		-\sin ((d-1)\theta )\\
	\end{array} \right]+O(r^2)\\&=&\frac{re^{-i(d-1)\theta }}{d-1}+O(r^2)=z_k(1-\frac{1}{d-1}\frac{r^d}{z_k^d})+O(r^2).
\end{eqnarray*}

Therefore,
$$F(z_k-w_k=\frac{1}{2}+\mathcal{R}(z_k^d(1-\frac{1}{d-1}\frac{r^d}{z_k^d})^d)+O(r^{d+1}),$$
and
$$F(z_k-w_k+O(r^2))-F(z_k)=\mathcal{R}(z_k^d[(1-\frac{1}{d-1}\frac{r^d}{z_k^d})^d-1])+O(r^{d+1}).$$

It is also easy to check that 
$$\frac{1}{3}<w_k,\nabla F(z_k)>=\frac{dr^d}{3(d-1)}+O(r^{d+1}).$$

Hence, it is reduced to showing that
$$B_k:=\mathcal{R}(z_k^d[(1-\frac{1}{d-1}\frac{r^d}{z_k^d})^d-1])+\frac{dr^d}{3(d-1)}+O(r^{d+1})<0.$$

We can divide by $r^d$ and use $z=re^{i\theta}$ to simplify the above expression to 
$$\frac{B_k}{r^d}=\mathcal{R}(e^{id\theta }[(1-\frac{e^{-id\theta}}{d-1})^d-1])+\frac{d}{3(d-1)}+O(r).$$

Defining $u=-\frac{e^{-id\theta}}{d-1}$, then 
\begin{eqnarray*}
e^{id\theta }[(1-\frac{e^{-id\theta}}{d-1})^d-1]&=&\frac{-1}{(d-1)u}[(1+u)^d-1]\\
&=&\frac{-1}{d-1}[d+\frac{d(d-1)}{2}u+\frac{d(d-1)(d-2)}{6}u^2+\ldots ],
\end{eqnarray*}
which yields
$$\mathcal{R}(e^{id\theta }[(1-\frac{e^{-id\theta}}{d-1})^d-1])+\frac{d}{3(d-1)}=\frac{-d}{d-1}\mathcal{R}[\frac{2}{3}+\frac{d-1}{2}u+\frac{(d-1)(d-2)}{6}u^2+\ldots].$$

We now estimate the first 3 terms in the RHS of the above expression:
$$M=\mathcal{R}[\frac{2}{3}+\frac{d-1}{2}u+\frac{(d-1)(d-2)}{6}u^2].$$

We claim that $M>1/10$. Indeed, we can write $u=e^{i\alpha }/(d-1)$ for some $\alpha \in [-3\pi /4,5\pi /4]$. Then 
\begin{eqnarray*}
M&=&\mathcal{R}[\frac{2}{3}+\frac{d-1}{2}u+\frac{(d-1)(d-2)}{6}u^2]\\
&=&\frac{2}{3}+\frac{d-1}{2}\frac{1}{d-1}\cos (\alpha )+ \frac{(d-1)(d-2)}{6}\frac{1}{(d-1)^2}\cos (2\alpha ). 
\end{eqnarray*}

If $-3\pi /4\leq \alpha \leq 3\pi /4$, then $\cos (\alpha)\geq -\sqrt{2}/2$, while $\cos(2\alpha )\geq -1$. Therefore,
\begin{eqnarray*}
M&\geq& \frac{2}{3}-\frac{d-1}{2}\frac{1}{d-1}\frac{\sqrt{2}}{2}-\frac{(d-1)(d-2)}{6}\frac{1}{(d-1)^2}\\
&\geq&\frac{2}{3}-\frac{\sqrt{2}}{4}-\frac{1}{6}\geq \frac{1}{10}. 
\end{eqnarray*}

On the other hand, if $3\pi/4 \leq \alpha \leq 5\pi /4$, then $\cos (\alpha )\geq -1$ and $\cos (2\alpha)\geq 0$. Hence
\begin{eqnarray*}
M&\geq& \frac{2}{3}-\frac{d-1}{2}\frac{1}{d-1}\\
&=&\frac{2}{3}-\frac{1}{2}\geq \frac{1}{10}. 
\end{eqnarray*}

From this, since $|u|=1/(d-1)$, we obtain
\begin{eqnarray*}
&&\mathcal{R}(e^{id\theta }[(1-\frac{e^{-id\theta}}{d-1})^d-1])+\frac{d}{3(d-1)}\\
&=&\frac{-d}{d-1}[M+\frac{(d-1)(d-2)(d-3)}{4!}u^3+\frac{(d-1)(d-2)(d-3)(d-4)}{5!}u^4+\ldots ]\\
&\leq&\frac{-d}{d-1}[M-\sum_{j=4}^{\infty}\frac{1}{j!}]\leq\frac{-d}{d-1}[\frac{1}{10}-(e-1-1-\frac{1}{6})]<0. 
\end{eqnarray*}

This implies that $B_k<0$ as wanted, and the proofs of Claim and also Lemma \ref{LemmaApproximateFormulaDynamics} are finished. 
\end{proof}

To summarize, near the critical point $z^*=0$ of order $d$, the update rule of BNQN is: 
$$z_{k+1}=z_k-\frac{r}{d-1}\left[ \begin{array}{ccc}
		\cos ((d-1)\theta ) \\
		-\sin ((d-1)\theta )\\
	\end{array} \right]+O(r^2),$$
where $z_k=re^{i\theta}$. 

Hence, the dynamical system associated to BNQN is given by:
\begin{equation}
\Phi (z)=\left[ \begin{array}{ccc}
		x \\
		y\\
	\end{array} \right]-\frac{r}{d-1}\left[ \begin{array}{ccc}
		\cos ((d-1)\theta ) \\
		-\sin ((d-1)\theta )\\
	\end{array} \right]+O(r^2),
\label{EquationFormulaBNQN}\end{equation}
where $z=(x,y)=re^{i\theta}$. Its linear part is 
\begin{equation}
\Phi _1(z)=\left[ \begin{array}{ccc}
		x \\
		y\\
	\end{array} \right]-\frac{r}{d-1}\left[ \begin{array}{ccc}
		\cos ((d-1)\theta ) \\
		-\sin ((d-1)\theta )\\
	\end{array} \right].
\label{EquationFormulaBNQNLinearization}\end{equation}

From the formula for $\Phi _1(z)$, it is clear that if $d\geq 3$ then $\Phi (z)$ is continuous  but not $C^1$ near $z^*=0$. Therefore, the usual results for stable/unstable manifolds are not applicable for $\Phi (z)$.

\subsection{The dynamics of the linearization $\Phi _1(z)$} As a warm-up, before considering the complicated dynamical system $\Phi (z)$ in (\ref{EquationFormulaBNQN}), let us study its linearization $\Phi _1(z)$ from (\ref{EquationFormulaBNQNLinearization}). $\Phi _1(z)$ is easier to deal with, while still displays the essence of the dynamics of $\Phi (z)$ and suggests a strategy for treating $\Phi (z)$. 

We obtain the following refinement of Theorem \ref{Theorem2}. The condition $d\geq 3$ is necessary, since for $d=2$ the map $\Phi _1$ does not have an isolated fixed point at $(0,0)$. 
\begin{thm} Let $d\geq 3$, and define ${C_j}$ to be the half-line $\{z\in \mathbf{C}: arg(z)=\frac{\pi j}{d}\}$, for $j=1,\ldots ,2d$. 

1) The map $\Phi _1:\mathbf{C}\rightarrow \mathbf{C}$ defined by (\ref{EquationFormulaBNQNLinearization}) has only 1 fixed point, which is $z^*=0$. 

2) For $j=1,\ldots d$, $C_{2j}$ are stable curves for the dynamics of $\Phi _1(z)$, and $C_{2j+1}$ are unstable curves for the dynamics of Analytic BNQN. 

More precisely, if $z_0\in \mathbf{C}$, and $\{z_n\}$ is the sequence constructed by the map $\Phi _1(z)$ with initial point $z_0$, then: 

2a) If $z_0\notin \bigcup _{j=1}^dC_{2j}$, then $\{z_n\}$ converges to $\infty$. If $z_0\in C_{2j}$, then $z_n\in C_{2j}$ for all $n$ and $\lim _{n\rightarrow \infty}z_n=z^*$. 

2b) The curves $C_{2j+1}$ are invariant by the dynamics of $\Phi _1(z)$, and the restriction of $\Phi _1(z)$ to $C_{2j+1}$ is repelling. 

2c) The curves $C_1,\ldots ,C_{2d}$ partition $\mathbb{C}$ into $2d$ open sets $U_1,\ldots ,U_{2d}$ each being bounded between one stable curve and one unstable curve, and has the following property: 

If $z_0\in U_j$, then  $z_n\in U_j$ for all $n$, comes closer to the unstable curve, and converge to $\infty$. See the figure \ref{fig3}.
\begin{figure}[!htb]
	\centering
\includegraphics[width=0.7\textwidth,height=0.4\textheight]{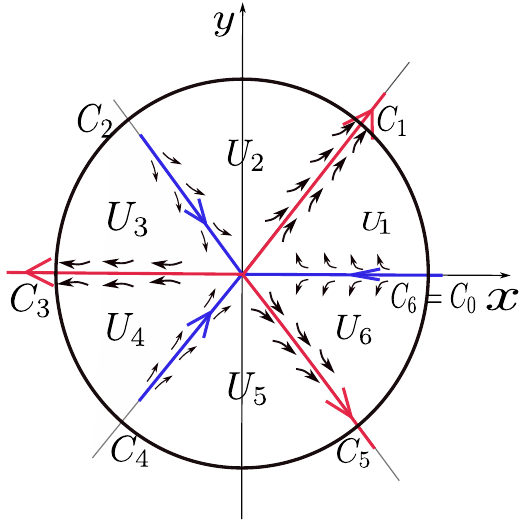}
	\caption{The orbits of the linearization $\Phi _1$ near the critical point for $d=3$.}
	\label{fig3}
\end{figure}
\label{Theorem2Linearization}\end{thm}
\begin{proof}
1) The map $\Phi _1(z)$ can also be written as: 
$$\Phi _1(z)=z-\frac{\overline{z}^{d-1}}{(d-1)r^{d-2}},$$
from which it is ready to see that $\Phi _1(z)$ has only 1 fixed point, which is $z^*=0$. 

2) 

2a) and 2b): 

We first check that the half-lines $C_j$'s are invariant by $\Phi _1$, and the restriction of $\Phi _1$ to each curve is either attracting or repelling. Since the curves are similar, it is enough to check with $C_d$ and $C_{2d}$. 

The curve $C_{2d}$ is precisely the halfline $\{z=(x,0):x\geq 0\}$, where $\theta =2\pi$. In this case, the restriction map $\Phi _1|_{C_{2d}}$ is 
$$\Phi _1|_{C_{2d}}(x)=\frac{d-2}{d-1}x,$$
which is attracting. 

The curve $C_{d}$ is precisely the halfline $\{z=(x,0):x\leq 0\}$ and $\theta =\pi $. Depending on whether $d$ is odd or even, the restriction map $\Phi _1|_{C_{d}}$ has different behaviour. If $d$ is odd, then $\cos ((d-1)\theta )=1$, and  the restriction map $\Phi _1|_{C_{d}}$ is 
$$\Phi _1|_{C_{d}}(x)=\frac{d}{d-1}x,$$
which is repelling. On the other hand, if $d$ is even, then $\cos ((d-1)\theta )=-1$, and  the restriction map $\Phi _1|_{C_{d}}$ is
$$\Phi _1|_{C_{d}}(x)=\frac{d-2}{d-1}x,$$
which is attracting. 

This finishes the proofs of 2a and 2b, except the assertion that if $z_0\notin \bigcup _{j=1}^dC_{2j}$, then the constructed sequence $\{z_n\}$ converges to $\infty$. This will be done in the proof of part 2c, which we now show. 

2c) For $j=1,\ldots ,2d-1$, we define $U_j$ to be the open set bound between the 2 curves $C_j$ and $C_{j+1}$. We also define $U_{2d}$ to be the open set bound between the 2 curves $C_{2d}$ and $C_1$. 

By symmetry, it suffices to consider the sector $U_{2d}$. It is the sector $\{z\in \mathbf{C}: z\not= 0, 0<arg(z)<\pi /d\}$. We need to show that: if $z_0=re^{i\theta }\in U_{2d}$, then $z_n\in U_{2d}$ for all $n$, and $\lim _{n\rightarrow\infty}|z_n|=\infty$.

To show that $z_n\in U_{2d}$ for all $n$, by induction it is enough to show that $z_1\in U_{2d}$. Writing $z_1=(x_1,y_1)=r_1e^{i\theta _1}$, it amounts to showing that $x_1>0$ and $0<\tan (\theta _1)<\tan (\pi /d)$. 

To see that $x_1>0$, we observe that since $0<\theta <\pi /d \leq \pi /2$, we have $0<\theta <(d-1)\theta <(d-1)\pi /d$. Since the cosine function is strictly decreasing on $[0,\pi /2]$ and non-negative there, and is negative on $(\pi /2,\pi ]$ it follows that $\cos (\theta )>\cos ((d-1)\theta )$. Therefore, 
\begin{eqnarray*}
x_1=r\cos (\theta )-\frac{r}{d-1}\cos ((d-1)\theta )>0. 
\end{eqnarray*}

We next check that $0<\tan (\theta _1)<\tan (\pi /d)$. This is equivalent to
\begin{eqnarray*}
0<\frac{y+r\sin ((d-1)\theta )/(d-1)}{x-r\cos ((d-1)\theta )/(d-1)}<\tan(\pi/d),
\end{eqnarray*}
since
\begin{eqnarray*}
\tan(\theta_1)&=&\frac{y_1}{x_1},\\
x_1&=&x-r\cos ((d-1)\theta )/(d-1),\\
y_1&=&y+r\sin((d-1)\theta )/(d-1).
\end{eqnarray*}

Since $z_0=(x_0,y_0)\in U_{2d}$, it follows that both $x_0,y_0>0$ and $\sin ((d-1)\theta )>0$. Hence $y_1>0$, one inequality $0<\tan (\theta _1)$ easily follows. For the other inequality $\tan (\theta _1)<\tan (\pi /d)$, we observe that it is equivalent to (after dividing by $r$):
$$\frac{\sin (\theta )+\sin ((d-1)\theta )/(d-1)}{\cos (\theta )-\cos ((d-1)\theta )/(d-1)}<\frac{\sin (\pi/d)}{\cos(\pi/d)}.$$
The denominators of the fractions on both sides are positive, hence after multiplying  them and simplifying (by using the usual trigonometric identities for sine and cosine of a sum) we obtain an equivalent inequality: 
$$v(\theta):=\sin(\theta-\frac{\pi}{d})+\frac{1}{d-1}\sin ((d-1)\theta +\frac{\pi}{d})<0,$$
for $0<\theta <\pi /d$. We can see that
\begin{eqnarray*}
v'(\theta )&=&\cos(\theta-\frac{\pi}{d})+\cos((d-1)\theta +\frac{\pi}{d})\\
&=&2\cos(\frac{d\theta }{2})\cos(\frac{(d-2)\theta -2\pi/d}{2}).
\end{eqnarray*}
The latter is easily seen to be positive for $\theta \in (0,\pi /d)$. Hence, the function $v(\theta )$ is strictly increasing on $(0,\pi /d)$, and 
$$v(\theta )<v(\pi /d)=0.$$

Therefore, we showed that if $z_0\in U_{2d}$, then $z_n\in U_{2d}$ for all $n$. 

To show that $\{z_n\}$ comes closer to the unstable curve $C_1=\{z\in \mathbf{C}: arg(z)=\pi /d\}$, it suffices to show that $\theta _n$ strictly increases, where $z_n=r_ne^{i\theta _n}$. Again, it suffices, by induction, to show that $\theta _1>\theta $ (see the figure \ref{fig1}).
\begin{figure}[!htb]
	\centering
\includegraphics[width=0.6\textwidth,height=0.35\textheight]{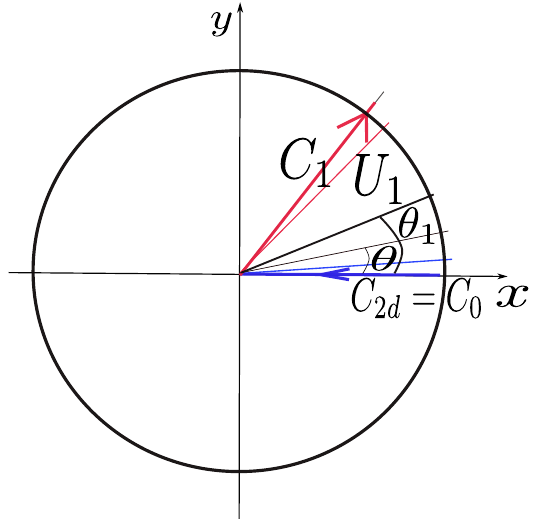}
	\caption{The angle in $U_1$, the sector bound by the stable ray $C_{2d}$ and the unstable ray $C_1$, is expanding under $\Phi _1$.}
	\label{fig1}
\end{figure}
We already knew that $0<\theta ,\theta _1<\pi /d$, hence $\theta _1>\theta $ is equivalent to that $\tan (\theta )<\tan (\theta _1)$. The latter is equivalent to that $y_0/x_0 < y_1/x_1$, which is:
$$\frac{\sin (\theta)}{\cos(\theta )}<\frac{\sin (\theta )+\sin ((d-1)\theta )/(d-1)}{\cos (\theta )-\cos ((d-1)\theta )/(d-1)},$$
which after multiplying with the product of the denominators, a positive number, and simplifying, becomes: $\frac{\sin (d\theta )}{d-1}>0$. The latter is true, since $0<d \theta <\pi$.  

$$y_{n+1}=y_n(1+\frac{\sin ((d-1)\theta _n)}{(d-1)\sin (\theta _n)}).$$
From the above arguments, we know that $\theta _n<\theta _{n+1}<\pi /d$ for all $n$. Therefore, $\theta _n$ converges to some $0<\theta ^*\leq \pi /d$. Then $(d-1)\theta _n$ converges to $(d-1)\theta ^*$, which lies between $0$ and $\pi$. Hence, $y_{n+1}>y_n>0$ for all $n$. Moreover, 
$$\lim _{n\rightarrow\infty}\frac{y_{n+1}}{y_n}=1+\frac{\sin ((d-1)\theta ^*)}{(d-1)\sin (\theta ^*)}>1.$$
This shows that $\lim _{n\rightarrow\infty}y_{n}=+\infty$, as wanted. The proof of Theorem \ref{Theorem2Linearization} is completed. 

\end{proof}

\subsection{Proof of Theorem \ref{Theorem2}} First, we will treat the case $d\geq3$, based on the insights from Theorem \ref{Theorem2Linearization}. Afterwards, we will treat the case $d=2$, which is fundamentally different. 

\subsubsection{The case $d\geq 3$} Here we prove Theorem \ref{Theorem2} for the case $d\geq 3$. 

Here is a short comment about the main strategy. We note that since both $\Phi$, $\Phi _1$ and $\Phi -\Phi _1$ are generally not $C^1$ functions near $z^*=0$, it may be difficult (or even impossible) to use a perturbation argument or to use a conjugate by a local invertible real analytic map to directly deduce Theorem \ref{Theorem2} from Theorem \ref{Theorem2Linearization}. Here we will instead work directly with $\Phi$, but with the goal of establishing similar properties like those achieved in the proof of Theorem \ref{Theorem2Linearization}. Since
the stable/unstable curves for $\Phi _1$ are precisely the halflines $\{z\in \mathbf{C}: arg(z)=j\pi /d\}$ for $j=1,\ldots ,2d$, our idea is to work in a small sector around these halflines to find the corresponding stable/unstable curves for $\Phi$. 

We will first work with a small sector around the halflines $\{z\in \mathbf{C}: arg(z)=0\}$ and $\{z\in \mathbf{C}: arg(z)=\pi/d\}$, to find the stable curve $C_{2d}$ and the unstable curve $C_1$ in Theorem \ref{Theorem2}. We will show afterwards how to work with other halflines. 

Therefore, we concentrate now on 2 small sectors $S_{r_0,\epsilon _0,2d}:=\{z\in \mathbf{C}: 0<|z|<r_0, |arg(z)|<\epsilon _0\}$ and $S_{r_0,\epsilon _0,1}:=\{z\in \mathbf{C}: 0<|z|<r_0, |arg(z)-\pi /d|<\epsilon _0\}$, where $r_0,\epsilon _0>0$ are small enough constants to be determined later. The proof consists of several steps. 

{\bf Step 1: Show that $\Phi $ is real analytic in $ S_{r_0,\epsilon _0,2d}$.} Since $f$ is analytic in the disk $\mathbf{D}(0,r_0)$ and $\tau \in 2\mathbf{N}$, it follows that the matrix $A=\nabla ^2F(z)+\delta _0||\nabla F(z)||^{\tau}Id$ (used in Algorithm \ref{table:alg0}) is real analytic for $z\in \mathbf{D}(0,r_0)$. In the sector $S_{r_0,\epsilon _0,2d}$, we have $x>0$ and $|y|/|x|$ is small,  hence $r=\sqrt{x^2+y^2}=x\sqrt{1+(y/x)^2}$ is real analytic and nowhere 0. 

Then the matrix $B(z)=A(z)/r^{d-2}$ is also real analytic in the same sector. Moreover, by Lemmas \ref{LemmaEigenHessian} and \ref{LemmaApproximateFormulaDynamics}, the  2 eigenvalues of $B$ are $\pm d(d-1)+O(r)$. By the integral formula for projections on eigenspaces of a linear map in the theory of resolvents \cite{kato}, as used in \cite{RefTT}, we can describe $w_k$ in Algorithm \ref{table:alg0}), as a function in the variable $z_k$, as follows:
$$w_k=\frac{1}{r^{d-2}}[-\int _{L_1}(B(z_k)-\zeta Id )^{-1}.\nabla F(z_k)d\zeta +\int _{L_2}(B(z_k)-\zeta Id )^{-1}.\nabla F(z_k)d\zeta],$$
where $L_1$ is a small circle in the complex plane with center $d(d-1)$, and $L_2$ is a small circle in the complex plane with center $-d(d-1)$. Refer to the proof of Lemma \ref{LemmaApproximateFormulaDynamics}, in fact: 
\begin{eqnarray*}
-\frac{1}{r^{d-2}}\int _{L_1}(B(z_k)-\zeta Id )^{-1}.\nabla F(z_k)d\zeta&=&\frac{1}{||U_1||^2}<\nabla F(z_k),U_1>U_1,\\
-\frac{1}{r^{d-2}}\int _{L_2}(B(z_k)-\zeta Id )^{-1}.\nabla F(z_k)d\zeta&=&\frac{1}{||U_2||^2}<\nabla F(z_k),U_2>U_2.
\end{eqnarray*}

The integrands in the above integrals are real analytic, which can be seen easily by using Taylor's series:
$$(B(z_k)-\zeta Id )^{-1}.\nabla F(z_k)=-\frac{1}{\zeta}\sum _{j=0}^{\infty} \frac{B(z_k)^j}{\zeta ^j},$$
which is a uniform convergence of real analytic sequences defined on the same domain. Since nonsingular integrations involving real analytic maps are real analytic, we conclude that $w_k$ depends real analytically in the variable $z_k$, in the sector $S_{r_0,\epsilon _0,2d}$.  

Therefore, since $|w_k|=O(r)$,  also $1/(1+|w_k|^2)^{1/2}$ is real analytic in the same sector. Thus, the same holds for $\widehat{w_k}=w_k$, and $z_{k+1}=z_k-w_k$, which is the update rule of BNQN near $z^*$, by virtue of Lemma \ref{LemmaApproximateFormulaDynamics}. 

In conclusion, $\Phi $ is real analytic in $ S_{r_0,\epsilon _0,2d}$. 

{\bf Step 2: Conjugate $\Phi|_{S_{r_0,\epsilon _0,2d}}$ to a real analytic map $\Psi _{2d}$ on a rectangle which has  $(0,0)$ as a fixed point.} 

Even though $\Phi $ is real analytic in $ S_{r_0,\epsilon _0,2d}$, it has no fixed point in that sector, and hence we cannot apply stable/unstable manifolds theorem. We will now transform it to another map for which this is applicable. 

\begin{figure}[!htb]
	\centering
\includegraphics[width=1\textwidth,height=0.35\textheight]{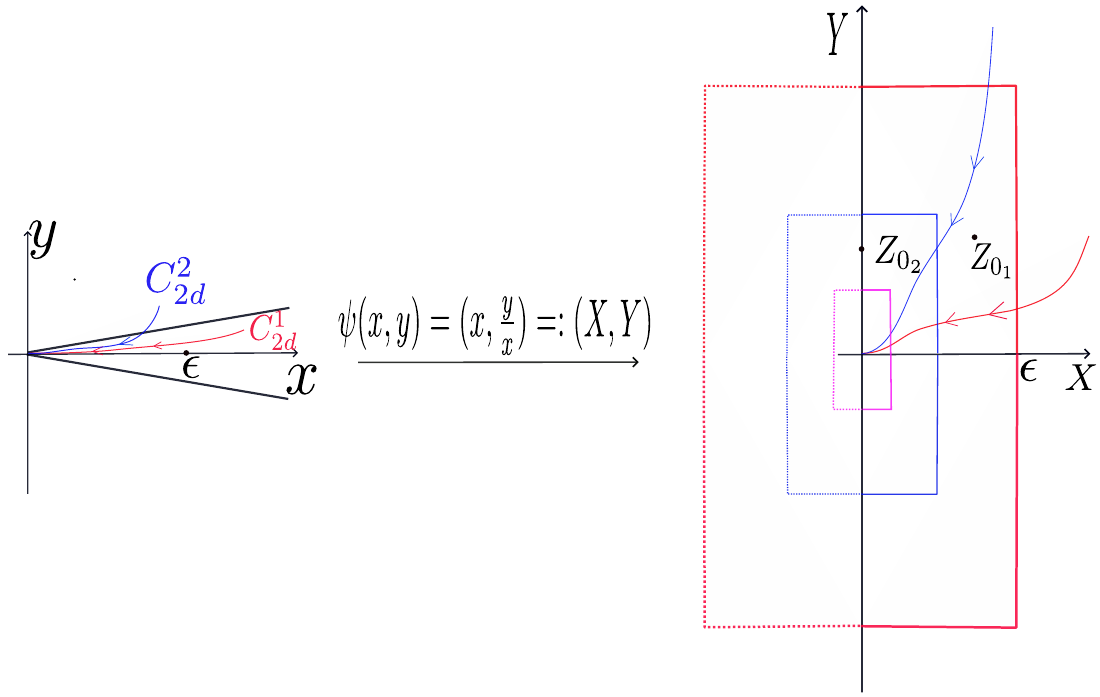}
	\caption{Transforming the small sector into a half of a rectangle containing $(0,0)$ inside. The curves $C_{2d}^1$, $C_{2d}^2$ in the sector are pullback of certain curves $Z_{0_1}$, $Z_{0_2}$ in the rectangle. To make sure that $C_{2d}^2$  only exits the   sector through the boundary on the right side (i.e. an arc of the circle), we need to shrink both the sector and the rectangle (from the big red rectangle to the smaller blue rectangle).}
	\label{fig4}
  \end{figure}

We consider the map $\psi :z=(x,y)\rightarrow w=(X,Y)=(x,\frac{y}{x})$, which is an invertible real analytic map from the set $\{x\not= 0\}$ onto its image (see the figure \ref{fig4})
Therefore, the map $\Phi \circ \psi ^{-1}$ is real analytic on the image $\psi (S_{r_0,\epsilon _0,2d})$. We now show that $\Phi \circ \psi ^{-1}$ can be extended real analytically to the whole rectangle $R_{r_0,\epsilon _1}$. We note that in $R_{r_0,\epsilon _1}$, both $r$ and $\theta $ are real analytic functions of $x,y$ (for example, they are the real and imaginary parts of the logarithm $\log (z)$ which is analytic in the domain). Moreover, $\tan(\theta )=Y$ shows that the coordinates $\theta$ and $Y$ are real analytic inverses of each other. The function $r=\sqrt{x^2+y^2}=X\sqrt{1+Y^2}$, which is originally defined for $0<X<r_0$, $|Y|<\epsilon _1$ can also be extended real analytically to $|X|<r_0$,  $|Y|<\epsilon _1$. Since the function $\Phi \circ \psi ^{-1}$ is real analytic in $r,\theta $, it can be extended real analytically to the rectangle $R_{r_0,\epsilon _1}$, as claimed. (For the reader's convenience, we provide a detailed proof as follows. The possibility to extend  $\Phi$ real analytically is the same as that $w_k$ in Step 1 can be extended real analytically to the whole rectangle. To this end, it suffices to show that for the matrix $B(z)$ in Step 1, after the change of coordinates $x=X$, $y=XY$, the matrix $B(X,XY)$ is extended real analytic to the domain $X\leq 0$. Note that we have then matrix $A(z)$ is real analytic in a neighbourhood of $(0,0)$, hence has a convergent power series $A(x,y)=\sum _{n_1,n_2\geq 0}a_{n_1,n_2}x^{n_1}y^{n_2}$. By Lemma \ref{LemmaApproximateFormulaDynamics}, we have that for all $a_{n_1,n_2}\not=0$, we must have $n_1+n_2\geq d-2$. Therefore, when substituting $x=X$, $y=XY$, the new series, denoted by $\widehat{A}(X,Y)$ is divisible by $X^{d-2}$ and has the same radius of convergence as before. Therefore, the new series for $B(x,y)$, denoted by $\widehat{B}(X,Y)=\widehat{A}(X,Y)/(X^{d-2}\sqrt{(1+Y^2)^{d-2}})$ is analytically extended to $X\leq 0$. The eigenvalues of $\widehat{B}(X,Y)$ on the rectangle are close to $\pm d(d-1)$, and hence using integral formula as in Step 1 shows that $w_k(X,XY)$, originally defined for $X>0$, can be analytically extended to the domain $X\leq 0$.)

Since $X=r/\sqrt{1+Y^2}$, it follows that $X$ is also a real analytic function in $r$ and $\theta$ (in the given domain), and vice versa $r$ is a real analytic function in $X$ and $\theta$. Then we can write 
$$\Phi (X,XY)=\left[ \begin{array}{ccc}
		X \\
		XY\\
	\end{array} \right]-\frac{X\sqrt{1+Y^2}}{d-1}\left[ \begin{array}{ccc}
		\cos ((d-1)\theta ) \\
		-\sin ((d-1)\theta )\\
	\end{array} \right]+\left[ \begin{array}{ccc}
		G_1(X,\theta ) \\
		G_2(X,\theta )\\
	\end{array} \right],$$
where $G_1,G_2$ are real analytic in $X,\theta$ and divisible by $X^2$. 

First consider $X\not= 0$.  As in the proof of Theorem \ref{Theorem2Linearization}, one can check that then the x-coordiate of $\Phi \circ \psi ^{-1}(X,Y)$ is non-zero. Therefore, we can apply $\psi $ to $(V,W)=\Phi \circ \psi ^{-1}(X,Y)$, and obtain: 
$$\psi \circ \Phi \circ \psi ^{-1}(X,Y)=(V,W/V)=(V,W_2/W_1),$$ 
where
\begin{eqnarray*}
V&=&X-\frac{X\sqrt{1+Y^2}}{d-1}\cos ((d-1)\theta ) +G_1(X,\theta ),\\
W_2&=&Y-\frac{\sqrt{1+Y^2}}{d-1}\sin((d-1)\theta ) + \frac{G_2(X,\theta )}{X},\\
W_1&=&1-\frac{\sqrt{1+Y^2}}{d-1}\cos((d-1)\theta )+\frac{G_1(X,\theta )}{X}.
\end{eqnarray*}
 From the properties of $G_1$ and $G_2$, both $V$, $W_1$ and $W_2$ are real analytic functions (also for $X=0$), and moreover $W_1$ is never 0 in the given domain. Therefore, $\psi \circ \Phi \circ \psi ^{-1}(X,Y)$ (originally defined on $\psi (S_{r_0,\epsilon _0,2d})$) can be extended to a real analytic map in the rectangle $R_{r_0,\epsilon _1}$, as claimed.

{\bf Step 3: Apply Theorems \ref{TheoremStableUnstableManifoldRealAnalytic} and \ref{TheoremNoJumpLemma} to $\Psi _{2d}$.} Denote by $H:R_{r_0,\epsilon _1}\rightarrow \mathbf{R}^2$ the map defined by $H(X,Y)=\psi \circ \Phi \circ \psi ^{-1}(X,Y)$. From Step 2, we know that $H(X,Y)$ is real analytic. It is also easy to check that $(0,0)$ is a fixed point of $H$. We now calculate the Jacobian of $H$ at $(0,0)$ and determine its eigenvalues and eigenvectors. 

From the formulas for $V,W_2,W_1$ and the properties of $G_1,G_2$, it can be checked that (for example, by using that $Y$ is approximately $\theta$ and also $\sin (\theta)$) near $(0,0)$, we have(recall that here we are treating the case $d\geq 3$): 
\begin{eqnarray*}
H(X,Y)=\left[ \begin{array}{ccc}
		\frac{d-2}{d-1}X+O(X^2+Y^2) \\
		\frac{2(d-1)}{d-2}Y+aX+O(X^2+Y^2)\\
	\end{array} \right],
\end{eqnarray*} 
for some real number $a$.

From this, it is easy to check that the Jacobian of $H$ at $(0,0)$ is: 
\begin{eqnarray*}
JH(0,0)=\left[ \begin{array}{ccc}
		\frac{d-2}{d-1}&0 \\
		a&\frac{2(d-1)}{d-2}\\
	\end{array} \right].
\end{eqnarray*}
 
Hence, the two eigenvalues are 
$$0<\frac{d-2}{d-1}<1<\frac{2(d-1)}{d-2},$$
with the corresponding eigenvectors $(1,0)$ and $(0,1)$. Therefore, by shrinking $r_0$ and $\epsilon _1$ (also $\epsilon _0$ so that still we have $\psi (S_{r_0,\epsilon _0,2d})\subset R_{r_0,\epsilon _1}$), we can assume that $(0,0)$ is the only fixed point of $H$ inside $R_{r_0,\epsilon _1}$ and both Theorems \ref{TheoremStableUnstableManifoldRealAnalytic} and \ref{TheoremNoJumpLemma} apply. (In this case, instead of Theorem \ref{TheoremStableUnstableManifoldRealAnalytic}, we can do as follows. First, complexify and use \cite[Theorem 6.4.1]{RefMNTU} to obtain complex stable and unstable manifolds, then use \cite{shub} to deduce that $C^1$ stable and unstable manifolds for the original map exist, therefore the intersection of the complex stable and unstable manifolds to the real plane are smooth real analytic curves.)

{\bf Step 4: Construct the stable curve $C_{2d}$ in $S_{r_0,\epsilon _0,2d}$}: 

The eigenvector corresponds to the eigenvalue $\lambda _1=\frac{d-2}{d-1}$ of $JH(0,0)$ is $e_1=(1,\frac{a}{\lambda _1-\lambda _2})$. Hence,  the stable manifold of $H$ can be defined by an equation $Y=G(X)$, where $G(0)=0$ and $G'(0)=\frac{a}{\lambda _1-\lambda _2}$. Therefore, the pullback of this stable curve, which a priori is defined in $\{0<|x|<r_0, |y_0|<\epsilon _0\}$ by the equation $y=xG(x)$, can be extended as a closed smooth analytic curve $\widehat{C_{2d}}$ in the whole $\mathbb{D}(0,r_0)$ (since $G$ is defined as a real analytic function on $\{|x|<r_0\}$). Since $G(0)=0$, it is then easy to check that $\widehat{C_{2d}}$ is tangent to the x-axis at $(0,0)$. We define $C_{2d}=\widehat{C_{2d}}\cap S_{r_0,\epsilon _0,2d}$.     

{\bf Step 5: Construct the unstable curve $C_1$, and exclude stable curves in $S_{r_0,\epsilon _0,1}$}: 

We now construct the unstable curve $C_1$ in the sector $S_{r_0,\epsilon _0,1}$. The main idea is to use a rotation to the setting considered in Steps 1--4. Let $z'$ be contained in the sector $S_{r_0,\epsilon _0,1}$. Consider the rotation $\phi (z)=e^{i\pi /d}z$. Then, there is $z\in S_{r_0,\epsilon _0,2d}$ so that $z'=\phi (z)$. Then $\Phi ^{\circ n}(z')=\phi \circ \widehat{\Phi}^{\circ n}(z)$, where $\widehat{\Phi }(z)=\varphi ^{-1}\circ \Phi \circ \phi (z)$ and ${\Phi}^{\circ n}$ is the n-th iteration of $\Phi$ (and similar for $\widehat{\Phi}^{\circ n}$). The linear part $\widehat{\Phi}_1$ of $\widehat{\Phi}$ can be readily computed. Since $\widehat{\Phi}_1(z)=\varphi ^{-1}\circ \Phi _1\circ \phi (z)$, where $\Phi _1(z)$ - the linear part of $\Phi (z)$ - has the formula: 
$$\Phi _1(z)=z-\frac{\overline{z}^{d-1}}{(d-1)r^{d-2}},$$
 we have
\begin{eqnarray*}
\widehat{\Phi}_1(z)&=&e^{-i\pi /d}\Phi _1(e^{i\pi /d}z)\\
&=&e^{-i\pi /d}[e^{i\pi /d}z-\frac{(\overline{e^{i\pi/d}z})^{d-1}}{(d-1)r^{d-2}}]\\
&=&e^{-i\pi /d}[e^{i\pi /d}z-\frac{(e^{-i\pi/d}\overline{z})^{d-1}}{(d-1)r^{d-2}}]\\
&=&z+\frac{\overline{z}^{d-1}}{(d-1)r^{d-2}}.
\end{eqnarray*}
 Therefore, in the polar coordinates $z=re^{i\theta }$, the formula for $\widehat{\Phi}$ is:
 $$\widehat{\Phi}(z)=\left[ \begin{array}{ccc}
		r\cos (\theta )+\frac{r}{d-1}\cos ((d-1)\theta )+O(r^2) \\
		r\sin (\theta )-\frac{r}{d-1}\sin ((d-1)\theta )+O(r^2)\\
	\end{array} \right].$$

Therefore, as in Step 2, there is a real analytic map $\widehat{H}:R_{r_0,\epsilon _1}\rightarrow \mathbf{R}^2$, which on $\psi (S_{r_0,\epsilon _0,2d})$ conjugates to $\widehat{\Phi}$. As in Step 3, we find
\begin{eqnarray*}
J\widehat{H}(0,0)=\left[ \begin{array}{ccc}
		\frac{d}{d-1}&0 \\
		b&0\\
	\end{array} \right].
\end{eqnarray*}
This has eigenvalues $0=\lambda _1<1<\frac{d}{d-1}=\lambda _2$, with eigenvectors $(1,0)$ and $(0,1)$.  hence Theorem \ref{TheoremStableUnstableManifoldRealAnalytic} provides an unstable manifold which is real analytic. (In this case, instead of Theorem \ref{TheoremStableUnstableManifoldRealAnalytic}, we can do as follows. First, complexify and use \cite[Theorem 6.4.3 case ii]{RefMNTU} to obtain the complex unstable manifold, then use \cite{shub} to deduce that $C^1$ unstable manifold for the original map exists, therefore the intersection of the unstable manifold to the real plane is a  smooth real analytic curve.)  Again, as in Step 3, the pullback of the unstable curve of the map $\widehat{H}$ can be extended to a closed smooth real analytic curve in $\mathbb{D}(0,r_0)$ which is tangent to the x-axis at $(0,0)$. Pullback further by the rotation $\phi$, we obtain the real analytic curves $\widehat{C_1}$ and $C_1$ in Theorem \ref{Theorem2}. 

Potentially, the pullback under the map $\psi \circ \phi$ of the stable curves of $\widehat{H}$ may be non-empty in $S_{r_0,\epsilon _0,1}$. This will make the dynamics of $\Phi $ more complicated (in particular, in relation to part 3c of Theorem \ref{Theorem2}). We now exclude this possibility, i.e. showing that there is no stable curve inside the sector $S_{r_0,\epsilon _0,1}$ when $r_0>0$ is small enough and $0<\epsilon _0\leq \pi /(2d)$ (independent of $r_0$). Indeed, we will show the stronger claim that $\Phi$ is repelling in the whole sector $S_{r_0,\epsilon _0,1}$.

To this end, from the formula for $\Phi (x,y)$ we find that: 
\begin{eqnarray*}
|\Phi (x,y)|^2&=&(x-\frac{r}{d-1}\cos ((d-1)\theta )+O(r^2))^2\\
&&+(y+\frac{r}{d-1}\sin ((d-1)\theta )+O(r^2))^2\\
&=&(x^2+y^2)+\frac{r^2}{(d-1)^2}-\frac{2r^2}{d-1}\cos (d\theta )+O(r^3).
\end{eqnarray*}
Since $z\in S_{r_0,\epsilon _0,1}$, we have $-\epsilon _0 <\theta -\pi/d <\epsilon _0$, and hence $\pi -d\epsilon _0 <d\theta <\pi +d\epsilon _0$. Under the assumption $0<\epsilon _0\leq \pi /(2d)$ we have $\pi /2 < \theta <3\pi /2 $. Therefore $\cos (d\theta )\leq 0$, and
$$|\Phi (x,y)|^2\geq r^2+\frac{r^2}{(d-1)^2}+O(r^3).$$
Therefore if $r_0>0$ is small enough then $|\Phi (z)|>\frac{(2d-1)}{(2d-2)}|z|$ for all $z\in S_{r_0,\epsilon _0,1}$, and is repelling. In particular, there is no stable curve inside $S_{r_0,\epsilon _0,1}$.

{\bf Step 6: Construct other stable and unstable curves}: 

By using rotations we can similarly construct all other stable and unstable curves. Note that if the angle of the rotation is an even multiple of $\pi /d$, then the linearization will be the same as that of $\Phi _1$, while if the angle is an odd multiple of  $\pi /d$, then the linearization will be the same as that of $\widehat{\Phi }_1$. Therefore, the even curves $C_2,\ldots C_{2d}$ are all stable curves and the local behaviour of $\Phi $ near these curves are the same, while the odd curves $C_1,\ldots C_{2d-1}$ are all unstable curves and the local behaviour of $\Phi$ near these curves are the same. 

{\bf Step 7: Describe the dynamics of $\Phi$ in the domain $U_{d}$ bound by the 2 consecutive stable curves $C_{2d}$ and $C_2$ and contains the unstable curve $C_1$}: 

We now describe the behaviour of the dynamics of $\Phi$ in the domain $U_{d}$ bound by the 2 stable curves $C_{2d}$ and $C_2$, and contains the unstable curve $C_1$ inside. We first consider the half of this domain bound by the stable curve $C_{2d}$ and the unstable curve $C_1$.  The behaviour of $\Phi $ on the other half of $U_{d}$ will be similar and can be seen by reflecting about the halfline $\{arg(z)=\frac{\pi}{d}\}$.  

Note first that since $C_1$ and  $C_{2d}$ are smooth and have different tangent directions at $(0,0)$, by shrinking $r_0$ we will have that $C_1\cap C_{2d}=(0,0)$. 

We now show the following generalization of an argument in the proof of Theorem \ref{Theorem2Linearization}. 
\begin{lem} For $z=(x,y)= re^{i\theta }$, we define $z_1=\Phi (z)=r_1e^{i\theta _1}$. Let $\sigma =3.5$. If $0<\epsilon _0\leq \theta  \leq \frac{\pi}{d}-\epsilon _0$, $0<r\leq r_0$, $\epsilon _0=r_0^{1/4}$ and $r_0>0$ is small enough, then 
\begin{eqnarray*}
\theta <\theta _1<\frac{\pi}{d}-\epsilon _0^{\sigma}.
\end{eqnarray*}
\label{LemmaIncreasedAnggle}\end{lem}
\begin{proof}
Adapting the proof of Theorem \ref{Theorem2}, we see that this is equivalent to 
\begin{eqnarray*}
\tan (\theta ) <\tan( \theta _1)<\tan (\frac{\pi}{d}-\epsilon _0^{\sigma}), 
\end{eqnarray*}
and is the same as the following two inequalities
\begin{eqnarray*}
\frac{\sin (\theta )}{\cos (\theta )}&<& \frac{\sin (\theta )+\frac{1}{d-1}\sin ((d-1)\theta )+O(r)}{\cos (\theta )-\frac{1}{d-1}\cos ((d-1)\theta )+O(r)},\\
\frac{\sin (\frac{\pi}{d}-\epsilon _0^{\sigma})}{\cos (\frac{\pi}{d}-\epsilon _0^{\sigma})}&>& \frac{\sin (\theta )+\frac{1}{d-1}\sin ((d-1)\theta )+O(r)}{\cos (\theta )-\frac{1}{d-1}\cos ((d-1)\theta )+O(r)}. 
\end{eqnarray*}
This is reduced to two simpler inequalities: 
\begin{eqnarray*}
\frac{\sin (d\theta )}{d-1}+O(r)&>0&,\\
\sin (\frac{\pi}{d}-\epsilon _0^{\sigma}-\theta )-\frac{1}{d-1}\sin(\frac{\pi}{d}-\epsilon _0^{\sigma}+(d-1)\theta )+O(r)&>&0. 
\end{eqnarray*}
Ignoring the $O(r)$ terms for a moment, we can see that the two functions on the LHS are strictly monotone for $\theta$ in $[\epsilon _0, \frac{\pi}{d}-\epsilon _0]$. Therefore, the two needed inequalities are equivalent to the following 4 inequalities (by plugging in $\theta =$ $\epsilon _0$ and $\pi /d$): 
\begin{eqnarray*}
\frac{\sin (d\epsilon _0)}{d-1}+O(r)&>&0,\\ 
\frac{\sin (\pi - d\epsilon _0)}{d-1}+O(r)&>&0,\\
\sin (\frac{\pi}{d}-\epsilon _0^{\sigma}-\epsilon _0)-\frac{1}{d-1}\sin(\frac{\pi}{d}-\epsilon _0^{\sigma}+(d-1)\epsilon _0 )+O(r)&>&0,\\
\sin (-\epsilon _0^{\sigma}+\epsilon _0)-\frac{1}{d-1}\sin(\pi-\epsilon _0^{\sigma}-(d-1)\epsilon _0 )+O(r)&>&0. 
\end{eqnarray*}
The first 3 inequalities can be easily checked to hold  when $\epsilon _0=r_0^{1/4}$ and $r_0$ is small enough. Using the formula $\sin (\pi -t)=\sin (t)$ and $\sin (t)=t-t^3/6+O(t^5)$ when $|t|$ is small, the function on the LHS of the 4-th inequality above can be written as follows: 
\begin{eqnarray*}
&&\sin (-\epsilon _0^{\sigma}+\epsilon _0)-\frac{1}{d-1}\sin(\pi-\epsilon _0^{\sigma}-(d-1)\epsilon _0 )+O(r)\\
&=&\sin (-\epsilon _0^{\sigma}+\epsilon _0)-\frac{1}{d-1}\sin(\epsilon _0^{\sigma}+(d-1)\epsilon _0 )+O(r)\\
&=&[-\epsilon _0^{\sigma}+\epsilon _0-\frac{(-\epsilon _0^{\sigma}+\epsilon _0)^3}{6}]\\
&&-\frac{1}{d-1}[\epsilon _0^{\sigma}+(d-1)\epsilon  _0 -\frac{(\epsilon _0^{\sigma}+(d-1)\epsilon _0 )^3}{6}]+O(\epsilon _0^4)+O(r), 
\end{eqnarray*}
which can also be seen to be positive under the given assumption, recall that $d-1\geq 2$. 
\end{proof}

The argument in the proof of Lemma \ref{LemmaIncreasedAnggle} also gives the following control on how much the orbit of $\Phi$ can cross the halfline $\{arg(z)=\frac{\pi}{d}\}$.  
\begin{lem} Let $r_0>0$ small enough and $\epsilon _0=r_0^{1/4}$. Let $z=re^{i\theta}$ with $0<r\leq r_0$ and $\theta = \frac{\pi}{d}+\epsilon $ where $|\epsilon |\leq \epsilon _0$. Let $z_1=\Phi (z)=r_1e^{i\theta _1}$ where $\theta _1=\frac{\pi}{d}+\epsilon _1$. Then $|\epsilon _1|= O(\epsilon _0^3)$. 
\label{LemmaControlCrossingUnstableCurve}\end{lem}
\begin{proof}
Since $\Phi$ is continuous in $\mathbb{D}(0,r_0)$ and the unstable curve $C_1$ is tangent to the halfline $\{arg(z)=\frac{\pi}{d}\}$, $\epsilon _1$ is small when $\epsilon _0$ is small. Applying the last part of the proof of Lemma Lemma \ref{LemmaIncreasedAnggle}, with $\epsilon ^{\sigma}$ replaced by $\epsilon _1$ and $\theta= \frac{\pi}{d}+\epsilon $, we see that 
$$-\frac{d}{d-1}\epsilon _1+\frac{(d-1)^2-1}{6}\epsilon ^3+O(r)+O(\epsilon _1^3)+O(\epsilon ^4)=0.$$
From this, the conclusion of the lemma follows immediately. 
\end{proof}

We also need the following simple geometric property. 
\begin{lem} Assume that $\epsilon _0=r_0^{1/4}$ and $r_0$ is small enough. Then the unstable curve $C_1$ lies above the halfline $\{arg(z)=\frac{\pi}{d}-\epsilon _0^3\}$, and lies below the halfline $\{arg(z)=\frac{\pi}{d}+\epsilon _0^3\}$. Similarly, the stable curve $C_{2d}$ lies above the halfline $\{arg(z)=-\epsilon _0^3\}$ and lies below the halfline $\{arg(z)=\epsilon _0^3\}$. 
\label{LemmaBoundTheStableCurves}\end{lem}
\begin{proof}
As in Step 5, we use the rotation of angle $\pi /d$ to bring the halfline $\{arg(z)=\frac{\pi}{d}\}$ to the positive x-axis. Under this rotation, $C_1$ becomes a curve of the form $C=\{y=g(x)\}$ where $g$ is real analytic and has order $N\geq 2$ at $x=0$. Also, under this rotation, the halfline $\{arg(z)=\frac{\pi}{d}-\epsilon _0^3\}$ becomes the halfline $L=$ $\{arg(z)=-\epsilon _0^3\}$. The first assertion of the Lemma is then equivalent to that the curve $C$ lies above the halfline $L$. Which is the same as: if $(x,y)\in L$, then $g(x)>y$. For $(x,y)\in L$, then $y=r\sin (-\epsilon _0^3)\sim - r r_0^{3/4}$, while $|g(x)|\sim |x|^N\sim r^N$. Since $-r^N>-rr_0^{3/4}$ for $0<r\leq r_0<1$, the proof is completed. The other assertions of the Lemma can be proved similarly.      
\end{proof}

We are now ready to describe the behaviour of $\Phi$ in the domain $U_{d}$. 

\begin{lem} Let $z_0\in U_{d}$, and construct the sequence $\{z_n\}$ by Analytic BNQN with initial point $z_0$. Then there is $N\geq 0$ so that $z_0,\ldots ,z_N\in U_{d}$, while $z_{N+1}\notin \mathbb{D}(0,r_0)$.  
\label{LemmaDynamicsInPetal}\end{lem}
\begin{proof}
The proof is divided into 2 steps. 

{\bf Step a}: 

We will show the following: If $z_1\in \mathbb{D}(0,r_0)$, then $z_1\in U_{d}$. 

We note that $U_{d}$ does not contain the two boundary curves $C_{d}$ and $C_2$. By symmetry, we can assume that $z_0$ is in the half of $U_{d}$ bound by the 2 curves $C_{2d}$ and $C_1$. 

There are 3 cases to consider: 

Case a.1: $z_0$ is above the curve $C_{2d}$, but inside the sector $\{|arg(z)|<\epsilon _0\}$. In this case, the dynamics of $\Phi$ conjugates to that of $H$ in Step 4. Since the Jacobian of $H$ at $(0,0)$ has two positive eigenvalues $0<a<1<A$, the ''No jump lemma" (Theorem \ref{TheoremNoJumpLemma}) applies. This translates to that $z_1$ is still above the curve $C_{2d}$. Moreover, the $x$-coordinate of $z_1$ is positive. Also, the angle of $z_1$ is clearly $<\pi /d$. Hence, $z_1\in U_{d}$.    

Case a.2: $z_0$ has angle between $\epsilon _0$ and $\frac{\pi}{d}-\epsilon _0$. In this case, Lemma \ref{LemmaIncreasedAnggle} gives that the angle of $z_1$ is between $\epsilon _0$ and $\frac{\pi}{d}-\epsilon _0^{3.5}$. Hence, $z_1\in U_{d}$. 

Case a.3: $z_0$ has angle between $\frac{\pi}{d}-\epsilon _0$ and $\frac{\pi}{d}+\epsilon _0$. In this case, Lemma \ref{LemmaControlCrossingUnstableCurve} gives that the angle of $z_1$ is between $\frac{\pi}{d}-\epsilon _0+O(\epsilon _0^3)$ and $\frac{\pi}{d}+\epsilon _0+O(\epsilon _0^3)$. Therefore, again $z_1\in U_{d}$. 

{\bf Step b}: By Step a and induction, if $z_0,\ldots ,z_{N}\in U_{d}$ and $z_{N+1}\in \mathbb{D}(0,r_0)$, then also $z_{N+1}\in U_{d}$. So, to prove the lemma, it suffices to show that the sequence $\{z_n\}$ eventually leaves $\mathbb{D}(0,r_0)$. Assume otherwise, then $z_n\in \mathbb{D}(0,r_0)$ for all $n$. By Theorem \ref{TheoremConvergence}, it follows that $z_n$ converges to $z^*=0$. 

We have several different cases: 

Case b.1: There is $n_0$ so that the angle of $z_{n_0}$ is between $\epsilon _0$ and $\frac{2\pi }{d}-\epsilon _0$. By Lemmas \ref{LemmaIncreasedAnggle} and \ref{LemmaControlCrossingUnstableCurve}, it follows that $arg(z_n)\in [\epsilon _0,\frac{2\pi }{d}-\epsilon _0]$ for all $n\geq n_0$. However, then the analysis so far shows that $z_n$ cannot converge to $0$. Thus we obtain a contradiction. 

Case b.2: In the remaining case, then either $arg(z_n)\leq \epsilon _0$ for all $n$, or $arg(z_n)\geq \frac{2\pi }{d}-\epsilon _0$ for all $n$. By symmetry, we can assume that $arg(z_n)\leq \epsilon _0$ for all $n$. Moreover, $z_n$ does not belong to the stable curve $C_{2d}$. This translates to that for $\psi$ and $H$, in Step 3, all $Z_n=\psi (z_n)$ belongs to the rectangle $R_{r_0,\epsilon _1}$ but no $Z_n$ belongs to the stable manifold of $H$. Thus $\{Z_n\}$ must have an accumulate point inside $R_{r_0,\epsilon _1}$ which must belong to the $Y$-axis, while the dynamics of $H$ dictates that this cannot happen.  
\end{proof}

{\bf Step 7: Describe the dynamics of $\Phi $ in the other domains $U_j$}: 

For $j=1,\ldots ,d-1$, the domain $U_j$ is that bound by the 2 stable curves $C_{2j}$ and $C_{2j+2}$ and contains the unstable curve $C_{2j+1}$. By a rotation by an angle which is a multiple of $2\pi /d$, we can reduce to the case of $U_{2d}$ and prove similar results like that in Step 6. 

{\bf Step 8: Finish the proof of Theorem \ref{Theorem2} for $d\geq 3$}: 

The construction of the curves $C_j$'s is given in Steps 4, 5 and 6, where parts 1, 2, 3b of Theorem \ref{Theorem2} are proven. 

Part 3a and 3c are proven in Steps 7 and 8.

\subsubsection{Proof for the case $d=2$} Now we prove Theorem \ref{Theorem2} for the case $d=2$. 

In this case, by Lemma \ref{LemmaApproximateFormulaDynamics} for $d=2$, it is clear that the map $\Phi$ is real analytic. The linearization of $\Phi$ is given by 
$$
\Phi _1(z)=\left[ \begin{array}{ccc}
		0 \\
		2y\\
	\end{array} \right],
$$
which has eigenvalues $0,2$. Therefore, Theorem \ref{TheoremStableUnstableManifoldRealAnalytic}, as well as arguments by rotations as in the proof of Theorem \ref{Theorem2} for $d\geq 3$, furnishes curves $C_1$, $C_2$, $C_3$ and $C_4$ which satisfy parts 1, 2, 3a, 3b of Theorem \ref{Theorem2}. 

We now proceed to proving part 3c of Theorem   \ref{Theorem2}. As in the case for $d\geq 3$, we work in $\mathbb{D}(0,r_0)$ with $r_0>0$ small enough, and choose $\epsilon _0=r_0^{1/4}$. By Lemma \ref{LemmaBoundTheStableCurves}, the stable curve $C_4$ is contained in the sector $S_4=\{|arg(z)|<\epsilon _0\}$, and the unstable curve $C_1$ is contained in the sector $\{|arg(z)-\frac{\pi}{2}|<\epsilon _0\}$. 

Apply the ''No jump lemma" Theorem \ref{TheoremNoJumpLemma} to $\Phi$, for which the unstable eigenvalue $2$ of $J\Phi (0,0)$ is positive, for $z\in S_4$ above the stable curve $C_4$, we have $\Phi (z)$ also stays above $C_4$. 

Consider now $z=(x,y)$ with $\epsilon _0\leq \arg (z)\leq \frac{\pi}{2}+\epsilon _0$. In this domain, there is $c>0$ so that $y=r\sin (\theta )\geq c\times rr_0^{3/4}$. Hence, if $z_1=(x_1,y_1)=\Phi (z)$, then $x_1=O(r^2)$ while $y_1\geq c\times rr_0^{3/4}+O(r^2)$. In particular, by the arguments in Lemma \ref{LemmaBoundTheStableCurves}, $z_1$ stays above the curve $C_4$.

We can argue similarly for the curve $C_2$. This shows that if $z_0\in U_2$, where $U_2$ is the domain bound by the curves $C_2$ and $C_4$, and contains the curve $C_1$, then $\Phi (z_0)\in U_2$. Then part 3c of Theorem \ref{Theorem2}, for $d=2$, can be completed by similar arguments as in the proof of Lemma \ref{LemmaDynamicsInPetal}. 

\section{Proofs of Theorems \ref{Theorem1},  \ref{Theorem1Bis} and \ref{Theorem3}} We provide in this section proofs of Theorems \ref{Theorem1},  \ref{Theorem1Bis} and \ref{Theorem3}. We finish the section, with some illustrating experiments and a discussion on how the results in this paper can be used to heuristically explain the apparent connections between BNQN, Voronoi's diagram and Newton's flow, as well as the boundary of basins of attraction for roots in BNQN dynamics seems to be smooth, in particular from many experiments done in \cite{RefFHTW}, in particularly to Newton's flow and Poincar\'e-Bendixon theorem. 

We first recall a fact, used in \cite{RefT}, that away from the poles, zeros  and critical points of the function $f$, the dynamics of BNQN is locally uniformly descending. For the reader's convenience, we give a sketch of proof for this fact. 

{\bf Fact 1}: Let $K$ be a compact subset of $\mathbf{C}\backslash (\mathcal{Z}(f)\cup \mathcal{C}(f)\cup \mathcal{P}(f))$. Then there is $a>0$ so that for all $z_0\in K$ and $z_1$ is constructed from $z_0$ by BNQN, we have $F(z_1)-F(z_0)\leq -a$. 

Sketch of proof of Fact 1: There are constants $b_1,B_1>0$ so that $b_1\leq ||\nabla F(z_0)||^{\tau}\leq B_1$ for all $z_0\in K$ (recall that critical points of $F$ are precisely the roots and critical points of $f$). Let $A_0=\nabla ^2F(z_0)+\delta _j ||\nabla F(z_0)||^{\tau}$ be the matrix chosen in Algorithm \ref{table:alg0}. Then if $\lambda _1,\lambda _2$ be the 2 eigenvalues of $A$, we have $|\lambda _1|,|\lambda _2|\geq \kappa ||\nabla F(z_0)||^{\tau}\geq \kappa b_1$. There is obviously another constant $B_2>0$ so that $|\lambda _1|, |\lambda _2|\leq B_2$. Let $u_1,u_2$ be the 2 corresponding eigenvectors, which are moreover orthonormal. From
$$w_0=\frac{1}{|\lambda _1|}<\nabla F(z_0),u_1>u_1+\frac{1}{|\lambda _2|}<\nabla F(z_0),u_2>u_2,$$
we obtain
\begin{eqnarray*}
\frac{|\nabla F(z_0)|^2}{B_2^2}&\leq& |w_0|^2=\frac{1}{|\lambda _1|^2}<\nabla F(z_0),u_1>^2+\frac{1}{|\lambda _2|^2}<\nabla F(z_0),u_2>^2\\&\leq& \frac{|\nabla F(z_0)|^2}{b_1^2},
\end{eqnarray*}
which infers the existence of constants $c_1,c_2>0$ for which $c_1\leq |w_0|\leq c_2$. Therefore, there are constants $\widehat{c_1}, \widehat{c_2}>0$ so that $\widehat{c_1}\leq |\widehat{w_0}|\leq \widehat{c_2}$. Therefore, since  $F$ is $C^1$, there is a lower bound $1/2^{n_0}$ for the $\gamma$ obtained after the second While loop in Algorithm \ref{table:alg0}. Thus, by the requirement of the second While loop in Algorithm \ref{table:alg0}, for $z_1=z_0-\gamma \widehat{w_0}$, we have 
\begin{eqnarray*}
F(z_1)-F(z_0)&\leq& -\gamma <\nabla F(z_0),\widehat{w_0}>\\
&=&-\frac{\gamma }{|\lambda _1|}<\nabla F(z_0),u_1>^2-\frac{\gamma }{|\lambda _2|}<\nabla F(z_0),u_2>^2\\
&\leq& -a,
\end{eqnarray*}
for a constant $a>0$ depending on the above constants. (Q.E.D.)

We also need the following fact that: 

{\bf Fact 2}: Assume that $\mathcal{C}(f)\not=\emptyset$.  Then, on any bounded set away from poles, BNQN is real analytic outside a finite number of real analytic curves. 

Proof of Fact 2: Let $\Omega$ be a bounded open set which has a positive distance to $\mathcal{P}(f)$. We choose for each $z^*\in \mathcal{C}(f)$ a disc $\mathbb{D}(z^*,r(z^*))$ for which Theorem \ref{Theorem2} is valid. We choose for each $z^*\in \mathcal{Z}(f)$ a disc $\mathbb{D}(z^*,r(z^*))$ for which Theorem \ref{Theorem3} is valid. We choose them small enough so that the closures of the discs do not intersect, and the union of them is closed in $\mathbf{C}$. Here is the list of curves: 

Curves of type 0: the boundaries of the discs, or the boundaries of the small sectors in these discs. By the proof of Theorems \ref{Theorem2} and \ref{Theorem3}, inside each of these discs, BNQN is real analytic outside curves of type 1. Since $\Omega $ is bounded, only a finite number of such points $z^*$, and hence of curves of type 1, which is relevant. 

Now we consider the real analyticity of BNQN on the complement $\Omega '=\Omega \backslash \bigcup _{z^*\in \mathcal{Z}(f)\cup \mathcal{C}(f)}\mathbb{D}(z^*,r(z^*))$. We know from the proof of Fact 1 that for $z\in K$, there are only a finite number of choices for $\gamma$ obtained after the second While loop of Algorithm \ref{table:alg0}, namely $\gamma \in \{1,\frac{1}{2}, \ldots ,\frac{1}{2^{n_0}}\}$, which independent of whether the $\delta _j$ we choose satisfies the condition that the minimum of the absolute values of the eigenvalues of the matrix $A_j(z)=\nabla ^2F(z)+\delta _j||\nabla F(z_0)||^{\tau}$ (i.e. $minsp(A_j(z))$) is $\geq \kappa ||\nabla F(z_0)||^{\tau}$ or not. We consider the following real analytic curves (if any of these sets has isolated points, we can replace them by small real analytic curves going through the isolated points) here $z\in \Omega '$: 

Curves of type 1: those for which $det(A_j(z))=0$, for $j=0,1,2$. By Lemma \ref{LemmaApproximateFormulaDynamics}, none of these sets can be the whole $\mathbf{C}\backslash \mathcal{P}(f)$, because near a critical point, but not at the critical point, $A_j(z)$ is invertible. Therefore, they have dimension at most $1$. 

Curves of type 2: those for which $minsp(A_{j}(z))=\kappa ||\nabla F(z)||^{\tau}$ for $j_1\not= j_2$, $j_1,j_2\in \{0,1,2\}$. Again, by  Lemma \ref{LemmaApproximateFormulaDynamics}, none of these sets can be the whole $\mathbf{C}\backslash \mathcal{P}(f)$. Therefore, they have dimension at most $1$. 

Curves of type 3: those for which $1=\theta |{w}|$, where $w$ is constructed from any of the $A_j$ for which $det(A_j)\not= 0$. Again, by Lemma \ref{LemmaApproximateFormulaDynamics}, none of these sets can be the whole $\mathbf{C}\backslash \mathcal{P}(f)$, since near a critical point $z^*$ we have $|w|=O(|z|)$. Therefore, they have dimension at most $1$. 

Curves of type 4: those for which $F(z-\gamma \widehat{w})-F(z)=-\frac{\gamma}{3}<\nabla F(z),\widehat{w}>$, where $w$ is constructed from any of the $A_j$ for which $det(A_j)\not= 0$ and $\gamma \in \{1,\frac{1}{2},\ldots , \frac{1}{2^{n_0}}\}$. Except for $\gamma =1$, none of these sets can be the whole $\mathbf{C}\backslash \mathcal{P}(f)$, since near a critical point $z^*$ we have by Lemma \ref{LemmaApproximateFormulaDynamics} that $F(z-\gamma \widehat{w})-F(z)<-\frac{1}{3}<\nabla F(z),\widehat{w}>$. If for $\gamma =1$ the corresponding set is the whole  $\mathbf{C}\backslash \mathcal{P}(f)$, then we can discard the set since it does not affect our arguments.  

Now, let $S$ be the union of the curves of type 0, 1, 2, 3, and 4, we observe that for $z_0\in K\backslash S$ there is no ''switch" in the dynamics of BNQN: if at $z_0$ the parameters in Algorithm \ref{table:alg0} are $\delta _j$ and $\gamma$, then in an open neighbourhood of $z_0$, BNQN still uses the same parameters. As in the proof of Theorem \ref{Theorem2}, we can show that when there is no switch then BNQN is real analytic, since the involved matrix $A(z)$ in Algorithm \ref{table:alg0} is real analytic and invertible. 

\subsection{Proof of Theorem \ref{Theorem1}} Recall the the sets $\mathcal{Z}(f)$, $\mathcal{C}(f)$ and $\mathcal{P}(f)$ are discrete. If $\mathcal{C}(f)$ is empty, then there is nothing to prove. Therefore, we assume from now that $\mathcal{C}(f)\not=\emptyset$.  

For each $z^*\in \mathcal{Z}(f)$, we choose a disc $\mathbb{D}(z^*,r_{z^*})$ - with center $z^*$ and radius $r^*\leq 1$ - in the basin of attraction for $z$ of the dynamics of BNQN. For each $z^*\in \mathcal{C}(f)$, we choose a disc $\mathbb{D}(z^*,r_{z^*})$ - with center $z^*$ and radius $r^*\leq 1$ -  for which Theorem \ref{Theorem2} is valid. We can choose these in such a way that $\bigcup _{z^*\in \mathcal{Z}(f)\cup \mathcal{C}(f)}\overline{\mathbb{D}(z^*,r_{z^*})}$ is a closed set in $\mathbf{R}^2$. 

If $z_0\in \mathcal{E}$, and $\{z_n\}$ is the sequence constructed by BNQN, then by Theorem \ref{Theorem2} we know that there is $z^*\in \mathcal{C}(f)$ and a positive integer $n_0$ such that $z_{n_0}\in \mathbb{D}(z^*,r_{z^*}/2)$ and lands on a stable curve therein. 

For  integers $k,n_1,n_2\geq 0$, we define $\mathcal{E}_{k,n_1,n_2}$,  be the set of $z_0\in \mathcal{E}$ for which: i) $z_k\in \mathbb{D}(z^*,r_{z^*}/2)$ for a $z^*\in \mathcal{C}(f)$ and lands on a stable curve therein, ii) $|z_0|\leq n_1$, and iii) $|f(z_0)|\leq n_2$. It is easy to see that $\mathcal{E}=\bigcup _{k,n_1,n_2}\mathcal{E}_{k,n_1,n_2}$. Hence if we can show that each $\mathcal{E}_{k,n_1,n_2}$ is contained in a countable union of real analytic curves, then the proof is completed. 

Indeed, we will show a  stronger assertion: 

{\bf Claim 1}: For each $k,n_1,n_2$, $\mathcal{E}_{k,n_1,n_2}$ is contained in a {\bf finite} union of real analytic curves. 

Proof of Claim 1: 

Let $\delta _0, \delta _1, \delta _2$ be the parameters used in Algorithm \ref{table:alg0}. Define $\Omega _0=\{z\in \Omega : det(\nabla ^2F(z)+\delta _0||\nabla F(z)||^{\tau})=0\}$, $\Omega _1=\{z\in \Omega : det(\nabla ^2F(z)+\delta _1||\nabla F(z)||^{\tau})=0\}$, and $\Omega _2=\{z\in \Omega : det(\nabla ^2F(z)+\delta _2||\nabla F(z)||^{\tau})=0\}$. For each $j$, either $\Omega _j$ is $\empty$, a closed real analytic curve or the whole $\Omega$. In the latter case, the matrix $ det(\nabla ^2F(z)+\delta _1||\nabla F(z)||^{\tau}$ will never be used in BNQN, and hence we can redefine it to be $\emptyset$. Then $\Omega _0$, $\Omega _1$, and $\Omega _3$ are closed real analytic curves in $\Omega$. 

As in Step 1 of the proof of case $d\geq 3$ of Theorem \ref{Theorem2}, one can show that the assignment $z_n\rightarrow \widehat{w}_n$ in Algorithm \ref{table:alg0} is real analytic, for $z\in \Omega \backslash (\Omega _0\cup \Omega _1\cup \Omega _2)$ and for a {\bf fixed choice}  of $\delta _j$ in $\{\delta _0, \delta _1,\delta _2\}$ (where we {\bf disregard the condition} that both eigenvalues of $\nabla ^2F(z)+\delta _j||\nabla F(z)||^{\tau}$ have absolute values $\geq \kappa ||\nabla F(z)||^{\tau}$), since all 3 matrices $\nabla ^2F(z)+\delta _0||\nabla F(z)||^{\tau}$, $\nabla ^2F(z)+\delta _0||\nabla F(z)||^{\tau}$ and $\nabla ^2F(z)+\delta _2||\nabla F(z)||^{\tau}$ are invertible and real analytic.  

We prove Claim 1 by induction on $k$. 

The base case $k=0$: In this case, $z_0$ is in on one of the stable curves belonging to a critical point $z^*\in \mathcal{C}(f)$. Since we chose the radius $r(z^*)\leq 1$ and $z_0$ belongs to a bounded set, such a $z^*$ must be bounded. Since there are at most finitely many critical points $z^*$ of $f$ in any bounded set, and for each $z^*$ only a finite number of stable curves each being real analytic (by Theorem \ref{Theorem2}), we conclude that $\mathcal{E}_{0,n_1,n_2}$ is contained in a finite union of real analytic curves.   

Assume now that Claim 1 is proven for $k$. We need to prove that Claim 1 holds for $k+1$.

We treat first the case the parameter $\theta $ in Algorithm \ref{table:alg0} is positive. Let $z_0\in \mathcal{E}_{k+1,n_1,n_2}$. In this case it follows that the vector $\widehat{w_k}$ in the algorithm satisfies $|\widehat{w_k}|\leq \frac{1}{\theta}$. Therefore, if $z_1$ is constructed from $z_0$ by BNQN, we have: 1) $|z_1|\leq |z_0|+\frac{1}{\theta}$, 2) $|f(z_1)|\leq |f(z_0)|$ (recall that BNQN, by the use of Armijo's Backtracking line search, has descend property). Therefore, $z_1\in \mathcal{E}_{k,n_1',n_2} $ for any positive integer $n_1'\geq n_1+\frac{1}{\theta}$. Hence, $\mathcal{E}_{k+1,n_1,n_2}$ is the preimage of $ \mathcal{E}_{k,n_1',n_2}$ by BNQN. By Fact 2, outside a finite number of real analytic curves, BNQN is real analytic and hence the preimage of $ \mathcal{E}_{k,n_1',n_2}$ is contained in a real analytic set. To finish the proof, it suffices to show that when the parameters $\delta _j, \theta , \tau $ and $\gamma $ in Algorithm \ref{table:alg0}  are {\bf fixed} over all the open set where the matrix $A(z)$ is invertible and real analytic (here without checking either of the two While loops), then outside a real analytic curve, the assignment $H(z)=z-\gamma \widehat{w}$ is invertible (then the dimension of a real analytic set is preserved under preimage of the assignment). To this end, we note that the condition that the assignment is not invertible, which is $\det (JH(z))=0$, while not a real analytic function, can be written as a fraction of two real analytic functions $q_1(z)/q_2(z)$.  We consider hence the bigger set $\{z\in \mathbf{C}\backslash (\mathcal{Z}(f)\cup \mathcal{C}(f)\cup \mathcal{P}(f)): q_1(z)=0\}$. Then, since near a critical point  of $f$ (but not at the critical point), $H(z)$ is invertible by Lemma  \ref{LemmaApproximateFormulaDynamics},  and the open set  $\mathbf{C}\backslash (\mathcal{Z}(f)\cup \mathcal{C}(f)\cup \mathcal{P}(f))$ is connected, we see that the set where $q_1(z)=0$ cannot be the whole $\mathbf{C}\backslash (\mathcal{Z}(f)\cup \mathcal{C}(f)\cup \mathcal{P}(f))$. Therefore, the set where $H(z)$ is not invertible is contained in a real analytic set. 

Adding the finite exceptional list from Fact 2, we see that  $\mathcal{E}_{k+1,n_1,n_2}$ is contained in a finite union of real analytic curves. 

For the case $\theta =0$, as in the case $\theta >0$, it suffices to show that  $|z_1|\leq |z_0|+b$, for some $b>0$ independent of $z_0\in \mathcal{E}_{k+1,n_1,n_2}$. We consider $2$ cases: $z_0$ belongs to one of the discs $\mathbb{D}(z_0^*,r(z_0^*))$ or not. If $z_0$ does not belong to any of these discs, then the proof of Fact 1 gives that $|z_1|\leq |z_0|+b$ for some $b>0$ independent of such $z_0$. Otherwise, $z_0$ belongs to a finite number of the discs, and by the proof of Theorems \ref{Theorem2} and \ref{Theorem3} we know that BNQN is continuous on each of these discs. Therefore, we showed the existence of such $b$, and finish the induction proof.

\subsection{Proof of Theorem \ref{Theorem1Bis}} For each $z^*\in \mathcal{Z}(f)$ we choose a disc $\mathbb{D}(z^*,r(z^*))$ for which Theorem \ref{Theorem3} applies. For each $z^*\in \mathcal{C}(f)$ we choose a disc $\mathbb{D}(z^*,r(z^*))$ for which Theorem \ref{Theorem2} applies.

$\Rightarrow$ If the exceptional set $\mathcal{E}$ is locally contained in a finite union of real analytic curves, then  $\mathcal{E}\cap \mathbb{D}(z^*,r(z^*))$ is contained in a finite union of real analytic curves for each of the discs. 

$\Leftarrow$ Assume that $\mathcal{E}\cap \mathbb{D}(z^*,r(z^*))$ is contained in a finite union of real analytic curves for each of the discs. We will show that if $K$ is a compact set away from the poles, then $\mathcal{E}\cap K$ is contained in a finite union of real analytic curves. Since $K$ is compact and away from the poles, the number $R=\max _{z\in K}|f(z)|$ is finite. Hence $K\subset B_R=\{z\in \mathbf{C}: |f(z)|\leq R\}$. Since $f$ has compact sublevels, the set $B_R$ is compact. Hence, it suffices to show that $\mathcal{E}\cap B_R$ is contained in a finite number of real analytic curves. 

Let $a$ be the positive number from Fact 1, for $B_R\backslash \bigcup _{z^*\in \mathcal{C}(f)\cup \mathcal{Z}(f)}\mathbb{D}(z^*,r(z^*))$ (note that only a finite number of such $z^*$'s is relevant, since $B_R$ is compact), and define $N_0$ to be a positive integer $\geq R^2/a$. We first show that for every  $z_0\in B_R$, at least one of the points $z_0,\ldots ,z_{N_0}$ (constructed by BNQN from $z_0$) belongs to  $\bigcup _{z^*\in \mathcal{C}(f)\cup \mathcal{Z}(f)}\mathbb{D}(z^*,r(z^*))$. Indeed, since $B_R$ is a sublevel set, and BNQN has the descent property, we have that the whole sequence $z_0,\ldots ,z_{N_0}$ belongs to $B_R$. If none of the points $z_0,\ldots ,z_N$ belongs to  $\bigcup _{z^*\in \mathcal{C}(f)\cup \mathcal{Z}(f)}\mathbb{D}(z^*,r(z^*))$, then Fact 1 applies and gives: $F(z_{j+1})-F(z_j)\leq -a $ for all $j=0,\ldots, N_0-1$. Therefore, $$F(z_{N_0})\leq F(z_0)-N_0a\leq R^2/2-R^2<0,$$
which contradicts with that $F(z)=|f(z)|^2/2\geq 0$ for all $z$. 

Now, let $z_0\in \mathcal{E}\cap B_R$. Choose $k\leq N_0$ be such that $z_k\in $ one of the discs $\mathbb{D}(z^*,r(z^*))$. Note that $z^*$ cannot be a root of $f$, since by the choice  $\mathbb{D}(z^*,r(z^*))$ belongs to the basin of attraction of $z^*$. Hence, $z_k\in \mathcal{E}\cap \mathbb{D}(z^*,r(z^*))$, the latter being contained in a finite union of real analytic curves. Therefore, $z_0$ belongs to the k-th preimage of $\mathcal{E}\cap \mathbb{D}(z^*,r(z^*))$ by BNQN, where $k\leq N_0$. By the same proof as in Theorem \ref{Theorem1}, we deduce that $\mathcal{E}\cap B_R$ is contained in a finite union of real analytic curves.

\subsection{Proof of Theorem \ref{Theorem3}} We now prove Theorem \ref{Theorem3}. So, we consider here $z^*=0$ and $f(z)=z^d+$ higher order terms, for $d\geq 1$. Again, we use both the Cartesian and  polar coordinates $z=(x,y)= re^{i\theta}$, where $r>0$ is small enough.  

If $f=u+iv$, where $u,v$ are the real and imaginary parts of $f$, then the calculations for the derivatives $u_x,u_y,u_{xx},\ldots $ are the same as in Section 3.1. The formula for $v$ is also the same as in Section 3.1, while the formula for $u$ is changed to $r^d\cos (d\theta )+O(r^{d+1})$.

$F(z)=\frac{1}{2}|f(z)|^2=\frac{1}{2}r^{2d}+O(r^{2d+1})$. 

Then, we find
$$\nabla F(z)=dr^{2d-2}z+O(r^{2d}),$$
and 
$$\nabla ^2F(z)=r^{2d-2}\left[ \begin{array}{ccc}
		d^2+d(d-1)\cos (2\theta)&d(d-1)\sin (2\theta) \\
		d(d-1)\sin (2\theta )&d^2-d(d-1)\cos (2\theta )\\
	\end{array} \right] + O(r^{2d-1}). $$

Both the trace and the determinant of $\nabla ^2F(z)$ are positive, hence both eigenvalues of  $\nabla ^2F(z)$ are positive and of the same size as $r^{2d-2}$. Therefore, the $\delta _j$ in Algorithm \ref{table:alg0} is $\delta _0$, and since $\tau \geq 1$ and $||\nabla F(z)||=O(r^{2d-1})$, the matrix $A$ used in that algorithm is $A=\nabla ^2F(z)+\delta _0||\nabla F(z)||^{\tau}Id=\nabla ^2F(z)+O(r^{2d-1})$ also has positive eigenvalues. 

Therefore, the $w$ in that algorithm is the same as the $v$, and is equal to: $A^{-1}.\nabla F(z)=\frac{z}{2d-1}+O(r^2)$. Hence, the $\widehat{w}$ in that algorithm is $\widehat{w}=w+O(r^2)=\frac{z}{2d-1}+O(r^2)$. Now, we show that the $\gamma$ after running the second While loop in Algorithm \ref{table:alg0} is $1$. This is the same as showing that: $F(z-\widehat{w})-F(z)\leq -\frac{1}{3}<\nabla F(z),\widehat{w}>$. From the above formulas, we have
\begin{eqnarray*}
F(z-\widehat{w})-F(z)&=&r^{2d}[(1-\frac{1}{2d-1})^{2d}-1]+O(r^{2d+1}),\\
-\frac{1}{3}<\nabla F(z),\widehat{w}>&=&-\frac{r^{2d}}{3(2d-1)}+O(r^{2d+1}).
\end{eqnarray*}

Therefore, the needed inequality is equivalent to 
$$(1-\frac{1}{2d-1})^{2d}-1\leq -\frac{1}{3(2d-1)}+O(r).$$
To this end, it suffices to show that for all $d\geq 1$ we have 
$$(1-\frac{1}{2d-1})^{2d}-1< -\frac{1}{3(2d-1)}.$$
We define $a=(1-\frac{1}{2d-1})^{2d}$. Then 
$$\frac{1}{a}=(1+\frac{1}{2d-2})^{2d}>1+\frac{2d}{2d-2}=\frac{2(2d-1)}{2d-2}.$$
Therefore, 
\begin{eqnarray*}
(1-\frac{1}{2d-1})^{2d}-1+\frac{1}{3(2d-1)}&<&\frac{2(2d-1)}{2d-2}-1+\frac{1}{3(2d-1)}\\
&=&\frac{6d-4}{12d-6}-1\\
&=&\frac{-6d+2}{12d-6},
\end{eqnarray*}
which is $<0$ for all $d\geq 1$, as wanted. 

From all the above, the update rule for BNQN near $z^*$ is 
$$z_{k+1}=z_k-\widehat{w_k}=\frac{2d-2}{2d-1}z_k+O(|z_k|^2),$$
and the proof of Theorem \ref{Theorem3} is completed. 

\subsection{Some (both rigorous and heuristic) explanations for the observed experiments on BNQN} There have been several experiments on BNQN, which display some remarkable features, in particularly those in the paper \cite{RefFHTW2}. In this subsection, we provide some heuritic explanations for these phenomena, based on the results proven in previous sections. 

{\bf Observation 1: The boundary of basins of attraction of the roots seems remarkably smooth, in particularly when compared to Newton's method.} 

While the boundary for basins of attraction of the roots for Newton's method (and also Relaxed Newton's method, Random Relaxed Newton's method) usually is fractal, the pictures obtained with BNQN seems much more smooth. Theorem \ref{Theorem1}, which says that the exceptional set $\mathcal{E}$ is contained in a countable union of real analytic curves, confirms  this experimental feature. 

Moreover, in all experiments so far, it seems that the exceptional set $\mathcal{E}$ is - away from poles - locally contained in a finite union of real analytic curves. This can be explained as follows. When drawing such a picture on computers, one usually has to make the following choices: 

- Choose initial points on a bounded grid, 

- Choose a maximum number k of how many iterations one will run at an initial point, 

- Choose an $\epsilon$ - threshold to decide whether the  point when the algorithm stops can be considered as being in the basin of a root. 

With these choices and assume that $f$ is a polynomial (or more generally has compact sublevels),  Theorem \ref{Theorem1Bis} can be used to rigorously explain the observations, provided one can show that for $z^*\in \mathcal{C}(f)$ and the disc $\mathbb{D}(z^*,r)$ for which Theorem \ref{Theorem2} applies, then $\mathcal{E}\cap \mathbb{D}(z^*,r)$ is contained in a finite union of real analytic curves. From Theorem \ref{Theorem2}, we know that the $d$ stable curves are real analytic and belong to $\mathcal{E}$. Hence, the difficulty lies at the fact that there could be ''bad" initial points $z_0\in \mathbb{D}(z^*,r)$ for which some $z_k$ leaves the disc, but then later some $z_{k'}$ re-enters the disc and lands on the stable curves. If one can show, for a given function $f$, that such ''bad" initial points do not exist, then one rigorously establishes that $\mathcal{E}$ is - away from poles of $f$ - contained in a finite union of real analytic curves.    

{\bf Observation 2: BNQN behaviour seems similar to Newton's flow and Poincar\'e-Bendixon theorem.} 

We recall that Newton's flow method  for finding roots of a meromorphic function $f$ is to find solutions $z(t)$ of the following ODE:
$$\frac{dz(t)}{dt}=-\frac{f(z(t))}{f'(z(t))},$$
with initial value $z(0)=z_0$. One would like to show that $\lim _{t\rightarrow\infty}z(t)$ is a root of $f(z)$. 

In \cite{RefFHTW2}, it has been found by many experiments that the basins of attraction for BNQN and Newton's method look very similar. Here, we will provide an extensive discussion on this similarity, based on the results in this paper.  

Usually, the dynamics of an ODE is better behaved than its discrete version. In our case, that Newton's flow (the ODE) behaves better than its discrete version (Relaxed Newton's method) is already known since the 1970s. Here, we recall some main facts about the behaviour of Newton's flow for finding roots of a meromorphic function $f$, following \cite{RefJongenEtAl} (the reader can find more detail in this paper and the survey \cite{RefBer}, as well as references therein), and compare with BNQN.

An important property of Newton's flow is that: $f(z(t))=e^{-t}f(z_0)$ (in the open interval - could be bounded -  where the solution exists). Based on this, even originally the right hand side  of the ODE for Newton's flow (i.e. $-f(z)/f'(z)$) is not real analytic (even not continuous), one can multiply it with a non-negative function (more precisely, $-(1+|f(z)|^4)^{-1}|f(z)|^2$) to obtain a new ODE: 
\begin{equation}
\frac{dz(t)}{dt}=g(z(t))
\label{LabelRealAnalyticNewtonFlow}\end{equation}
where $g(z)$ is a real analytic function, and solutions of this ODE and Newton's flow have the same trajectories. Therefore, studying the limit behaviour of Newton's flow is the same as studying the limit behaviour of the new ODE.  
 
Based on this, one can prove many facts, for example: A limit point of Newton's flow is either a root of $f$, a critical point of $f$, or $\infty$. This fact is known also for BNQN, see \cite{RefTT}. 

Moreover, since $g(z(t))$ is real analytic, the classical Poincar\'e-Bendixon theorem applies to the ODE (\ref{LabelRealAnalyticNewtonFlow}).  In particular, it follows that the behaviour of Newton's flow near a critical point of order $d$ of $f$ (but is not itself a root of $f$) is that of $d$-fold saddle points. This is exactly the same as that described in Theorem \ref{Theorem2} in the current paper. 

Theorem \ref{Theorem1} says that the set of initial points $z_0$ for which BNQN does not converge to a root of $f$ or $\infty$ is locally a countable union of real analytic curves. A similar (and stronger) property is known for Newton's method. 

We see therefore that indeed there are many similarities between BNQN (an iterative method) and Newton's flow (a continuous method). This is the more striking given that a fixed discrete version of Newton's flow (e.g. Relaxed Newton's method when $\gamma$ is fixed or the explicit Runger-Kutta with a fixed set of parameters) is known to not have these properties, thanks to \cite{RefMc}. Remark that when doing experiments in practice, one needs to use one of the discrete versions of Newton's flow, for which it is not guaranteed (by the mentioned reason) that the pictures obtained reflect well the theoretical results for Newton's flow. It is also worth noting that for the Relaxed Newton's method to be a reasonable approximation of Newton's flow, one needs to choose $\gamma$ to be small, which may make the running time to be very long (and this is actually observed in practice).    

%This can  be explained as follows. By Theorem \ref{Theorem2}, near any critical point $z^*$ but not a root of $f$ (which is allowed to be any meromorphic function), and of order $d$, there are precisely $d$ stable curves for the dynamics of Analytic BNQN. For Newton's flow, there is a similar fact, see e.g. the discussion on critical lines in \cite[Section 5]{RefHK} for the case of rational functions.  

%\textcolor{red}{Add: BNQN is discrete but satisfies Bendixon-Poincare theorem for flows. Newton's flow cannot be exactly solved in practice. Any discrete versions of Newton's flow, like Relaxed Newton's method or Runge-Kutta, do not have guarantee to find roots of all polynomials, by McMullen's result.}

{\bf Observation 3: The pictures by BNQN are resemblant of Voronoi's diagrams.} 

Recall (see \cite{RefV1}, \cite{RefV2}) that if $z_1^*,z_2^*,\ldots \subset \mathbf{C}$ is a discrete set, the Voronoi's diagram of the sequence consists of Voronoi's cells, each consisting of points closer to one point $z_j^*$ than to other points in the sequence. 

It seems from many experiments that the pictures obtained by BNQN for basins of attraction of the roots of $f$ are similar to the Voronoi's diagram of the roots. Here we offer an explanation. For simplicity, we restrict to the case $f$ is a polynomial of degree $4$, with $4$ distinct roots $z_1$. By Gauss-Lucas theorem, the critical points $c_1^*$, $c_2^*$ and $c_3^*$ of $f$ lie inside the convex hull of $z_1^*$, $z_2^*$, $z_3^*$ and $z_4^*$.

\begin{figure}
    \centering
    \includegraphics[width=10cm]{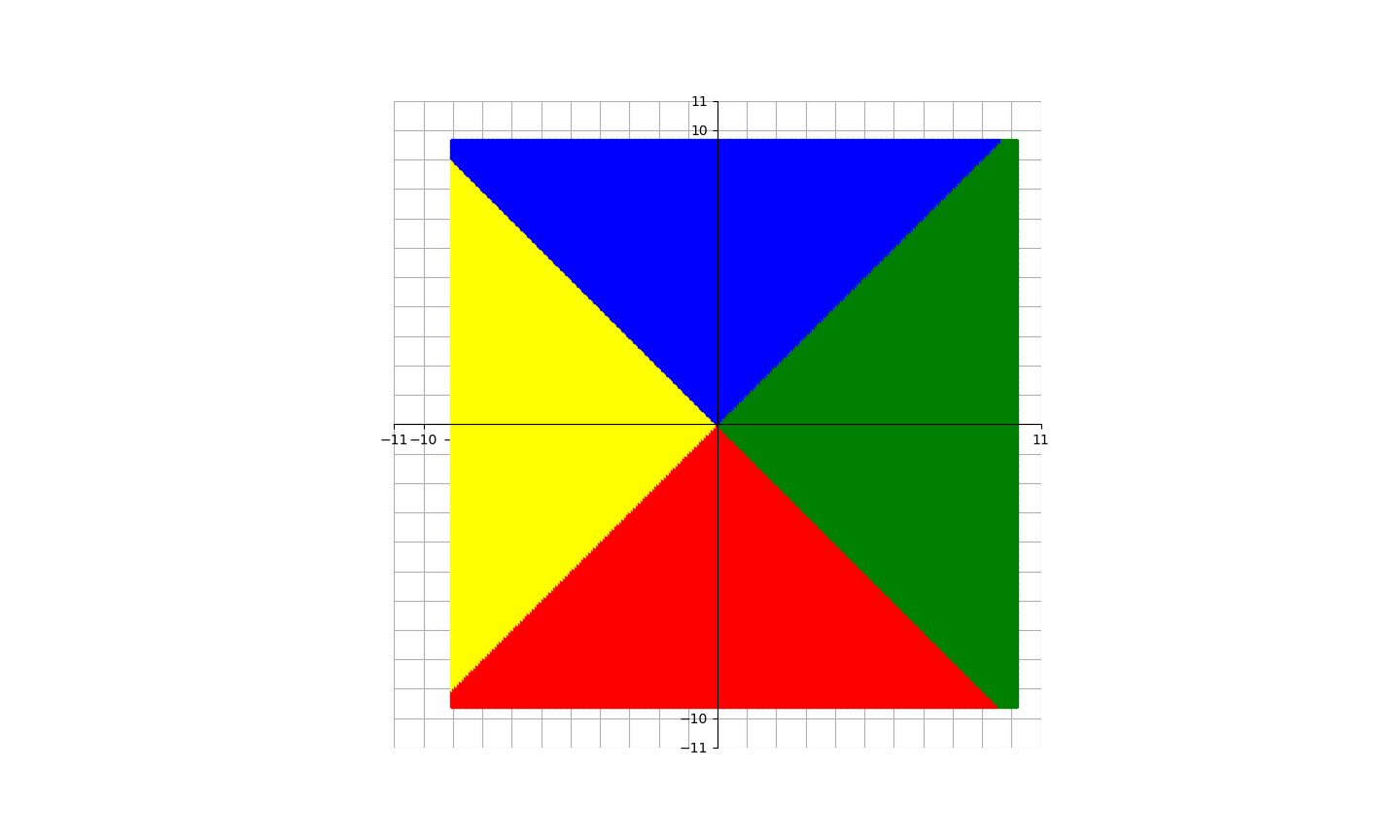}
    \includegraphics[width=10cm]{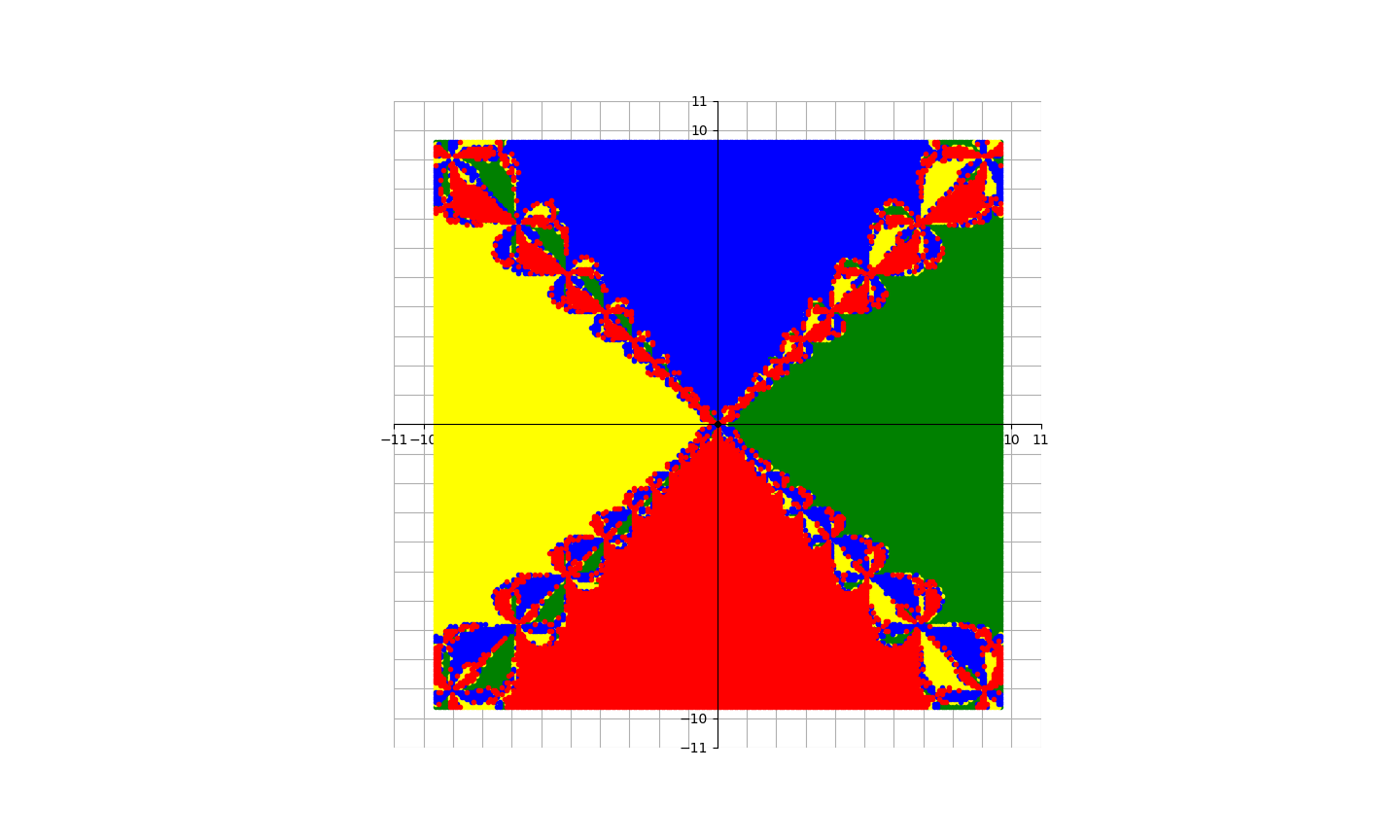}
    \includegraphics[width=10cm]{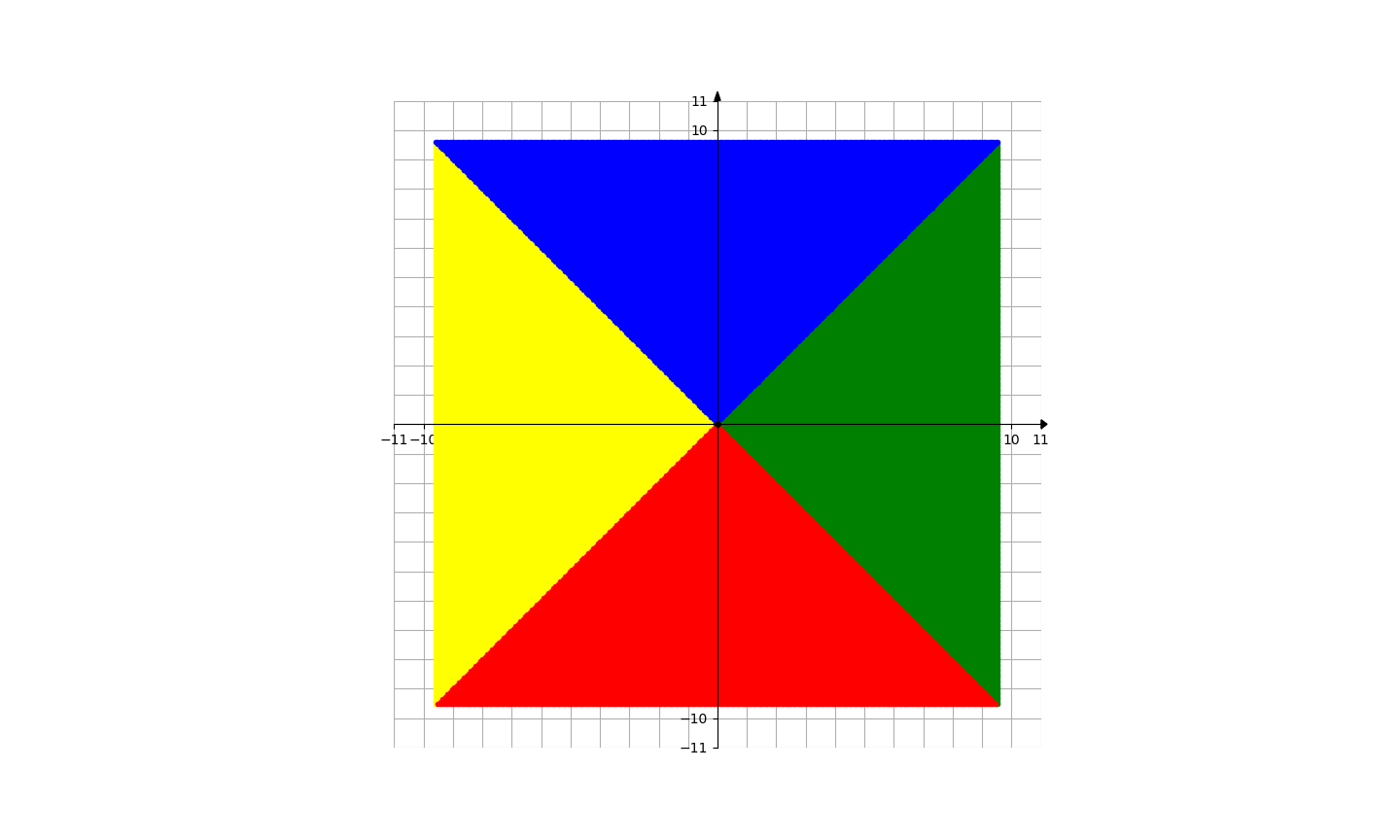}
    \caption{Basins of attraction for finding the roots of $f(z)=1-z^4$. Top picture: Voronoi's diagram of the 4 roots, center picture: Newton's method, bottom picture: BNQN. Points of the same colour  belong to basin of attraction of the same root. Here the parameter $\theta$ in Algorithm \ref{table:alg0} is $1$.}
    \label{fig:FF1}
\end{figure}

There are 4  different types of geometric configurations: 

Case 3a: The roots form a convex quadrilateral. This is the case where the above mentioned similarity is most perfect. In particular, for highly symmetric functions such as $1-z^4$ or $1+z^4$, then the picture produced by BNQN and the Voronoi's diagram look exactly the same. Consider for example the case $f(z)=1-z^4$ (see Figure \ref{fig:FF1}). In this case, we know from Theorem \ref{Theorem2} that the stable curves at the unique critical points are tangent to the boundary halflines of Voronoi's diagram. Why the stable curves seem to be precisely these halflines can be caused by the high symmetry of $f$ under rotations. 

\begin{figure}
    \centering
    \includegraphics[width=10cm]{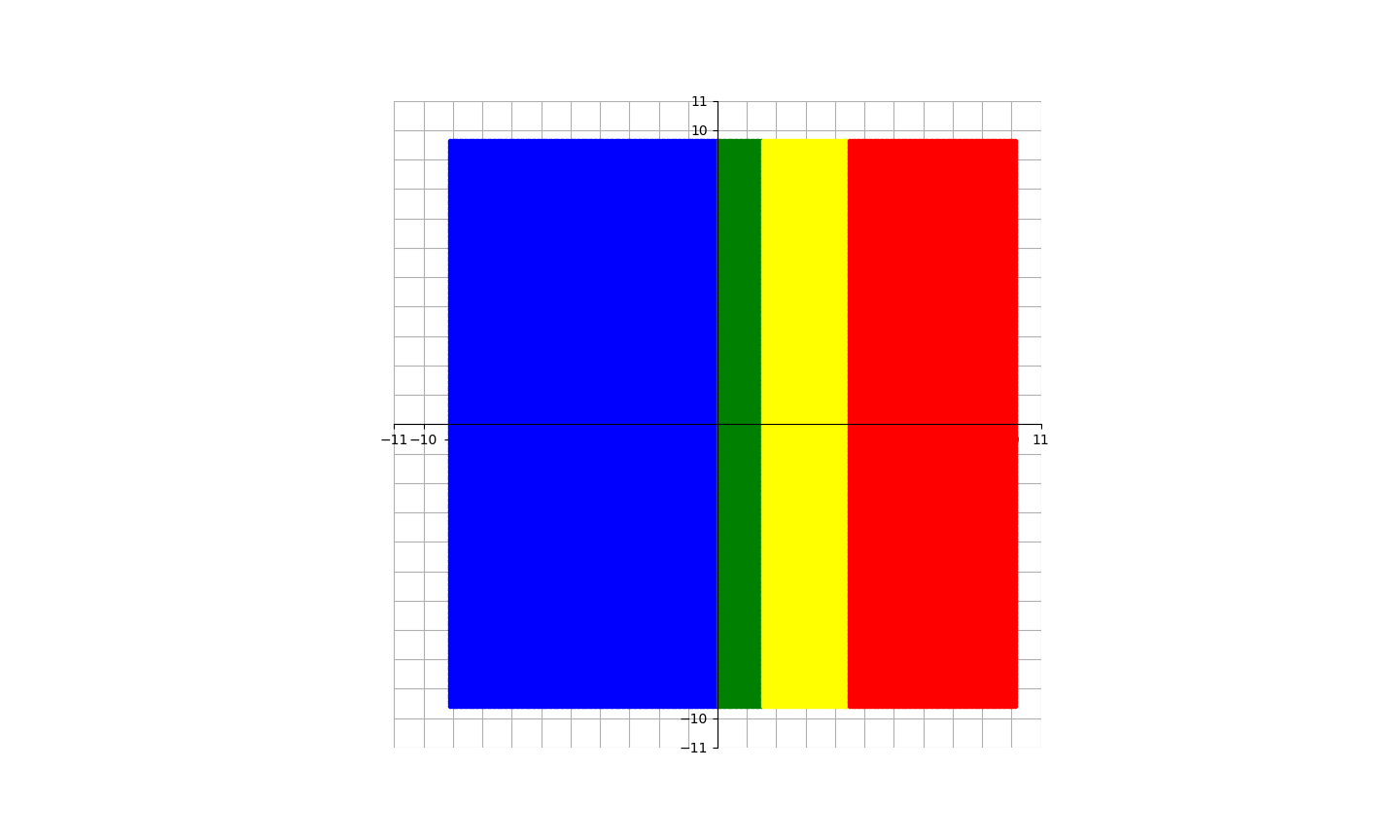}
    \includegraphics[width=10cm]{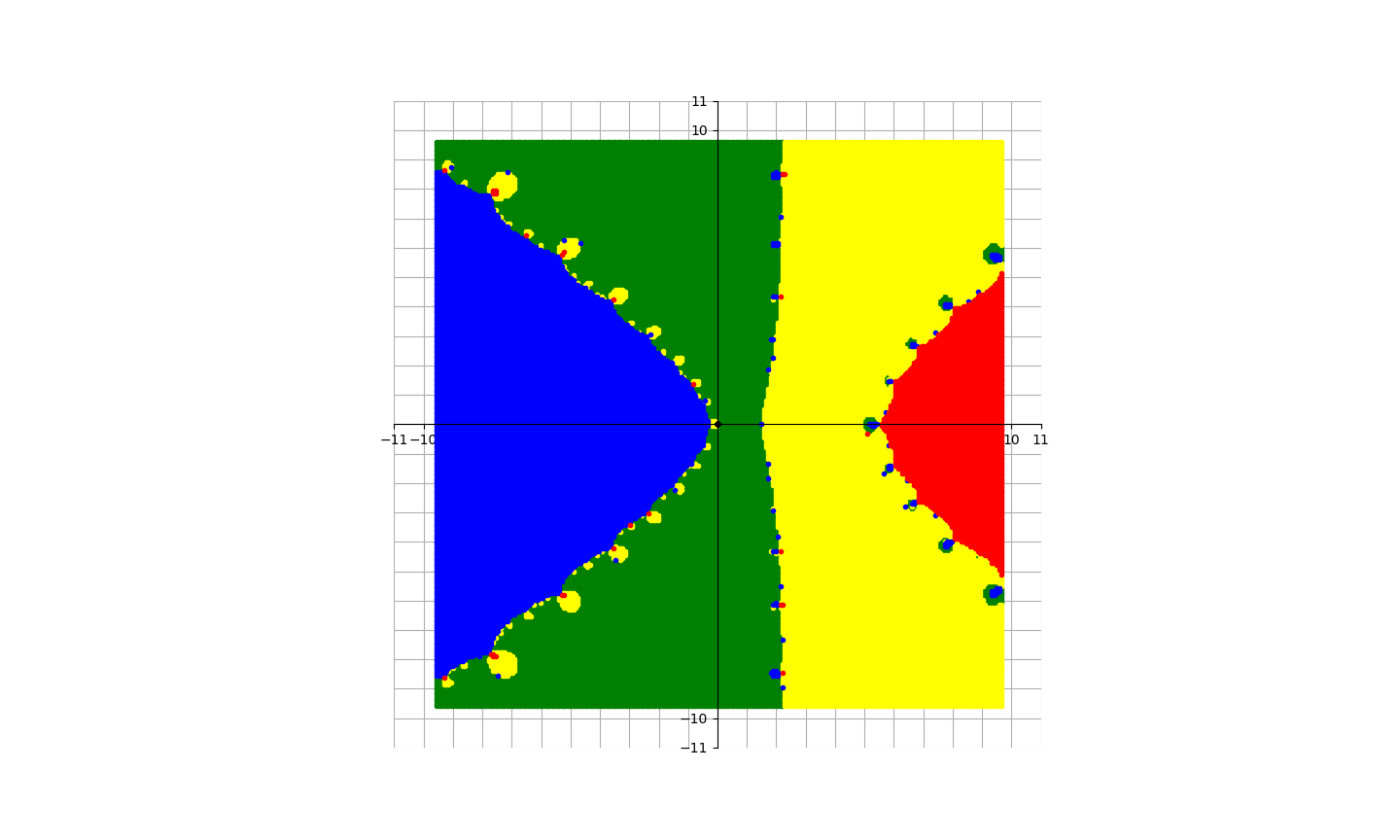}
    \includegraphics[width=10cm]{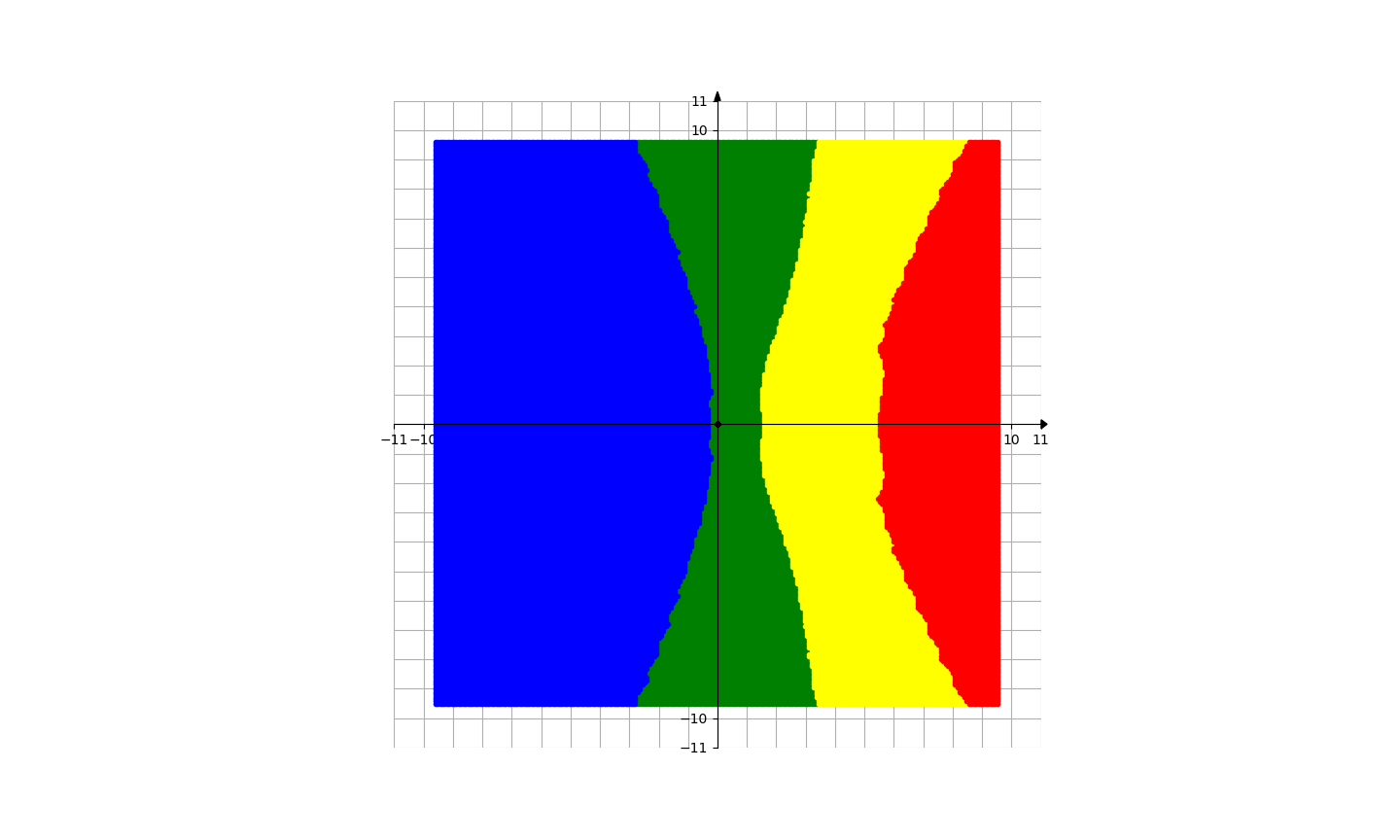}
    \caption{Basins of attraction for finding the roots of $f(z)=(z-1)(z-2)(z+1)(z-7)$. Top picture: Voronoi's diagram of the 4 roots, center picture: Newton's method, bottom picture: BNQN. Points of the same colour  belong to basin of attraction of the same root.  Here the parameter $\theta$ in Algorithm \ref{table:alg0} is $1$.}
    \label{fig:FF2}
\end{figure}

 Case 3b: $z_1^*$, $z_2^*$, $z_3^*$ and $z_4^*$ lie on the same line. In this case, the Intermediate Value Theorem is applicable, and implies that between any two consecutive roots there is a critical point. Pictures drawn for BNQN are a bit curved but look similar to Voronoi's diagram, except that the boundary curves do not start at the middle points of two consecutive roots but at the critical point lying between them. For an explicit example, see Figure \ref{fig:FF2}.

\begin{figure}
    \centering
    \includegraphics[width=10cm]{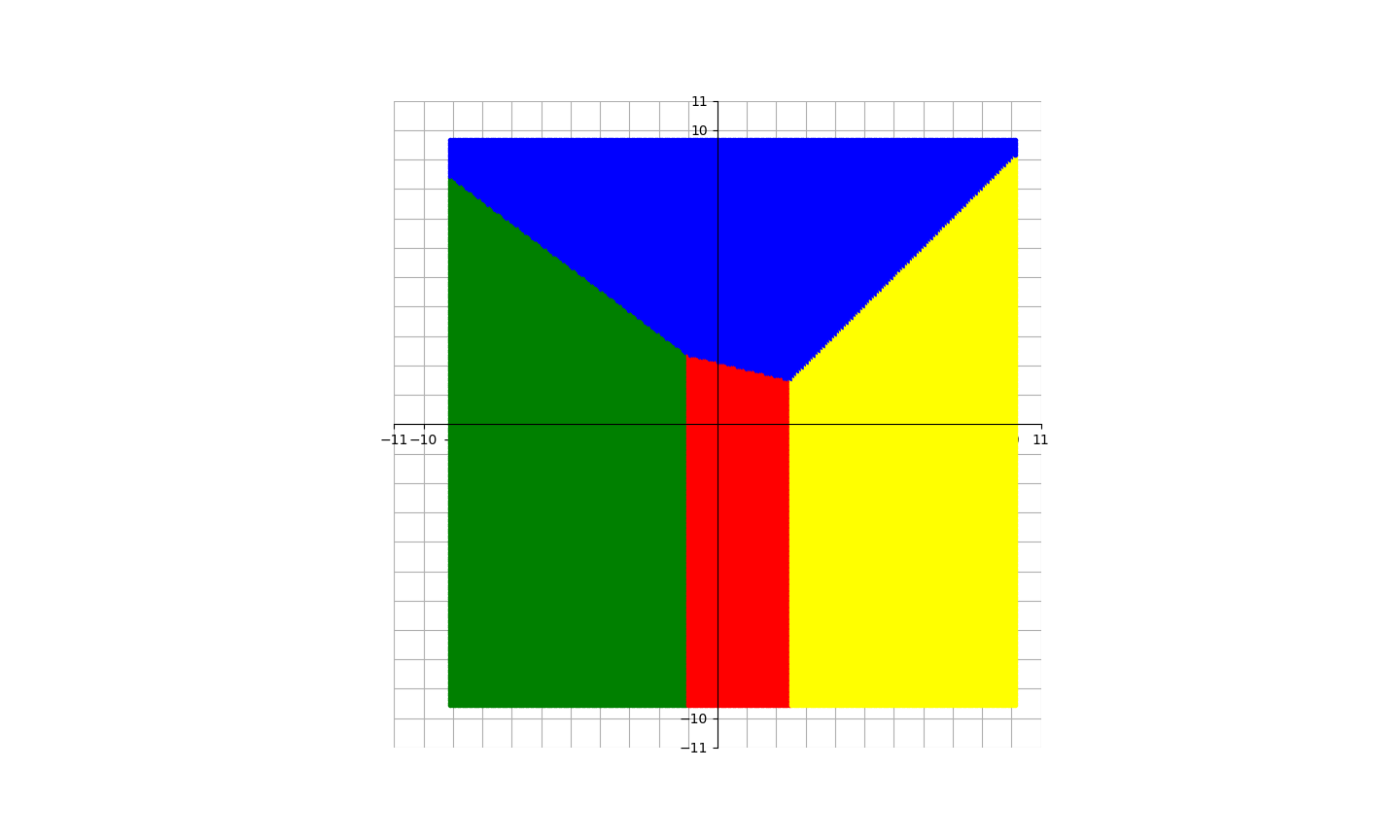}
    \includegraphics[width=10cm]{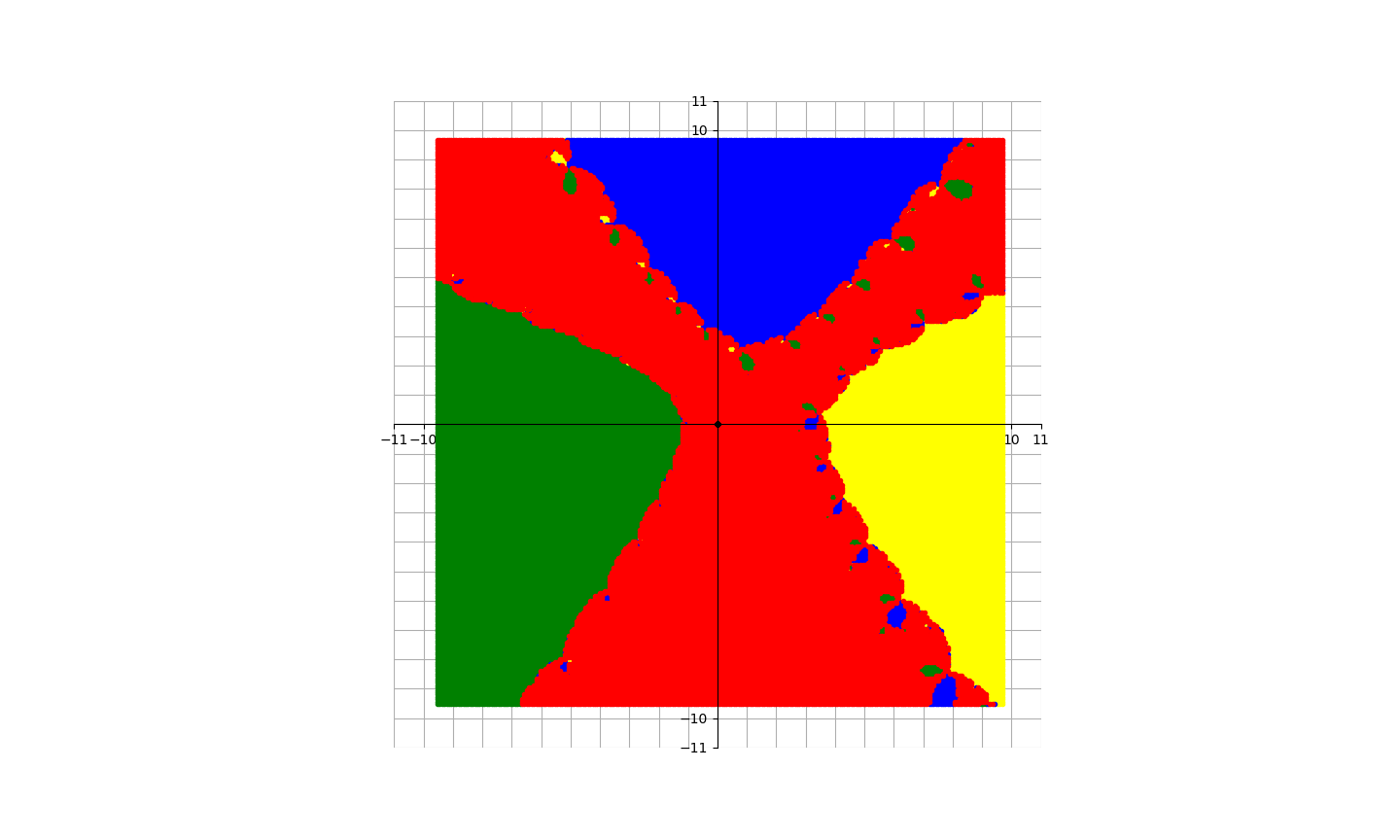}
    \includegraphics[width=10cm]{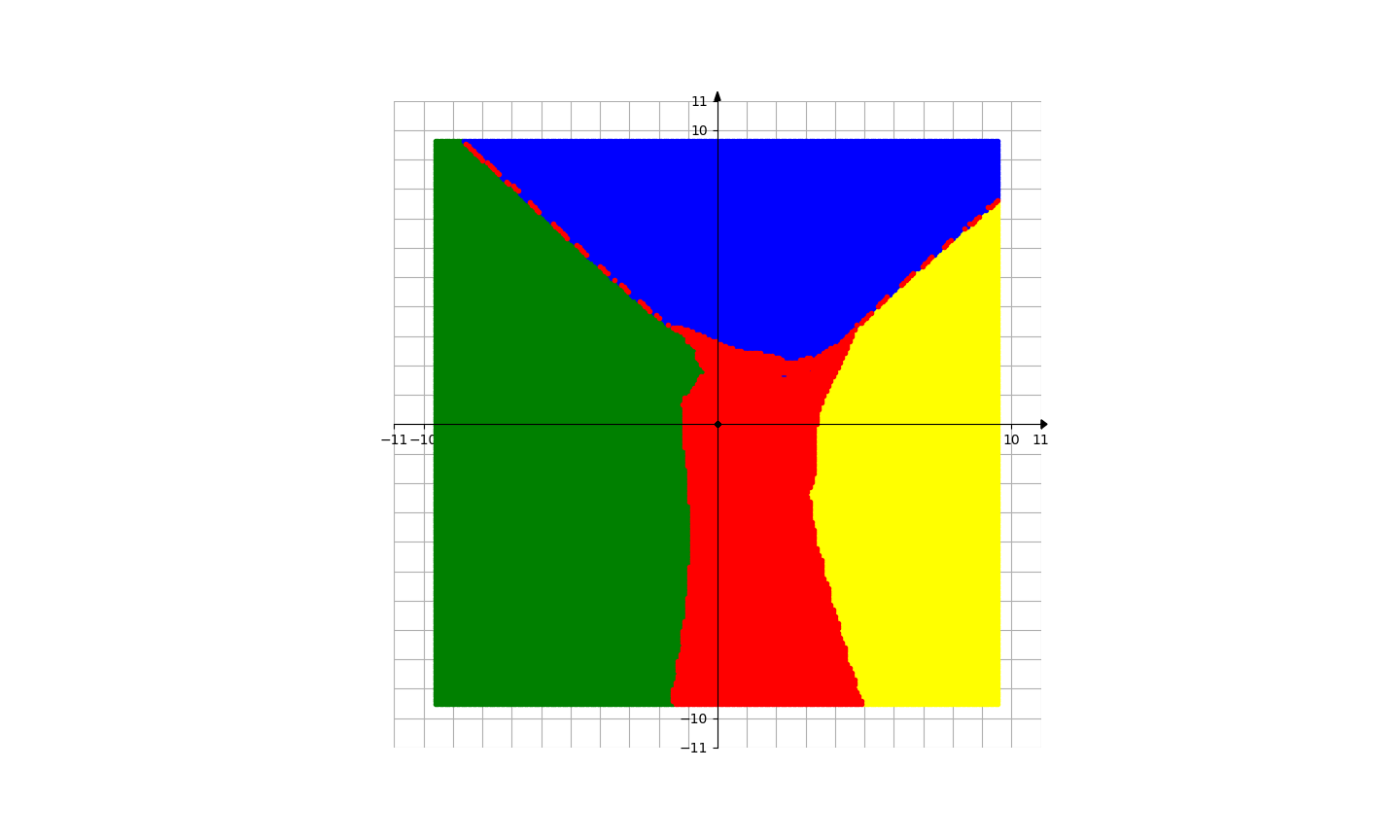}
    \caption{Basins of attraction for finding the roots of $f(z)=(z+2)(z-5)(z-1-4i)z$. Top picture: Voronoi's diagram of the 4 roots, center picture: Newton's method, bottom picture: BNQN. Points of the same colour  belong to basin of attraction of the same root.  Here the parameter $\theta$ in Algorithm \ref{table:alg0} is $1$.}
    \label{fig:FF3}
\end{figure}

Case 3c: $z_4^*$ is on a side of the triangle with 3 vertices $z_1^*$, $z_2^*$, $z_3^*$. Such an example is presented in Figure \ref{fig:FF3}. In this case, the picture for BNQN looks similar to Voronoi's diagram, except that there is a sequence of small open red domains (in the basin of attraction for the root $z_4*=0$) going to infinity.  

\begin{figure}
    \centering
    \includegraphics[width=10cm]{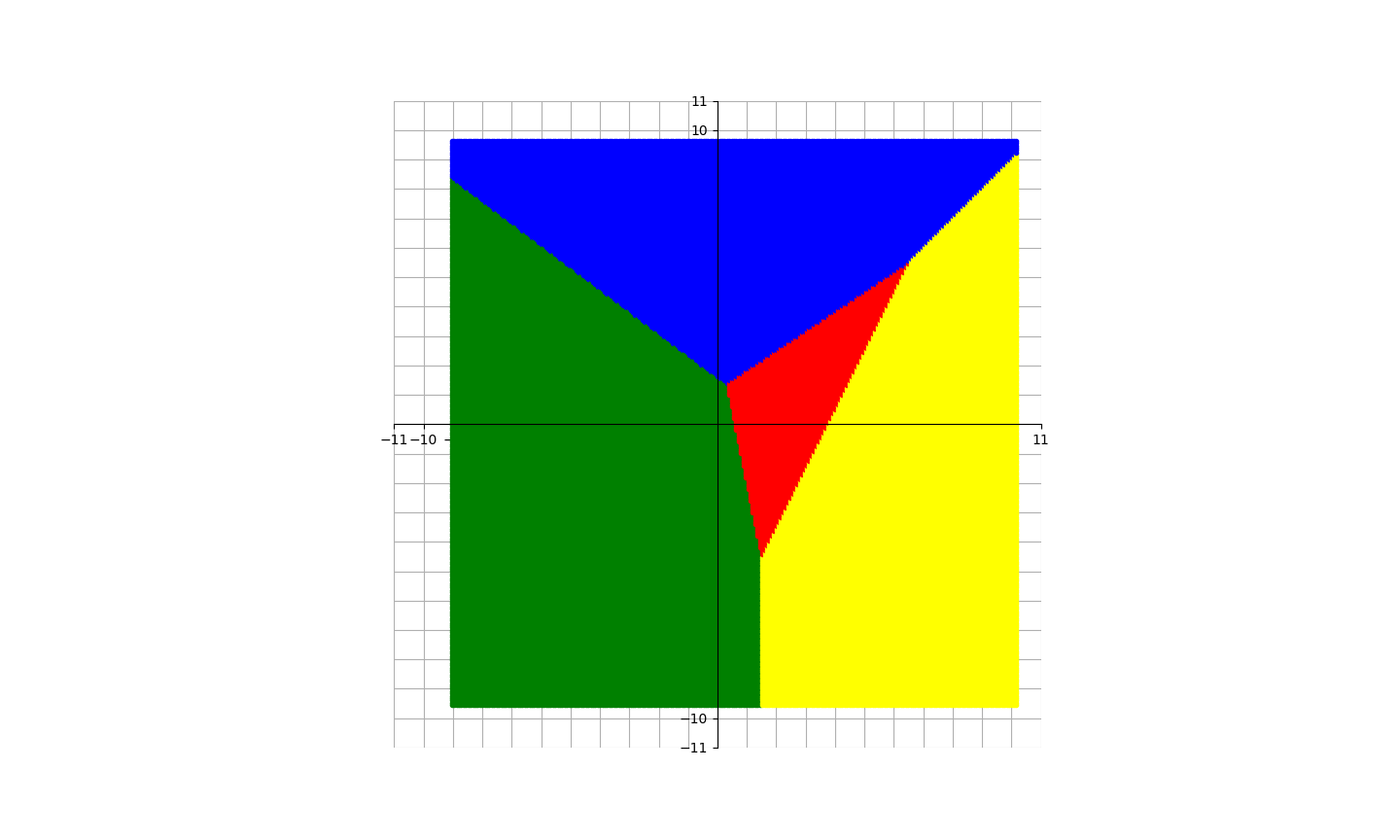}
    \includegraphics[width=10cm]{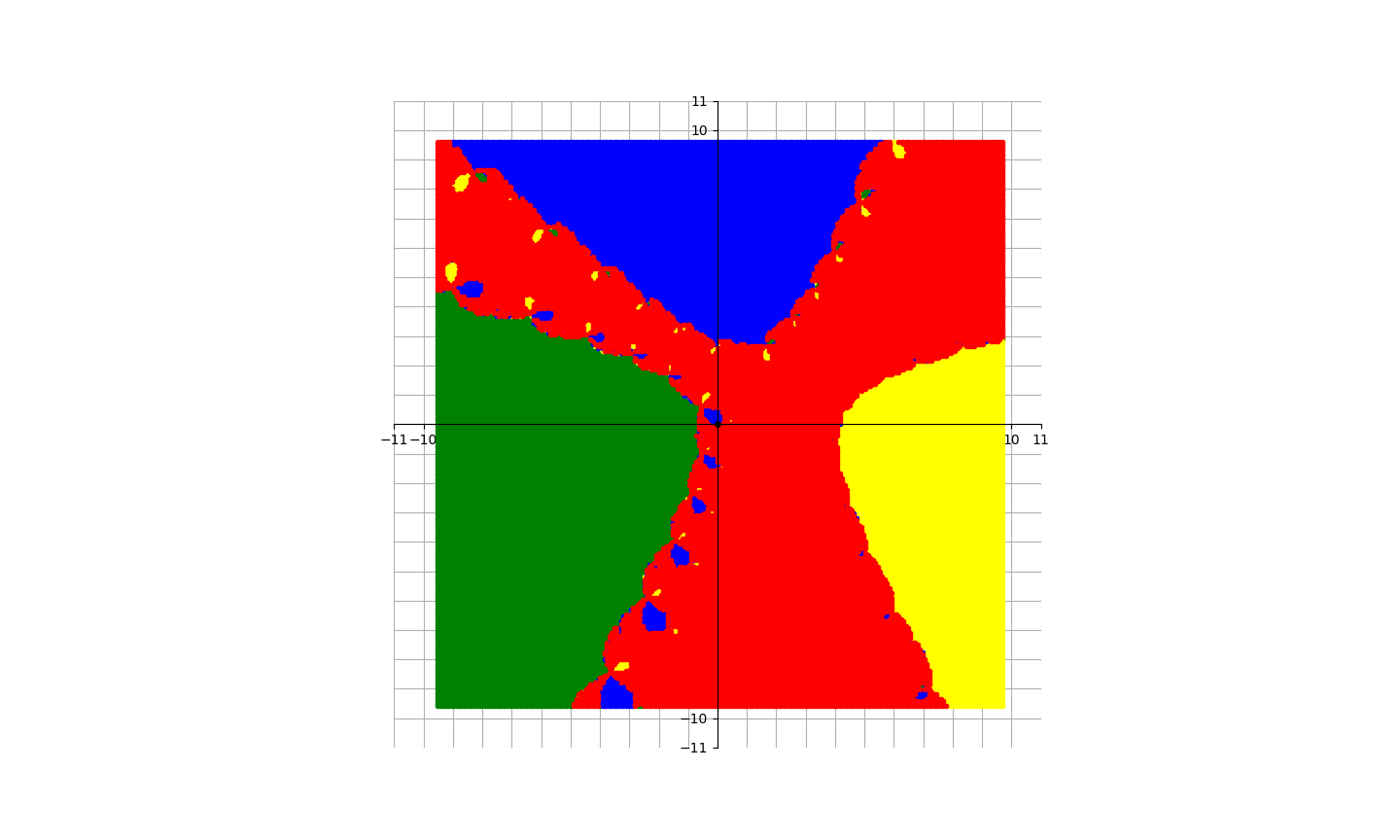}
    \includegraphics[width=10cm]{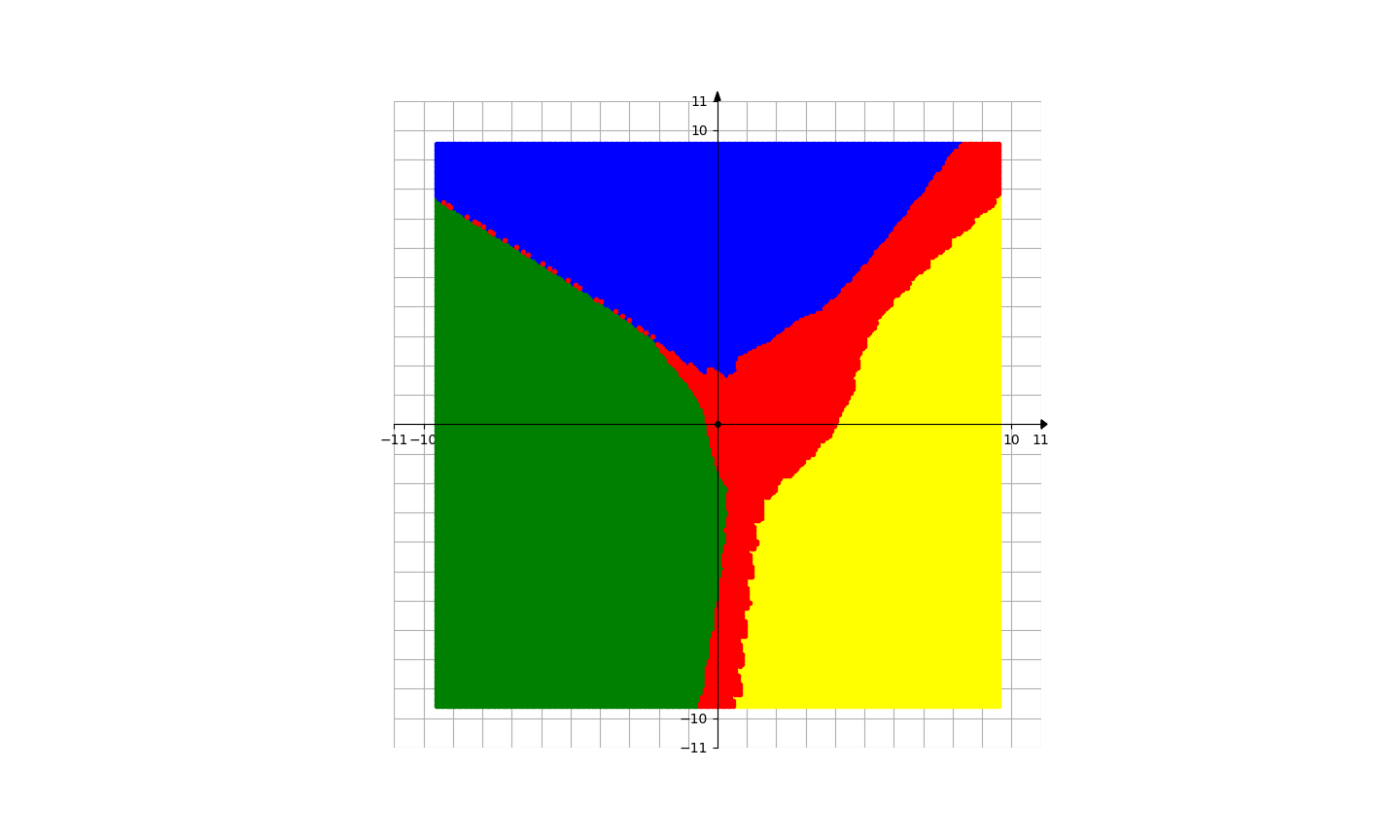}
    \caption{Basins of attraction for finding the roots of $f(z)=(z+2)(z-5)(z-1-4i)(z-3-i)$. Top picture: Voronoi's diagram of the 4 roots, center picture: Newton's method, bottom picture: BNQN. Points of the same colour  belong to basin of attraction of the same root. Here the parameter $\theta$ in Algorithm \ref{table:alg0} is $1$.} 
    \label{fig:FF4}
\end{figure}

Case 3d: $z_4^*$ is inside the triangle with 3 vertices $z_1^*$, $z_2^*$, $z_3^*$. Such an example is presented in Figure \ref{fig:FF4}. In this case the picture by BNQN is less similar to Voronoi's diagram: while Voronoi's cell for the root $z_4^*=3+i$ is completely bounded inside the Voronoi's cells of the other 3 roots, the basin of attraction of $z_4^*$ for BNQN opens up to $\infty$.  

{\bf Observation 4: There seems to be a ''channel to infinity" for every root.} 

In the experiments in the paragraph on Observation 3 above, as well as other experiments in \cite{RefT} and \cite{RefFHTW2}, it is apparent that in the case $f$ is a polynomial, then the immediate basin of attraction of every root opens up to $\infty$, or in other words, there is a ''passage to infinity" or a ''channel to infinity". A possible explanation is that the boundary curves do not close up, which leave open passages to infinity. We note that for Newton's method, such ''passages to infinity" of every roots of a polynomial are proven theoretically to exist, please see \cite{RefPrzytycki}, \cite{RefManning}, \cite{RefHSS} for more detail. 

\begin{figure}
    \centering
    \includegraphics[width=10cm]{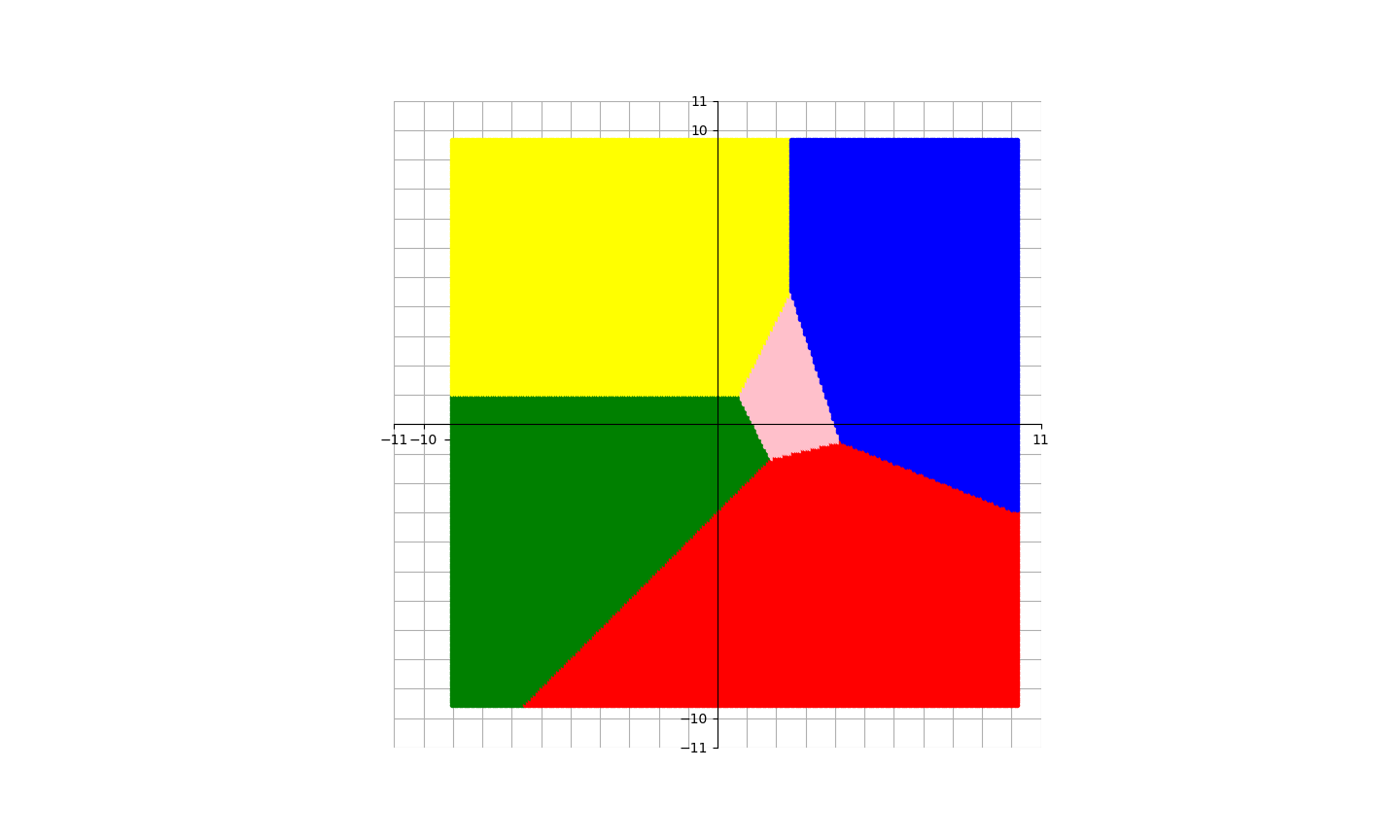}
    \includegraphics[width=10cm]{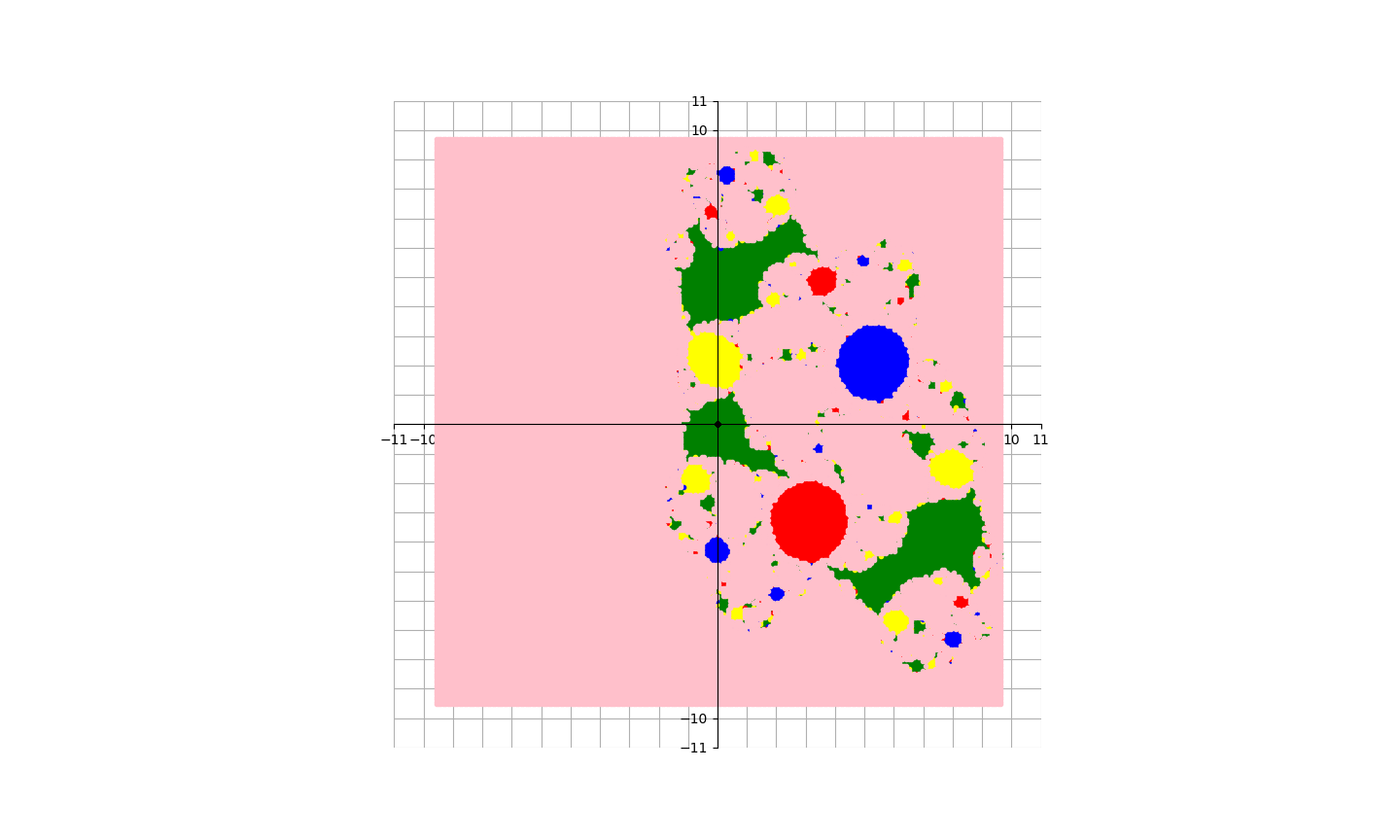}
    \includegraphics[width=10cm]{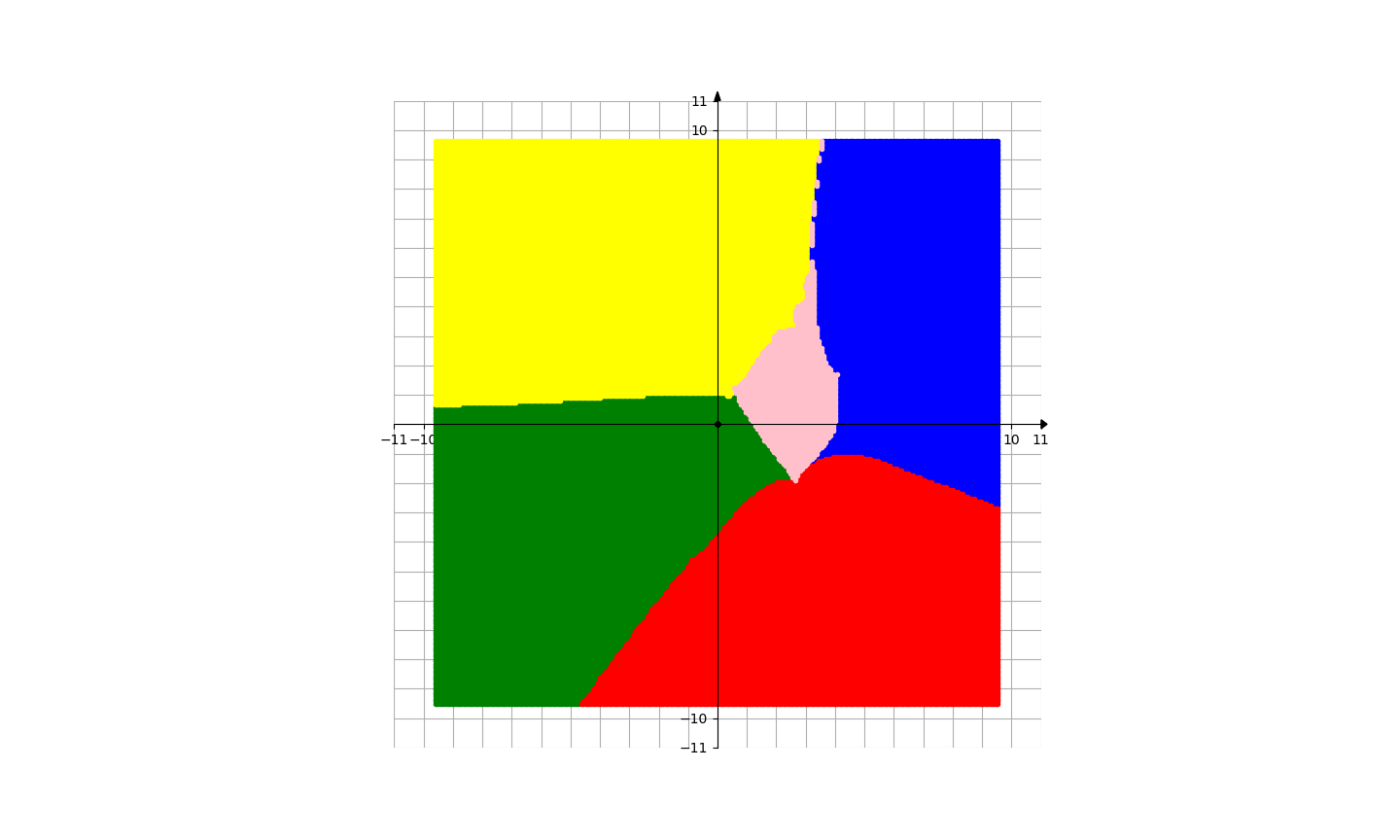}
    \caption{Basins of attraction for finding the roots of $f(z)=P(z)/P'(z)$, where $P(z)=z(z-2i)(z-5-2i)(z-3+3i)(z-2-i)$. Top picture: Voronoi's diagram of the 5 roots, center picture (from \cite{RefFHTW2}): Newton's method, bottom picture: BNQN. Points of the same colour  belong to basin of attraction of the same root. Here the parameter $\theta$ in Algorithm \ref{table:alg0} is $1$: note that the picture for BNQN here looks more smooth than that in \cite{RefFHTW2}, where $\theta =0$ was used.} 
    \label{fig:F17DividedByPrime}
\end{figure}

Experiments in \cite{RefFHTW2} display some interesting examples of functions with compact sublevels $f$ of the form $f=P(z)/P'(z)$,  where $P(z)$ is a polynomial.  See Figure \ref{fig:F17DividedByPrime} for the case $P(z)=z(z-2i)(z-5-2i)(z-3+3i)(z-2-i)$. Here, Newton's method only has ''channel to infinity" for one of the roots. On the other hand, BNQN seems to have ''channels to infinity" for 4 of the roots (it seems the basin of attraction of the remaining root contains a sequence of small open sets running off to infinity). Moreover, the geometric configurations in Voronoi's diagram and the basins of attraction for Newton's method seem reverse to each other: in the picture for Newton's method, the root $2+i$ (whose basin of attraction is in pink colour) is inside the convex hull of the other 4 roots, but its basin of attraction is the only one to have ''channel to infinity". On the other hand, the picture for BNQN still reflects well the geometric configuration of the roots and still has similarities to Voronoi's diagram: the basin of attraction of the root $2+i$, while not bounded, is still kept between the basins of attraction of other roots. 

\begin{thebibliography}{}

\bibitem{RefAlex}	A. S. Alexander,  A history of complex dynamics, Aspects of Mathematics, 1994.

\bibitem{RefAr} L. Armijo, {\it Minimization of functions having Lipschitz continuous first partial derivatives}, Pacific J. Math. 16, no. 1, 1--3, 1966. 

\bibitem{RefB} A. F. Beardon,  Iteration of rational functions, Springer-Verlag, New York, 1991.

\bibitem{RefBendixon} I. Bendixon, Sur les courbes d\'efinies par des \'equations diff\'erentielles, Acta Mathematica 1, 1--88, 1901.  

\bibitem{RefBer} W. Bergweiler,  Iteration of meromorphic functions, Bulletin of the American Mathematical Society, vol 29, number 2, October 1993, pp. 151--188,  1993.

\bibitem{RefCG} L. Carleson and T. W. Gamelin, Complex dynamics. Springer-Verlag, New York, 1993.  

\bibitem{RefD} 	 R. Devaney, An introduction to chaotic dynamical systems, CRC Press, 2023.

\bibitem{RefDDS}  T. C. Dinh, R. Dujardin and N. Sibony,  On the dynamics near infinity of some polynomial mappings in $\mathbf{C}^2$, Mathematische Annalen vol 333, pp.703--739, 2005.

%\bibitem{RefFatou} P. Fatour, Sur les \'equations fonctionnelles, Bulletin de la SMF 47, pp. 161--271, 1919.   

\bibitem{RefJongenEtAl}  H. Th. Jongen, P.Jonker and F. Twilt, The continuous Newton method for meromorphic functions, in Proceedings of the conference Geometric approaches to Differential equations, Scheveningen, The Netherlands 1979, Lecture notes in Mathematics, volume 810, edited by R. Martini, 1980.

\bibitem{RefFHTW2}  J. E. Forn\ae ss, M. Hu, T. T. Truong, T. Watanabe, 
Backtracking New Q-Newton's method, Newton's flow, Voronoi's diagram, and Stochastic root finding, Complex Analysis and Operator Theory, vol 18, article 112, 2024. 

\bibitem{RefFHTW}  J. E. Forn\ae ss, M. Hu, T. T. Truong, T. Watanabe, 
Backtracking New Q-Newton's method, Schr\"oder's theorem, and Linear conjugacy, arXiv:2312.12166. 

\bibitem{RefHK} F. von Haeseler, H, Kriete, The relaxed Newton's method for rational functions, Random Comput. Dynam., 3, 71--92, 1995.

\bibitem{RefHSS} J. Hubbard, D.Schleicher and S.Sutherland, How to find all roots of complex polynomials by Newton's method, Inventiones mathematicae, vol 146, 1--33, 2001.

%\bibitem{RefJulia} G. Julia, Sur quelques problemes relatifsa l'it\'eration des fractions rationnelles, CR Acad. Sci. Paris 168, pp. 147--149, 1919.  

\bibitem{kato} T. Kato, {\it Perturbation Theory for Linear Operators}, Originally publised as Vol 132 of the Grundlehren der mathematischen Wissenschaften, Springer-Verlag Berlin Heidelberg, 623 pages, 1995. 

%\bibitem{RefLeau} L. Leau, \'Etude sur les \'equations fonctionnelles \`a une ou \`a plusieurs variables, Annales de la Facult\'e des sciences de Toulouse: Math\'ematiques, vol 11.3.1897, E25--E110. 

\bibitem{RefManning} A. Manning, How to be sure of finding a root of a complex polynomial using Newton's method, Boletim da Sociedade Brasileira de Matematica, vol 22, no 2, 157--177, 1992. 

\bibitem{RefMc} C. McMullen, Families of rational maps and iterative root-finding algorithms, Ann. of Math. (2), 125, no 3, 467–493 , 1987.

\bibitem{RefMe} H.-G. Meier, The relaxed Newton-iteration for rational functions: the limiting case, Complex Variables Theory Appl., 16, 239--260, 1991.
     
\bibitem{RefM} J. Milnor, Dynamics in One Complex Variable. Princeton University Press, Princeton, 2006. 

\bibitem{RefMNTU} S.Morosawa, Y.Nishimura, M.Taniguchi and T.Ueda, Holomorphic Dynamics, Cambridge Studies in Advanced Mathematics, series number 66, revised version, 2000.

\bibitem{RefN} Wikipedia page on Newton's fractal, https://en.wikipedia.org/wiki/Newton$\_$fractal

%\bibitem{RefPei} H.-O. Peitgen, M. Prüfer, K. Schmitt,: 
%Global aspects of the continuous and discrete Newton method: a case study. 
%Acta Appl. Math. {\bf 13}, no.1-2, 123–202 (1988).

\bibitem{RefPoincare} H. Poincar\'e, Sur les courbes d\'efinies par un \'equation diff\'erentielle, Oeuvres, vol 1, Paris, 1892.   

\bibitem{RefPrzytycki} F. Przytycki, Remarks on the simple connectedness of basins of sinks for iterations of rational maps, in: Dynamical systems and Ergodic theory, ed. K. Krzyzewski, Polish scientific publishers, Warszawa, 229--235, 1989. 

\bibitem{RefSE2}  E. Schr\"{o}der, Ueber iterirte Functionen. Math. Ann., {\bf 3}, 296-322, 1871.

\bibitem{shub} M. Shub, {\it Global stability of dynamical systems}, Springer Science and Business Media, 1987.

\bibitem{RefS} H. Sumi, Negativity of Lyapunov exponents and convergence of generic random polynomial dynamical systems and random relaxed Newton’s method, Communications in Mathematical Physics, vol 384, pp.1513--1583, 2021.
      
\bibitem{RefThuanTuyen} T. Q. Tran and T. T. Truong, The Riemann hypothesis and dynamics of Backtracking New Q-Newton's method, arXiv:2405.05834.

\bibitem{RefT} T. T. Truong, Backtracking New Q-Newton's Method: A Good Algorithm for Optimization and Solving Systems of Equations (2023)	arXiv:2209.05378
		
\bibitem{RefTT} T. T. Truong, T. D.  To, H.-T. Nguyen, T.  H. Nguyen, H.  P. Nguyen, M. Helmy, A fast and simple modification of Newton’s method avoiding saddle points. J Optim Theory Appl, vol 199, pp.805--830, 2023.

\bibitem{RefV1} G. Voronoi, Nouvelles applications des parametres continus a la theorie des formes quadratiques, Premier memoire, Sur quelques proprietes des formes quadratiques positive parfaites, Journal fur die Reine und Angwandte Mathematik  (133), pp. 97--178, 1908. 

\bibitem{RefV2} G. Voronoi, Nouvelles applications des parametres continus a la theorie des formes quadratiques, Deuxieme memoire, Recherches sur les paralleloedres primitifs, Journal fur die Reine und Angwandte Mathematik  (134), pp. 198--287, 1908. 


\bibitem{RefVW}  Wikipedia page on Vornoi's diagrams and applications  https://en.wikipedia.org/wiki/Voronoi$\_$diagram$\#$Applications

\bibitem{RefY} Y. Yamagishi, On the local convergence of Newton's method to a multiple root, J. Math. Soc. Japan 55, 897--908, 2023.
        
\end{thebibliography}
\end{document}